\documentclass[a4paper,english]{article}
\usepackage[T1]{fontenc}
\usepackage[latin9]{inputenc}
\usepackage{color}
\usepackage{babel}
\usepackage{array}
\usepackage{textcomp}
\usepackage{amsthm}
\usepackage{amsmath}
\usepackage{amssymb}
\usepackage{graphicx}
\usepackage{subscript}
\usepackage[unicode=true,
 bookmarks=true,bookmarksnumbered=true,bookmarksopen=true,bookmarksopenlevel=4,
 breaklinks=false,pdfborder={0 0 1},backref=false,colorlinks=true,pdfpagemode=FullScreen]
 {hyperref}
\hypersetup{pdftitle={Latin Polytopes},
 pdfauthor={Miguel G. Palomo},
 pdfsubject={Latin squares, Latin triangles, Combinatorics},
 pdfkeywords={Latin squares, Latin triangles, Combinatorics},
 linkcolor=blue,citecolor=blue,filecolor=blue,urlcolor=blue}

\makeatletter

\pdfpageheight\paperheight
\pdfpagewidth\paperwidth

\newcommand{\lyxmathsym}[1]{\ifmmode\begingroup\def\b@ld{bold}
  \text{\ifx\math@version\b@ld\bfseries\fi#1}\endgroup\else#1\fi}

\providecommand{\tabularnewline}{\\}

\numberwithin{equation}{section}
\numberwithin{table}{section}
\theoremstyle{plain}
\newtheorem{thm}{\protect\theoremname}[section]
\theoremstyle{definition}
\newtheorem{example}[thm]{\protect\examplename}
\theoremstyle{definition}
\newtheorem{defn}[thm]{\protect\definitionname}
\ifx\proof\undefined
\newenvironment{proof}[1][\protect\proofname]{\par
\normalfont\topsep6\p@\@plus6\p@\relax
\trivlist
\itemindent\parindent
\item[\hskip\labelsep\scshape #1]\ignorespaces
}{%
\endtrivlist\@endpefalse
}
\providecommand{\proofname}{Proof}
\fi
\theoremstyle{plain}
\newtheorem{lem}[thm]{\protect\lemmaname}

\date{}

\renewcommand\paragraph{\@startsection{paragraph}{4}{\z@}%
   {-3.25ex\@plus -1ex \@minus -.2ex}%
   {.75ex \@plus .2ex}%
   {\normalfont\normalsize\bfseries}}
\usepackage{caption}
\renewcommand\subparagraph{\@startsection{subparagraph}{4}{\z@}%
                                     {-3.25ex\@plus -1ex \@minus -.2ex}%
                                     {.75ex \@plus .2ex}%
                                     {\normalfont\normalsize\bfseries}}

\usepackage{hyperref}

\usepackage[margin=10pt,font=small,labelfont=bf,labelsep=period]{caption}

\AtBeginDocument{

}

\makeatother

\providecommand{\definitionname}{Definition}
\providecommand{\examplename}{Example}
\providecommand{\lemmaname}{Lemma}
\providecommand{\theoremname}{Theorem}

\begin{document}

\title{Latin Polytopes}

\author{Miguel G. Palomo}
\maketitle
\begin{abstract}
\noindent \addcontentsline{toc}{section}{Abstract}Latin squares are
well studied combinatorial objects. In this paper we generalize the
concept and propose new objects like \emph{Latin triangles}, \emph{free
Latin squares}, \emph{Latin tetrahedra}, \emph{free Latin cubes},
etc. We start with a classic definition of Latin squares followed
by one based on the concept of \emph{latinized board}. A Latin square
appears then as a combinatorial design whose points are geometric
and whose lines (the rows and columns) are invariant under the symmetries
of the square. The generalization that follows proceeds by 1. broadening
this geometric symmetry 2. considering more general configurations
of points and 3. admitting lines that intersect more freely. The resulting
concept is the \emph{Latin board}. Finally, we particularize Latin
boards to define \emph{Latin polytopes}, \emph{Latin polygons} and
\emph{Latin polyhedra}.

\textbf{\footnotesize{}Keywords: }{\footnotesize{}asterism, biregular
polygon, board, Board Symmetrization Problem, class-symmetric board,
combinatorial design, hypergraph coloring, kaleidoscopic source, labeling,
latinization, Latin board, Latin polygon, Latin polyhedron, Latin
polytope, Latin puzzle, Latin square, source of symmetry, symmetric
board, symmetric Latin board, warp class, weft class, woven board.} 
\end{abstract}

\section{Introduction}

\label{def:A-Latin-square}
\begin{figure}
\begin{centering}
\includegraphics[height=2.5cm]{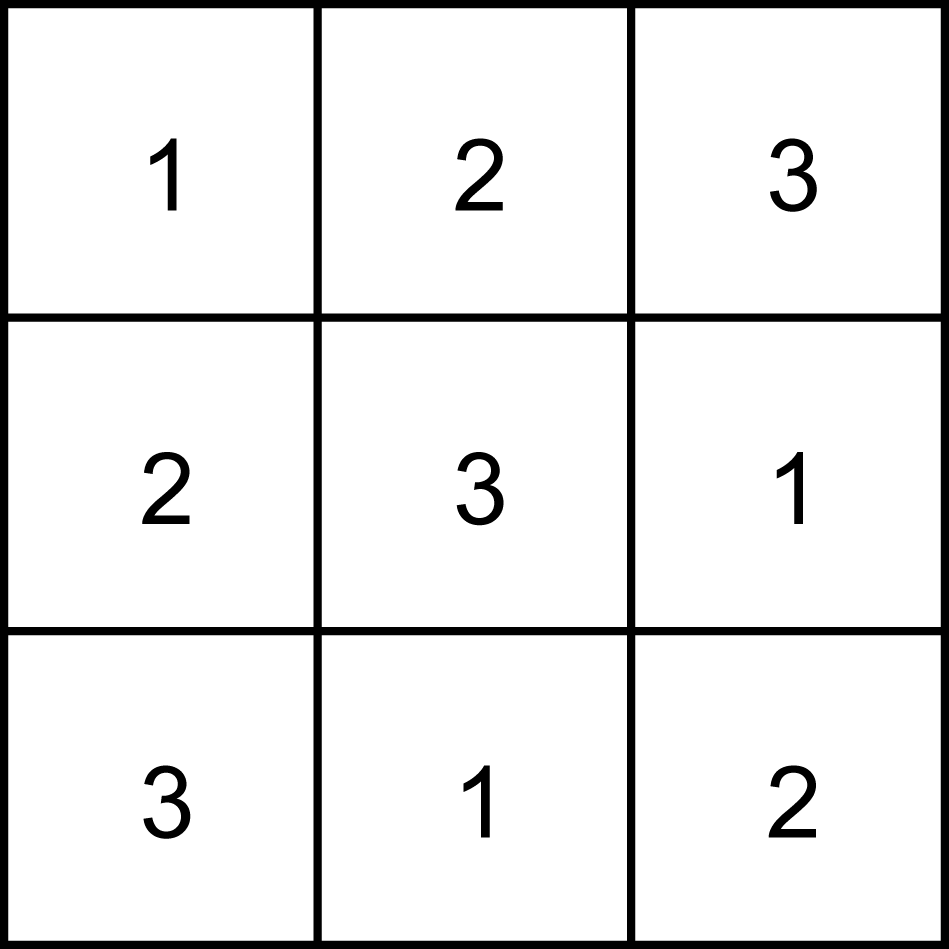}
\par\end{centering}

\protect\caption{\label{fig:latinSquareOrder3}Latin square}
\end{figure}

A Latin square of order $n$ is an $n\times n$
array filled with $n$ different labels, each occurring exactly once
in each row and exactly once in each column \cite{Latin squares and their applications}.
An example of order three with labels $1,2$ and $3$ is shown in
Fig. \ref{fig:latinSquareOrder3}\footnote{The set of labels is usually $\left\{ 1,2\text{,\dots},n\right\} $
because they double as indices for the rows and columns in the square.}.

Latin squares have a long and rich history. They are mentioned in
relation to magic squares as early as the $17$\textsuperscript{th}
century, although they seem to have been used much earlier in amulets
\cite{Handbook of Combinatorial Designs}. They are called ``Latin''
because the $18$\textsuperscript{th} century Swiss mathematician
Leonhard Euler used latin letters in his paper \emph{On Magic Squares}
\cite{De quadratis magicis}. 

\noindent In the $19$\textsuperscript{th} century the connection
between Latin squares and finite projective geometries was studied,
as well as their link with Algebra\footnote{Latin squares are multiplication tables of quasi-groups \cite{Latin squares and their applications}.}
and Combinatorics. Already in the $20$\textsuperscript{th} century,
applications to experimental design were formalized, especially in
agricultural engineering \cite{Handbook of Combinatorial Designs}. 

More recently, the advent of computers and digital networks has brought
codes based on Latin squares for error detection and correction in
data telecommunication\footnote{For additional historical facts about Latin squares see \cite{The history of Latin squares Andersen,An illustrated philatelic introduction,Some comments on Latin squares and on Graeco-Latin squares}.}.

\subsubsection*{Geometric Symmetry as a Source of New Objects}

\noindent Looking at the Latin square in Fig. \ref{fig:latinSquareOrder3},
it is easy to see that the result of a reflection across its vertical
symmetry axis results in another Latin square. Other geometric transformations
share this property: for example a $90\text{º}$ counter-clockwise
rotation around its center, or a reflection across any of the two
diagonals. 

\begin{figure}
\begin{centering}
\includegraphics[height=5cm]{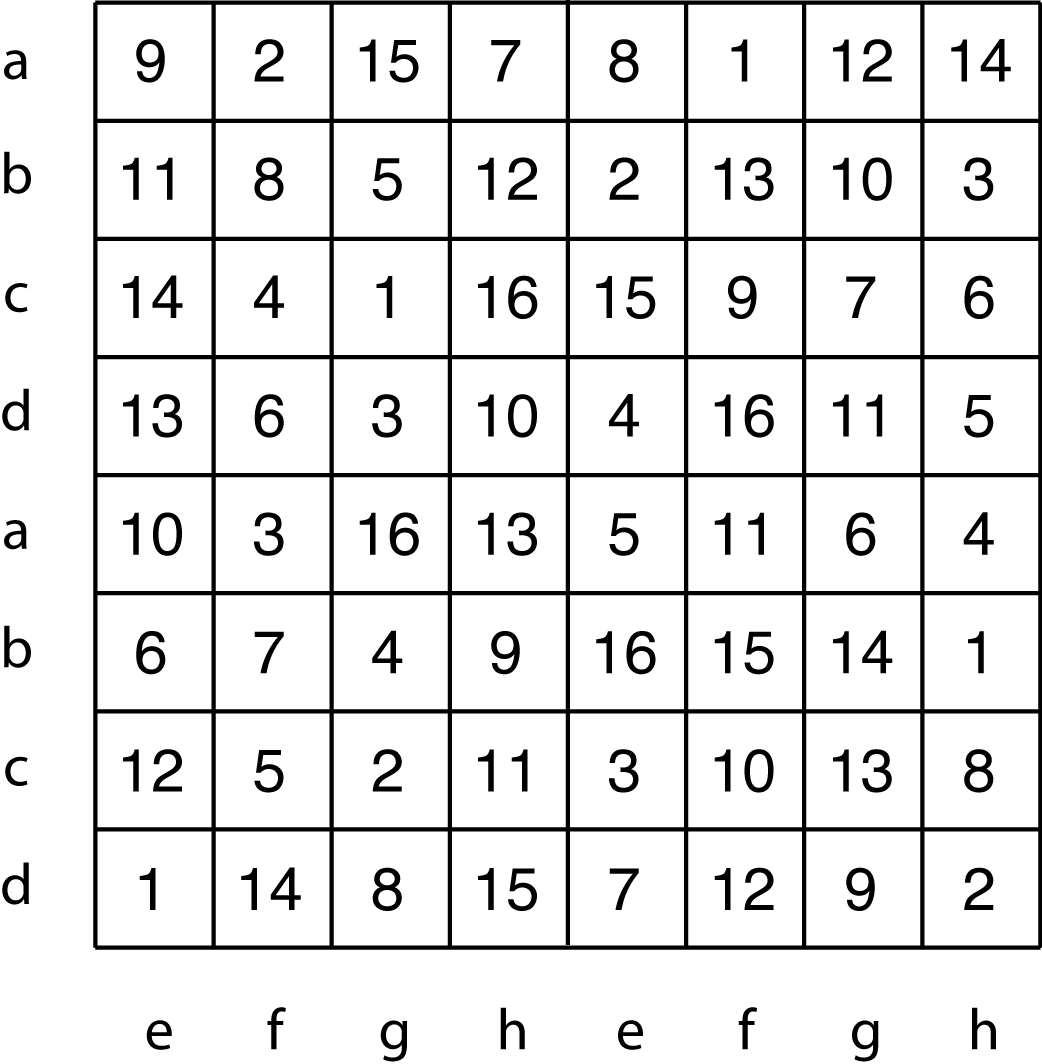}
\par\end{centering}

\protect\caption{\label{fig:free Latin square example}Free Latin square}
\end{figure}

In fact any symmetry transformation of the square produces another
Latin square. Geometric symmetry alone produces then new Latin squares
from a given one. The examples that follow show other objects for
which symmetry is also the source of new similar objects.
\begin{example}
Fig. \ref{fig:free Latin square example} shows a square mesh of cells
in which each pair of rows with the same letter holds all numbers
from $1$ to $16$, something that also applies to pairs of columns
similarly labeled. 

This property still holds if, keeping the letters in place, we apply
as before any symmetry of the square, like a reflection across the
horizontal symmetry axis or a $180\text{º}$ counter-clockwise rotation
around its center. In this paper we call objects like this one \emph{free
Latin squares.}
\end{example}

\begin{figure}
\begin{centering}
\includegraphics[height=5cm]{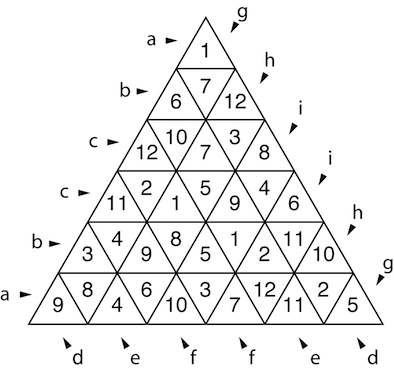}
\par\end{centering}

\protect\caption{\label{fig:Latin triangle example}Latin triangle}
\end{figure}

\begin{example}
Fig. \ref{fig:Latin triangle example} shows a triangular mesh of
cells in which each pair of stripes of triangles pointed to by the
same letter has all numbers from $1$ to $12$. 

This condition still holds if, keeping the letters in place, we apply
any symmetry transformation of the equilateral triangle, like a $120\text{º}$
counter-clockwise rotation around its center. This is an example of
\emph{Latin triangle}.
\end{example}

\begin{figure}
\begin{centering}
\includegraphics[height=5cm]{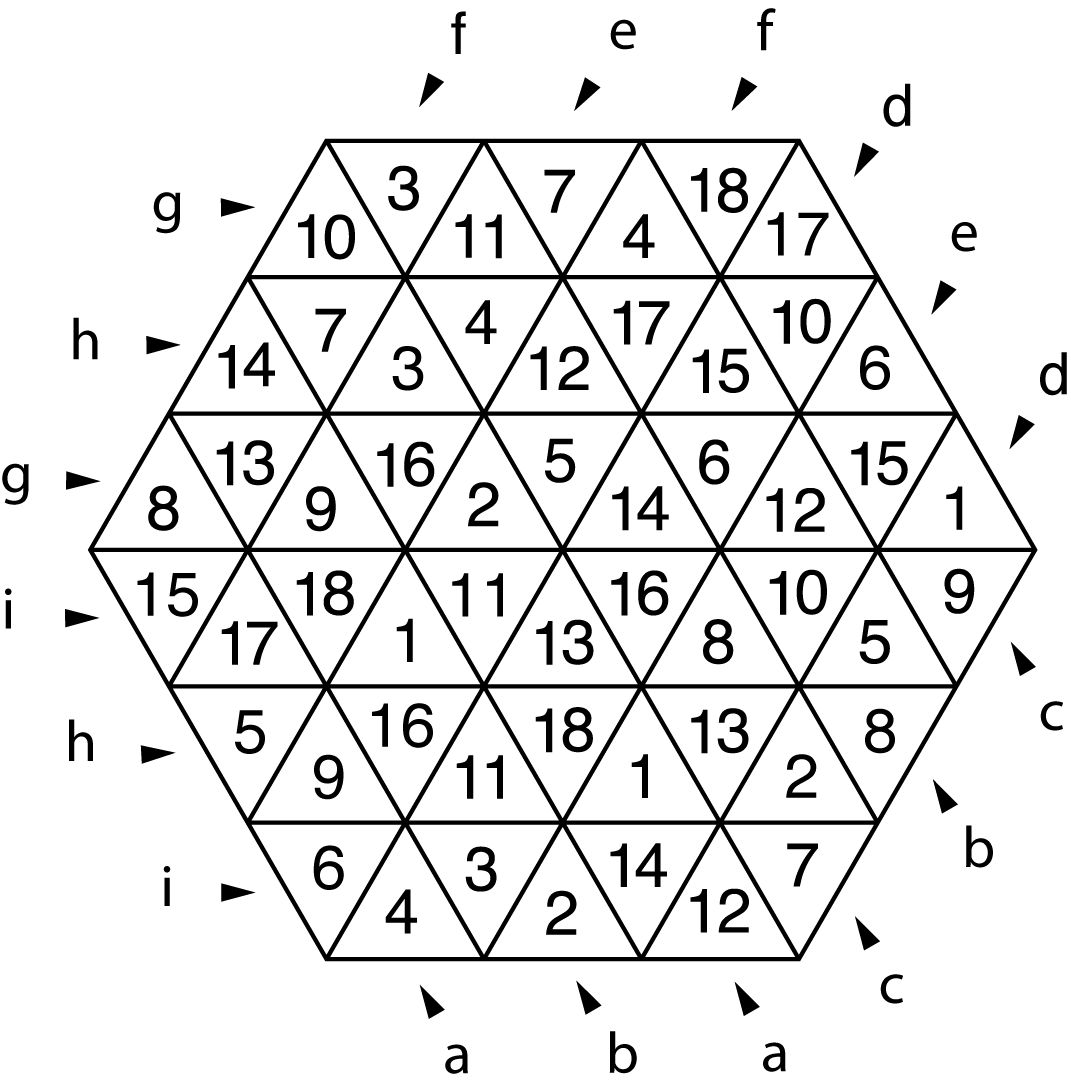}
\par\end{centering}

\protect\caption{\label{fig:Latin hexagon example}Latin hexagon}

\end{figure}

\begin{example}
Fig. \ref{fig:Latin hexagon example} shows a triangular mesh in which
each pair of stripes of triangles pointed to by the same letter contains
all numbers from $1$ to $18$. 

Here, the property also holds after any symmetry transformation of
the regular hexagon, for example a $60\text{º}$ counter-clockwise
rotation around its center. We call objects like this one \emph{Latin
hexagons}.
\end{example}
\noindent In this paper we generalize these objects by first defining
\emph{boards} as a pair of sets: one containing geometric points,
the other subsets of the points. Next, we define\emph{ Latin boards
}as specific labelings of the subsets. In both cases geometric symmetry
is not taken into account\emph{.} 

Next we define \emph{symmetric boards} as boards that are invariant
under the action of a symmetry group. We finally introduce \emph{Latin
polytopes}, a particularization of Latin boards\emph{ }that, just
like the examples mentioned, produce new ones under symmetry transformations. 

The \hyperref[TOC]{Table of Contents} at the end of the document
contains a detailed view of this plan.

\pagebreak{}

\section{Boards}

\noindent In this section and the next, a generalization of Latin
squares that does not involve geometric symmetry is proposed.

\subsection{Definition of Board}
\begin{defn}
\label{def:board}A \textit{board }is a pair $(P,C)$ where \textit{$P$}
is a finite set of \textit{\emph{geometric points}} and \textit{$C$}
is a set of subsets of $P$ whose union is\textit{\emph{ }}\textit{$P$.}
\end{defn}
\noindent It follows from the definition that intersection among subsets
is allowed. The set $C$ is the \emph{constellation} of the board.
Each of its subsets is an \emph{asterism}\footnote{We have borrowed these terms from Astronomy. They seem less ambiguous
in this context than others used in related objects (see Sect. \ref{sub:Incidence Structures-Related-to-boards}).}.
\begin{example}
Fig. \ref{fig:fanoPlaneUnsymmetric} (left) shows a set of seven $\mathbb{R}^{2}$
points. On the right we have formed seven asterisms, each one with
three points connected by two segments of the same color. The resulting
board is a geometric representation of the\emph{ Fano Plane} \cite{The Fano Plane},
a combinatorial structure with nice properties of balance: it has
seven points and seven asterisms, with three points in each asterism
and three asterisms sharing each point; furthermore, every two points
define an asterism and every two asterisms share a point.
\begin{figure}
\begin{centering}
\includegraphics[height=2.5cm]{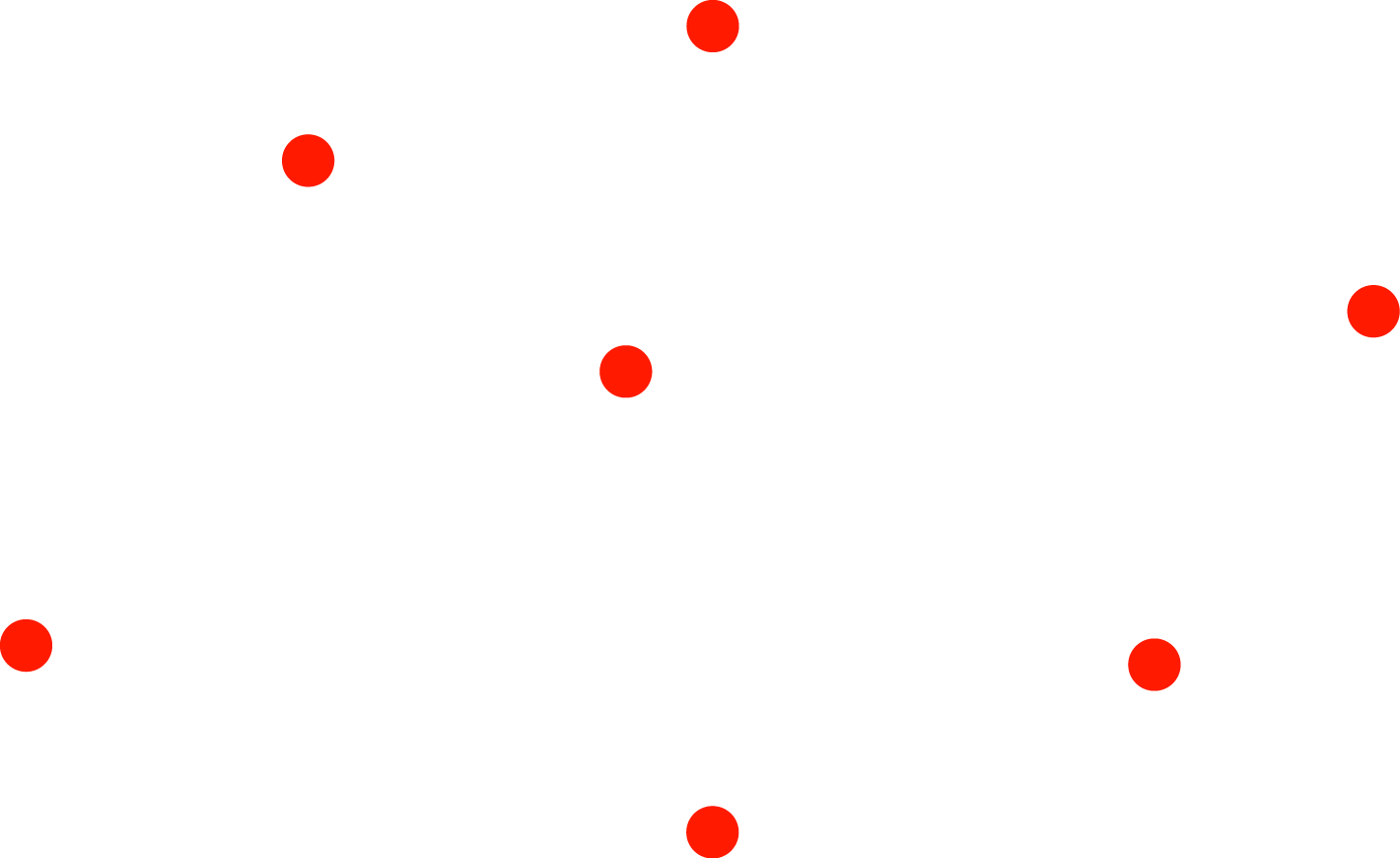}~~~~~~~\includegraphics[height=2.5cm]{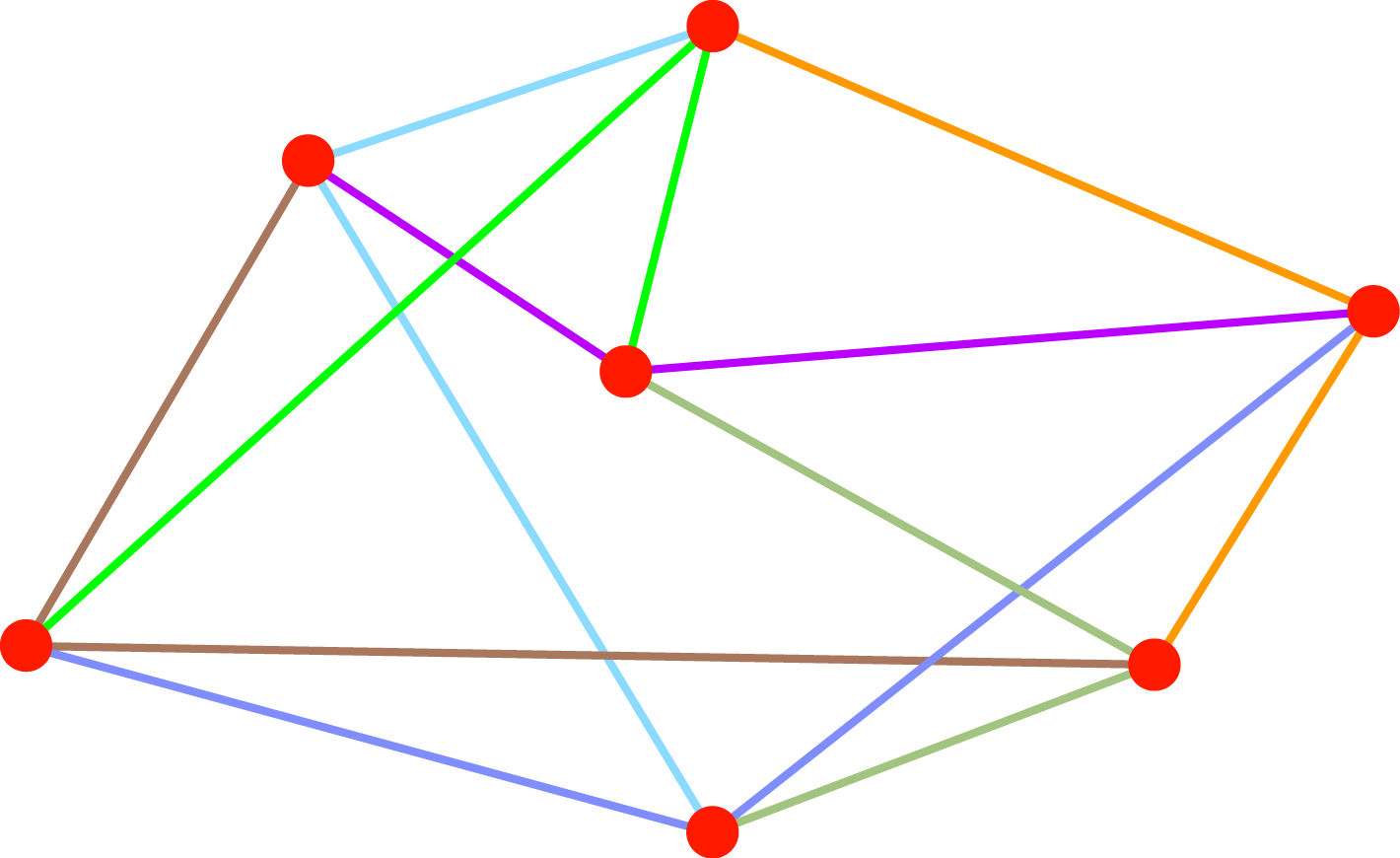}
\par\end{centering}

\protect\caption{\label{fig:fanoPlaneUnsymmetric}Points and asterisms in the Fano
Plane}
\bigskip{}
\end{figure}

\end{example}

\begin{example}
\noindent Let \textit{$P$} be the set of nine euclidean points in
Fig. \ref{fig:3x3points and vert and hor asterisms}. Let\textquoteright s
make an asterism out of each set of three points connected by segments
of the same color (right), and let \textit{$C$} be the constellation
formed by all asterisms. As every point is in \emph{$C$} via an asterism
we have the board $B=\left(P,C\right)$\label{board 3x3}.
\begin{figure}
\begin{centering}
\includegraphics[height=2cm]{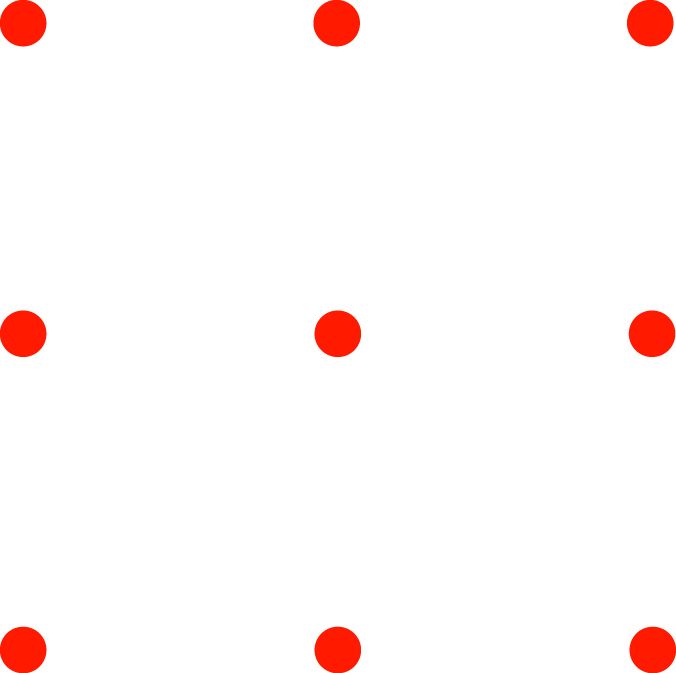}~~~~~~~~~~\includegraphics[height=2cm]{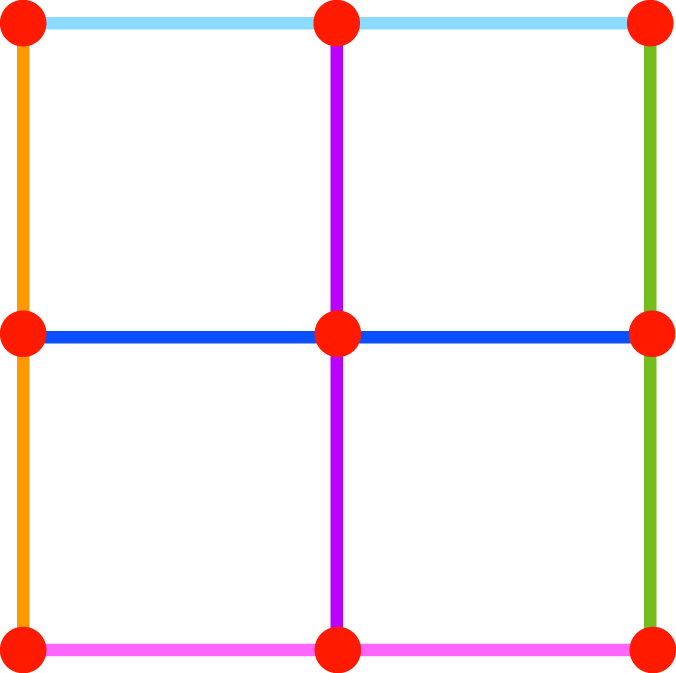}
\par\end{centering}

\protect\caption{\label{fig:3x3points and vert and hor asterisms}Board with nine points
and six asterisms}
\end{figure}

\end{example}

\subsubsection{Points and Cells in Boards\label{par:Points-and-Cells-in-boards}}

\noindent We will often represent a board with a grid of cells. To
conciliate this with Def. \ref{def:board} we will adopt the following
convention: ``cell'' will mean either \emph{centroid of the cell}
or \emph{point,} depending on the context. For example ``write a
label in the cell'' will mean \emph{label the centroid of the cell},
whereas ``a board with $n$ cells'' will be \emph{a board with $n$
points}.
\begin{example}
\noindent If we remove the numbers inside the cells in Fig. \ref{fig:latinSquareOrder3},
the resulting object is a board whose points are the (centroids of
the) cells\footnote{In this specific case, the centers of the small squares.}.
Each row and column of (centroids of the) cells is an asterism of
this board.
\end{example}

\subsection{Incidence Structures Related to Boards\label{sub:Incidence Structures-Related-to-boards}}

An abstract system consisting of two types of objects and a single
relation between these types of objects is called an \emph{incidence
structure} \cite{Handbook of Combinatorial Designs}. Following Def.
\ref{def:board}, boards are incidence structures in which the two
types of objects are the board points and the asterisms. The incidence
relation here is the membership of points in asterisms. We see next
how boards relate to other existing incidence structures.

\subsubsection{Families of Sets\label{sub:Families-of-Sets}}
\begin{defn}
A collection $C$ of subsets of a given set $P$ is called a \emph{family
of subsets} of $S$, or a \emph{family of sets} over $S$. 
\end{defn}
\noindent Identifying respectively the points and the constellation
of a board with the elements and the collection of subsets, a board
is a family of sets in which no set is repeated in the family.

\subsubsection{Finite Geometries\label{sub:Finite geometries}}

Finite geometries are geometric systems with a finite number of \emph{points}
and a collection of subsets of the points called \emph{lines}. The
incidence relation here is geometric incidence. The most studied finite
geometries are \emph{finite projective} \emph{geometries} (in which,
among other conditions, any two lines have a point in common) and
\emph{finite affine} \emph{geometries} (in which a necessary condition
is: if a point $p$ is not on a line $l_{1}$ then there is a unique
line $l_{2}$ through $p$ having no points in common with $l_{1}$).
As a finite geometry, the Fano Plane in Fig. \ref{fig:fanoPlaneUnsymmetric}
(right) is the smallest finite projective plane.

Identifying respectively the points and lines of the geometry with
the points and asterisms of a board we see that certain boards are
also finite geometries.

\subsubsection{Simple Combinatorial Designs\label{sub:Simple-Combinatorial-Designs}}
\begin{defn}
\label{def: simple-combinatorial-design}A \emph{simple combinatorial
design} is a pair $\left(P,L\right)$ where $P$ is a set of \emph{points}
and $L$ is a set of subsets of $P$ called \emph{lines} \cite{Stinson Combinatorial Designs}.
\end{defn}
\noindent Boards are then simple combinatorial designs\emph{ }whose
points are geometric and contained inside lines.

\subsubsection{Hypergraphs\label{sub:Geometric-Hypergraphs}}
\begin{defn}
A \emph{hypergraph} is a pair $\left(V,E\right)$ where $V$ is a
finite set of elements called \emph{vertices} and \emph{E} is a family
of non-empty subsets of \emph{V} called \emph{hyperedges} whose union
is \emph{V.}
\end{defn}
\noindent Hypergraphs become \emph{graphs} when hyperedges have just
one or two vertices \cite{hypergraphs}. When the vertices of a hypergraph
are geometric points we have a \emph{geometric hypergraph}. If we
call the points and asterisms of a board \emph{vertices} and \emph{hyperedges}
respectively, boards are also geometric hypergraphs.
\begin{example}
\noindent Let \emph{V} be the nine euclidean points in Fig. \ref{fig:3x3points and vert and hor asterisms}
(left). Let\textquoteright s make a hyperedge out of the three points
connected by segments of the same color (right). Then the pair $\left(V,E\right)$
is a geometric hypergraph.
\end{example}
\noindent Hypergraphs, finite geometries, combinatorial designs and
families of sets are then similar structures, the differences lie
in the focus of interest. In hypergraphs the focus is on graph theoretic
questions, like connectivity and \emph{colorability} \cite{hypergraphs}.
In finite geometries it is on geometric aspects like parallelism,
collinearity and incidence. The incidence relation here is also much
stricter than in other incidence structures. In combinatorial design
the interest lies on design theoretic issues, like combinatorial symmetry
and balance. In families of sets the focus is on set theoretic questions,
like \emph{Sperner} \emph{theory} for example (see Sect. \ref{sub:Asymmetric-Boards-with-unknown automorphism group}).

\subsection{Board Automorphisms\label{sub:Board-Automorphisms}}
\begin{defn}
\label{def:isomorphic boards}Two boards \textit{$B=(P,C)$} and \textit{$B\text{\textquoteright}=(P\text{\textquoteright},C\text{\textquoteright})$}
are \emph{isomorphic} if and only if there is a bijection from $P$
to $P\lyxmathsym{\textquoteright}$ that preserves the constellations.
Such bijection is called an \emph{isomorphism}.
\end{defn}

\begin{defn}
An \textit{automorphism} of a board \textit{$B=(P,C)$} is an isomorphism
from $B$ to itself.
\end{defn}
\noindent As the next example shows, isomorphisms are specific permutations
of the set of points.
\begin{example}
\noindent If we number the points in the board in Fig. \ref{fig:3x3points and vert and hor asterisms}
as indicated in Fig. \ref{fig:numberedPointsB1}, we have\renewcommand{\arraystretch}{1.5}
\textit{
\[
\begin{array}{c}
P=\left\{ 1,2,3,4,5,6,7,8,9\right\} \\
C=\left\{ \left\{ 1,4,7\right\} ,\left\{ 2,5,8\right\} ,\left\{ 3,6,9\right\} ,\left\{ 1,2,3\right\} ,\left\{ 4,5,6\right\} ,\left\{ 7,8,9\right\} \right\} 
\end{array}
\]
}\renewcommand{\arraystretch}{1}
\begin{figure}
\begin{centering}
\includegraphics[height=2.5cm]{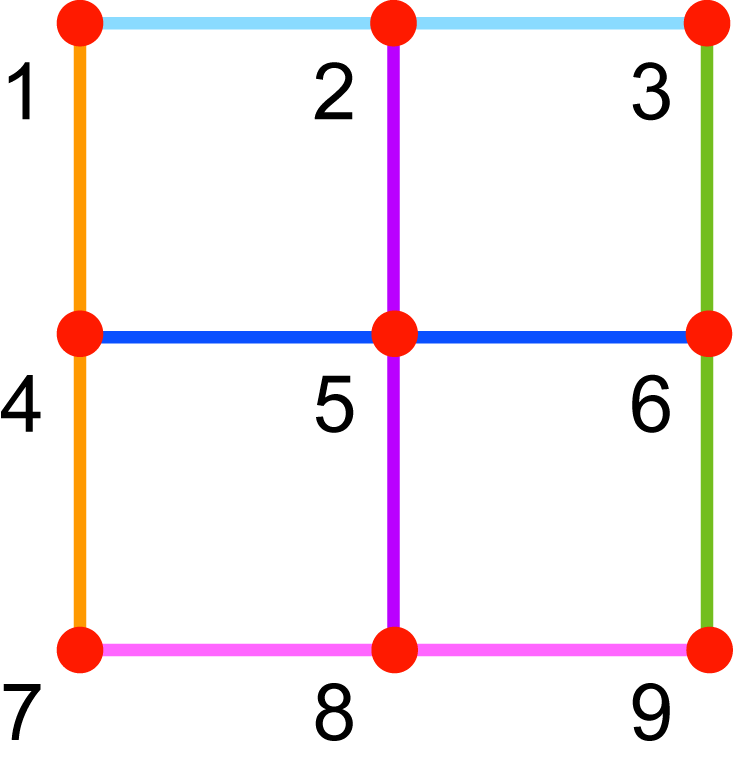}
\par\end{centering}

\protect\caption{\label{fig:numberedPointsB1}Board with numbered points}
\end{figure}

\noindent Using two-line notation for permutations, in which points
at the bottom are destinations for those at the top, we have for example
that:
\begin{equation}
\left(\begin{array}{ccccccccc}
1 & 2 & 3 & 4 & 5 & 6 & 7 & 8 & 9\\
2 & 1 & 3 & 5 & 4 & 6 & 8 & 7 & 9
\end{array}\right)\label{eq:a column permutation}
\end{equation}
swaps horizontally the points in the first and second columns, thus
leaving \emph{$C$} invariant: the permutation is then an automorphism
of the board. 
\end{example}
\noindent The automorphisms of a board form a group:
\begin{defn}
A group obtained under composition of automorphisms is an \textit{automorphism
group.} The group of all automorphisms of a board \textit{B} is the
\textit{full automorphism group} of \textit{B}, denoted by \emph{$Aut(B)$}.
\end{defn}
\noindent The set of all possible permutations of the elements of
a set also forms a group:
\begin{defn}
\label{def:The-symmetric-group}The \textit{symmetric group} on a
finite set \textit{$F$}, \emph{$Sym(F)$}, is the group whose elements
are all bijective functions from \textit{F} to \textit{F} with function
composition as the group operation. A subgroup of \emph{$Sym(F)$}
is called a \textit{permutation group}.
\end{defn}
\noindent As\emph{ $Sym(F)$} is the group of all permutations of
the elements of \textit{$F$,} its order is $|F|!$. This applies
to the set of points $P$ of a board $B$. Since every automorphism
of the board is a permutation of \textit{P}, it also belongs to \emph{Sym}(\textit{P}),
so $Aut(B)$ is a permutation group:
\begin{equation}
Aut(B)\subseteq Sym(P)\label{eq:Aut(B) included in Sym(P)}
\end{equation}

\subsection{Combinatorial Properties\label{sub:Combinatorial-Properties}}

\noindent Borrowing from the theory of Combinatorial Design \cite{Handbook of Combinatorial Designs},
we define now some concepts for boards. They will we used later in
the sections dealing with symmetric boards (Sects. \ref{sec:Symmetric-Boards}
and  \ref{sub:Latin-Polytopes}), where we will se that geometric
symmetry also induces combinatorial symmetry and balance.
\begin{defn}
A board is\textit{ k-uniform }if and only if all asterisms have size\textit{
k.}\end{defn}
\begin{example}
\noindent As all six asterisms in the board in Fig. \ref{fig:3x3points and vert and hor asterisms}
have size three, the board is $3$-uniform. The Fano Plane (see Fig.
\ref{fig:fanoPlaneUnsymmetric}) is also $3$-uniform.\end{example}
\begin{defn}
\label{def:parallel-class}A \textit{parallel class} or \textit{resolution
class }in a board \textit{$B=(P,C)$} is a set of asterisms that partitions
\textit{$P$.}
\end{defn}

\begin{defn}
A parallel class is \textit{$k$-uniform} if and only if all its asterisms
have size \textit{$k$}.
\end{defn}
\noindent A \textit{$k$}-uniform parallel class has then \textit{$\frac{|P|}{k}$}
asterisms and is a subset of the board constellation. A board can
thus have a \textit{$k$}-uniform parallel class only if \textit{$k$}\textit{\emph{
divides $|P|$,}} as the next example shows.
\begin{example}
\noindent \label{exa:H and V}In the board in Fig. \ref{fig:3x3points and vert and hor asterisms},
we see that each set $H$ and $V$, with
\[
\begin{gathered}H=\left\{ \left\{ 1,2,3\right\} ,\left\{ 4,5,6\right\} ,\left\{ 7,8,9\right\} \right\} ,\,\,\,\end{gathered}
V=\left\{ \left\{ 1,4,7\right\} ,\left\{ 2,5,8\right\} ,\left\{ 3,6,9\right\} \right\} 
\]
partitions the set of points and have asterisms of size three: each
is thus a $3$-uniform parallel class. It is immediate to see that
the Fano Plane (see Fig. \ref{fig:fanoPlaneUnsymmetric}) has no parallel
classes.\end{example}
\begin{defn}
In a board, two parallel classes \textit{$PC_{1}$, $PC_{2}$} are
said to be \textit{orthogonal} if and only if for any two asterisms
$a_{1}$, $a_{2}$, with $a_{1}\in PC_{1}$ and $a_{2}\in PC_{2}$,
\textit{\emph{$|a_{1}\cap a_{2}|=1$}}.\end{defn}
\begin{example}
\noindent As every asterism in \textit{$H$} in Example \ref{exa:H and V}
intersect at a single point every other in \textit{$V$,} both parallel
classes are orthogonal.\end{example}
\begin{defn}
\label{def:resolvable board}A \textit{resolution} of a board is a
set of parallel classes that partition the constellation. A board
with a resolution is said to be \textit{resolvable.}
\end{defn}

\begin{defn}
A resolution is \textit{$k$-uniform} if and only if each of its parallel
classes is \textit{k}-uniform.\end{defn}
\begin{example}
\noindent \textit{\emph{The board in}}\textit{ }Fig. \ref{fig:3x3points and vert and hor asterisms}
is resolvable, as it has a $3$-uniform resolution with two orthogonal
parallel classes. 
\end{example}
\noindent In $\mathbb{R}^{2}$ orthogonality among geometric lines
requires both intersection and equality of angles. Parallelism and
orthogonality are then special cases of a concept more appropriate
to boards:
\begin{defn}
\label{def:SIN}Given a board\textit{ $B=(P,C)$} we call $\left\{ |a\cap a\lyxmathsym{\textquoteright}|:a,a\lyxmathsym{\textquoteright}\in C,a\neq a\lyxmathsym{\textquoteright}\right\} $
its \textit{set of intersection numbers }or \textit{SIN}.
\end{defn}
\noindent It is trivial to form one or more parallel classes from
a general set of points, but not so if we demand a particular SIN
in the resulting board. 
\begin{example}
\noindent In the board in Fig. \ref{fig:3x3points and vert and hor asterisms},
the tiled square provides an interesting intersection pattern from
the outset: parallelism (within each parallel class) and orthogonality
(among classes), i.e. a pattern with a SIN $=\left\{ 0,1\right\} $.
As any two asterisms in the Fano Plane either have one point in common
or none at all (see Fig. \ref{fig:fanoPlaneUnsymmetric}), its SIN
is also $\left\{ 0,1\right\} $.
\end{example}

\begin{figure}
\begin{centering}
\includegraphics[height=4cm]{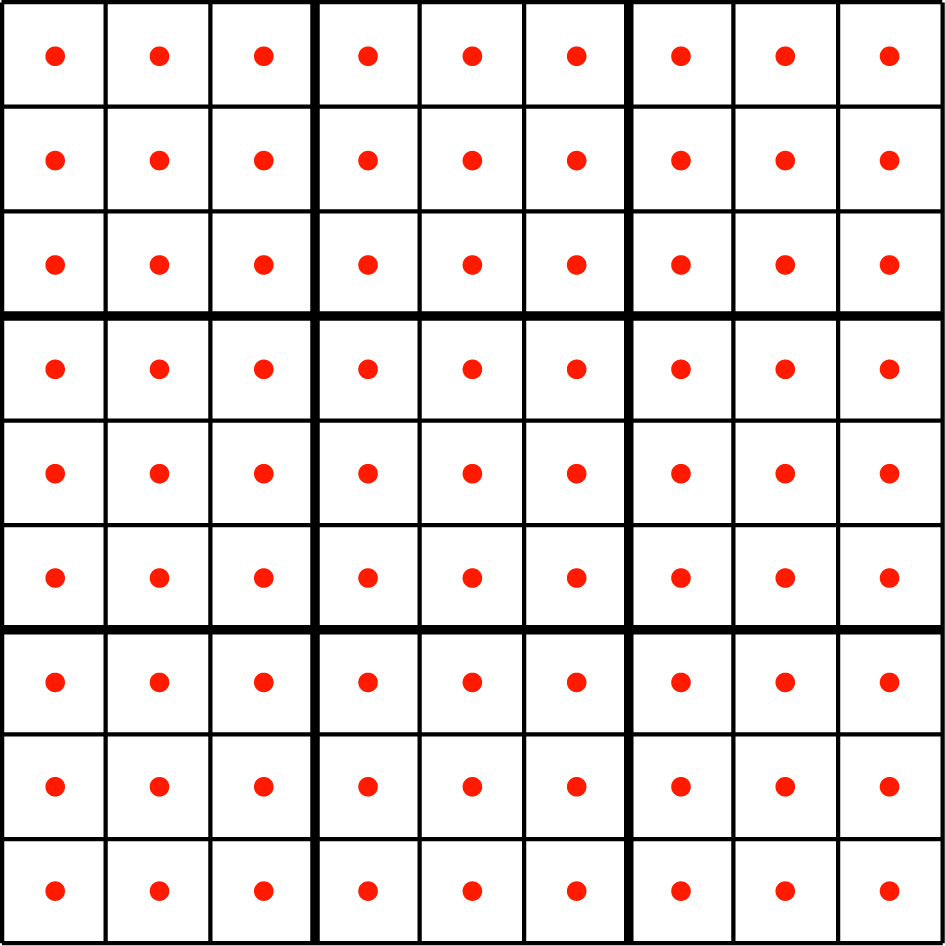}
\par\end{centering}

\protect\caption{\label{fig:sudoku design pts-1}Board with eighty one points and twenty
seven asterisms}
\end{figure}

\begin{example}
Fig. \ref{fig:sudoku design pts-1} shows a board with $81$ points
and $27$ asterisms, each with points lying on each\end{example}
\begin{itemize}
\item row ($9$ asterisms that we group in set \textit{$H_{S}$}) 
\item column ($9$ asterisms that we group in set \textit{$V_{S}$}) 
\item $3\times3$ sub-square ($9$ asterisms that we group in set \textit{$Q_{S}$}) 
\end{itemize}

\noindent As \textit{$H_{S}$}, \textit{$V_{S}$} and \textit{$Q_{S}$}
partitions each the set of points, each is a $9$-uniform parallel
class and the board is resolvable. \textit{$H_{S}$} and \textit{$V_{S}$}
are orthogonal, whereas each asterism in \textit{$Q_{S}$} intersects
in three points three asterisms\textit{ }\textit{\emph{in}}\textit{
$H_{S}$} and three in \textit{$V_{S}$}, having an empty intersection
with the rest.

\section{Latin Boards\label{sec:Latin-Boards}}

In this section we assign labels in a particular way to the points
of a board to create a new type of object.

\subsection{Definition of Latin Board}
\begin{defn}
\label{def:Latin board definition} Let $B=\left(P,C\right)$ be a
board, \emph{L }a $k$-multiset of \emph{labels}, $n_{i}$\emph{ }the
size of asterism $a_{i}$, and $q_{i},r_{i}$ the respective quotient
and remainder of $\frac{n_{i}}{k}$. A\emph{ Latin board} is a tuple
$\left(B,L,F\right)$ where $F$ is a function $P\rightarrow L'$,
with $L'\subseteq\left\{ l,l\in L\right\} $, such that the points
of every asterism $a_{i}$ in \emph{B} are labeled with $q_{i}$ whole
instances of $L$ plus any $r_{i}$ elements of it.
\end{defn}
\noindent In simpler terms: each asterism has as many whole copies
of the multiset as it can hold, plus some elements of the multiset
to pad the remaining cells. We also say that a Latin board is the
result of \emph{latinizing} or \emph{labeling} a board with a multiset\footnote{\noindent Multisets generalize sets in that any element count may
be more than one.}. 

According to this definition, and in contrast with Latin squares (see
Def. \ref{def:A-Latin-square}), Latin boards have no constraints
on topology, number of asterisms, equality of asterism sizes, uniqueness
of labels or number of labels relative to asterism sizes. The leeway
provided by the definition in the labeling of the remainder points
in every asterism is illustrated by the following example. This flexibility
can be used to define types of Latin boards, for example by specifying
the remaining labels on a per asterism basis.
\begin{example}
Let \{P, E, A, C, E\} be our multiset of labels and let $3,5,7,10,11$
be the respective number of points in the asterisms of a given board.
For the purpose of the example we will order the points on each asterism
in a sequence. Then the following sequences, or their distinct permutations,
are some of the respective possible latinizations of each asterism:\end{example}
\begin{itemize}
\item {[}P, E, A{]}, {[}A, E, E{]} (any $3$ labels of the multiset)
\item {[}A, C, E, E, P{]} (the full multiset)
\item {[}P, E, A, C, E, E, C{]}, {[}E, P, A, C, E, E, E{]} (one copy of
the multiset plus any two of its labels)
\item {[}P, E, A, C, E, P, E, A, C, E{]} (two copies of the multiset)
\item {[}P, E, A, C, E, P, E, A, C, E, P{]}, {[}P, E, A, C, E, E, C, A,
E, P, E{]} (two copies of the multiset plus any of its labels)\end{itemize}
\begin{example}
Fig. \ref{fig:non-unifit Latin board} (left) shows a board with seven
points and a constellation with three asterisms, each one formed by
points joined by segments of the same color. The board is not uniform
because not all asterisms have the same size: the ``green'' asterism
has seven points, the orange three and the blue four. 

If we want to latinize this board with the multiset $\left\{ 1,1,2,2\right\} $
we have to label the green asterism with one whole multiset plus any
three labels (since $\frac{7}{4}=1+\frac{3}{4}$), for example with
multiset $\left\{ 1,1,2,2,1,2,2\right\} $ or $\left\{ 1,1,2,2,1,1,2\right\} $.
Similarly, we have to label the blue asterism with just one whole
multiset ($\frac{4}{4}=1+\frac{0}{4}$); and the orange one with any
three labels (since $\frac{3}{4}=0+\frac{3}{4}$), for example with
$\left\{ 1,2,2\right\} $ or $\left\{ 1,1,2\right\} $. Fig. \ref{fig:non-unifit Latin board}
(center and right) shows two example Latin boards. 
\begin{figure}
\begin{centering}
\includegraphics[height=2.5cm]{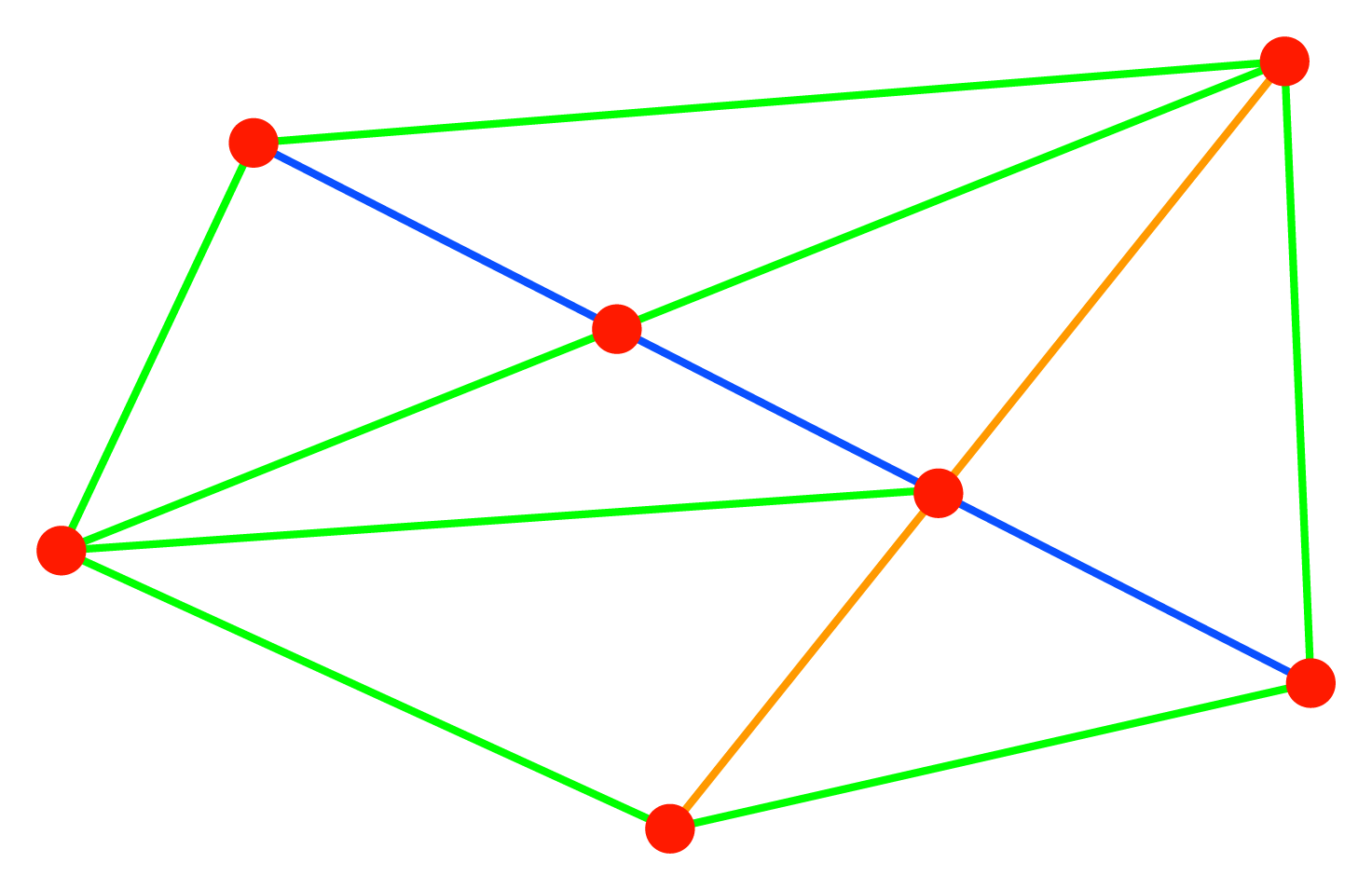}~~\includegraphics[height=2.5cm]{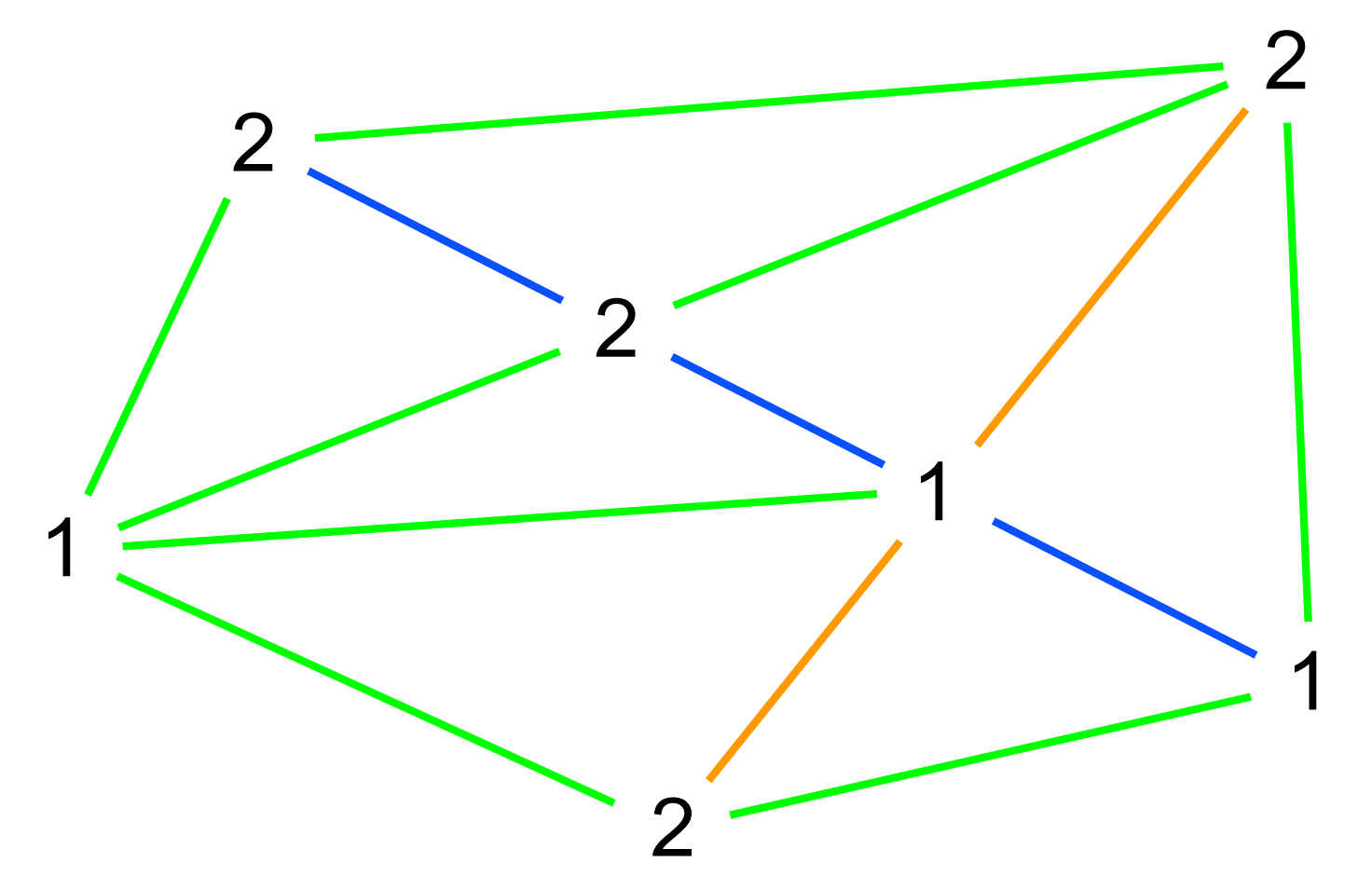}~~\includegraphics[height=2.5cm]{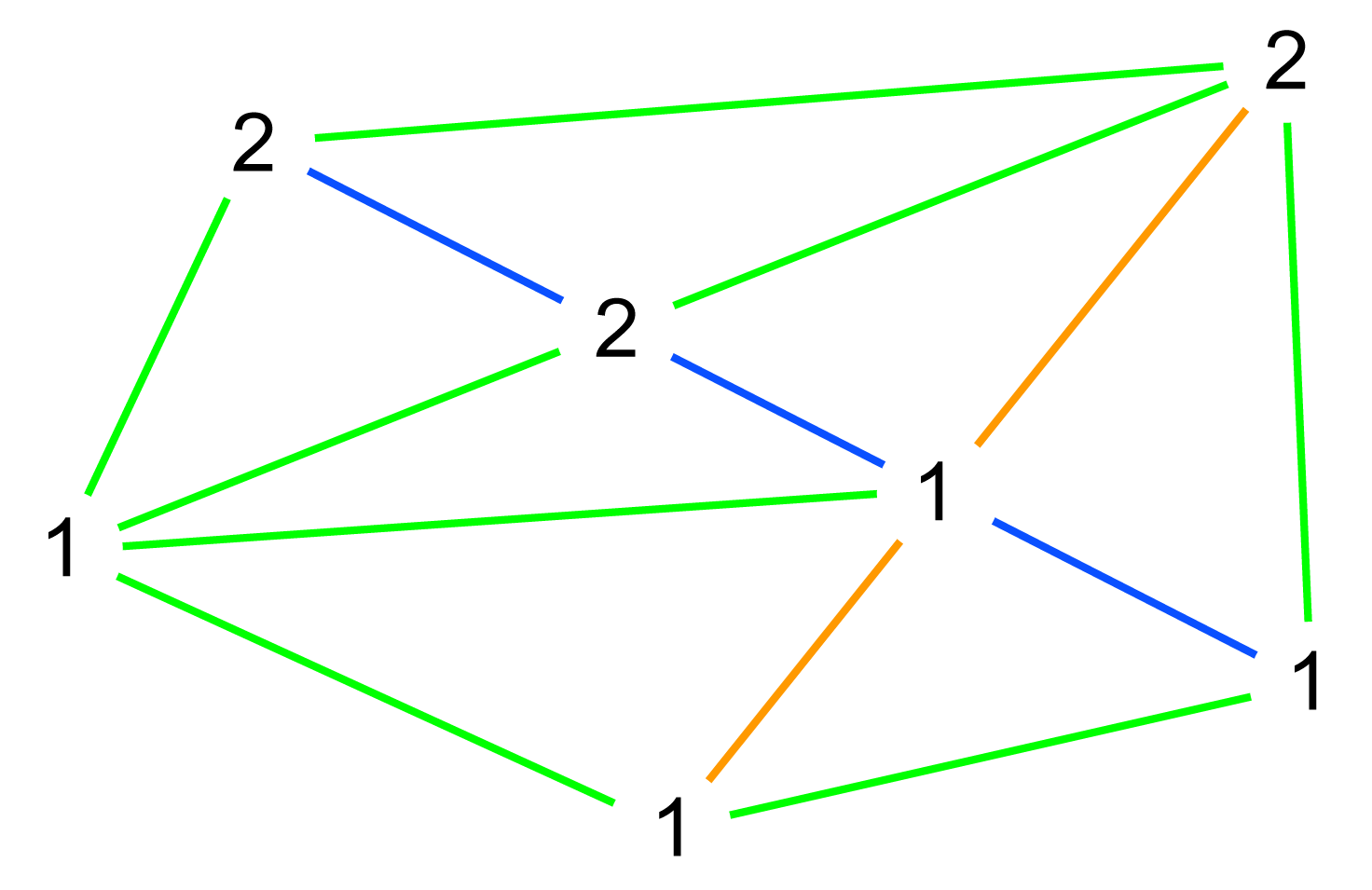}
\par\end{centering}

\protect\caption{\label{fig:non-unifit Latin board}A board and two Latin boards}
\end{figure}

\end{example}

\subsection{Incidence Structures Related to Latin Boards\label{sub:Incidence Structures Related to Latin Boards}}

We see here how Latin boards relate to the incidence structures mentioned
in Sect. \ref{sub:Incidence Structures-Related-to-boards}.

\subsubsection{Families of Sets}

In Sect. \ref{sub:Families-of-Sets} we identified the elements and
subsets in a family of sets with the points and asterisms of a board.
Latinization can be viewed then as the problem of finding an additional
set of subsets that partitions the elements of the family, with the
property that each one intersects each initial subset in the family
at exactly the same number of elements; with the number of elements
in the intersection being the count of the corresponding label in
the multiset.

\subsubsection{Finite Geometries}

Latinizing a board that is also a finite geometry means to find subsets
of points that meet the requirement imposed by the multiset. These
new subsets will not be lines in the geometry in general.

\subsubsection{Simple Combinatorial Designs}

Once the correspondence between boards and simple combinatorial designs
have been established as described in Sect. \ref{sub:Simple-Combinatorial-Designs},
latinizing a board is akin to finding specific parallel classes (see
Def. \ref{def:parallel-class}) in the design.

\subsubsection{Hypergraphs}

Vertices on an hypergraph may be assigned labels called \emph{colors}.
Erdös and Hajnal defined classic hypergraph coloring to generalize
graph coloring \cite{On chromatic number of graphs and set systems}:
\begin{defn}
Let $\left\{ 1,2,\ldots,n\right\} $ be a set of \emph{colors}. A
\emph{proper $n$-coloring} of a hypergraph is a labeling of its vertices
with the colors in such a way that every hyperedge with a size of
at least two has at least two vertices colored differently.\label{def:classic hg coloring}
\end{defn}
\noindent If we choose $\left\{ 1,2,3\right\} $ as a set of colors,
then the Latin square in Fig. \ref{fig:latinSquareOrder3} is a valid
classic coloring of the hypergraph in Fig. \ref{fig:3x3points and vert and hor asterisms}
(right)\footnote{\noindent It is immediately seen that not any coloring is a Latin
square though.}. We also notice that if the multiset used in a particular Latin board
has at least two different labels, then latinization is a more general
and demanding task than classic hypergraph coloring.

A more general view of hypergraph coloring admits any possible partition
of the set of vertices into classes of the same color \cite{Coloring Mixed Hypergraphs}.
As Latin boards are geometric hypergraphs, and since latinization
induces a partition of the points into classes holding points with
the same label, latinization is an instance of this general hypergraph
coloring.

\subsection{Unifit Latin Boards\label{sub:Unifit-Latin-Boards}}

\noindent Unless otherwise stated, we will only deal here with boards
of a specific type:
\begin{defn}
Let $B=\left(P,C\right)$ be a $k$-uniform board and \emph{L }a $k$-multiset
of \emph{labels}. A\emph{ Latin board} is a tuple $\left(B,L,F\right)$
where $F$ is a bijective function $P\rightarrow L$.
\end{defn}
\noindent This is a restricted version of Def. \ref{def:Latin board definition}.
An equivalent simpler definition is:
\begin{defn}
A \emph{k-unifit Latin board} is a $k$-uniform board latinized with
a $k$-multiset.\end{defn}
\begin{example}
\noindent \label{exa:3x3 free Latin squares}Fig. \ref{fig:another Latin square and a free one with multiset}
(left) shows a 3-unifit Latin board resulting from the latinization
of the board in Fig. \ref{fig:3x3points and vert and hor asterisms}
with $\left\{ 1,2,3\right\} $. This is identical to the Latin square
in Fig. \ref{fig:latinSquareOrder3}, so this particular Latin square
is a Latin board. In fact it is simple to prove that any Latin square
is a Latin board. 

Other latinizations of the same board are possible: if we write any
permutation of the labels in the top row, then shift them by one cell
to the left in successive rows, we obtain Latin boards that are also
Latin squares (see Fig. \ref{fig:another Latin square and a free one with multiset}
center). Fig. \ref{fig:another Latin square and a free one with multiset}
(right) shows a Latin board with multiset \{1, 1, 2\}. This proves
that not all latinizations of the board in Fig. \ref{fig:3x3points and vert and hor asterisms}
are Latin squares.
\begin{figure}
\begin{centering}
\includegraphics[height=2.5cm]{latinSquareOrder3}~~~~~~~~~~\includegraphics[height=2.5cm]{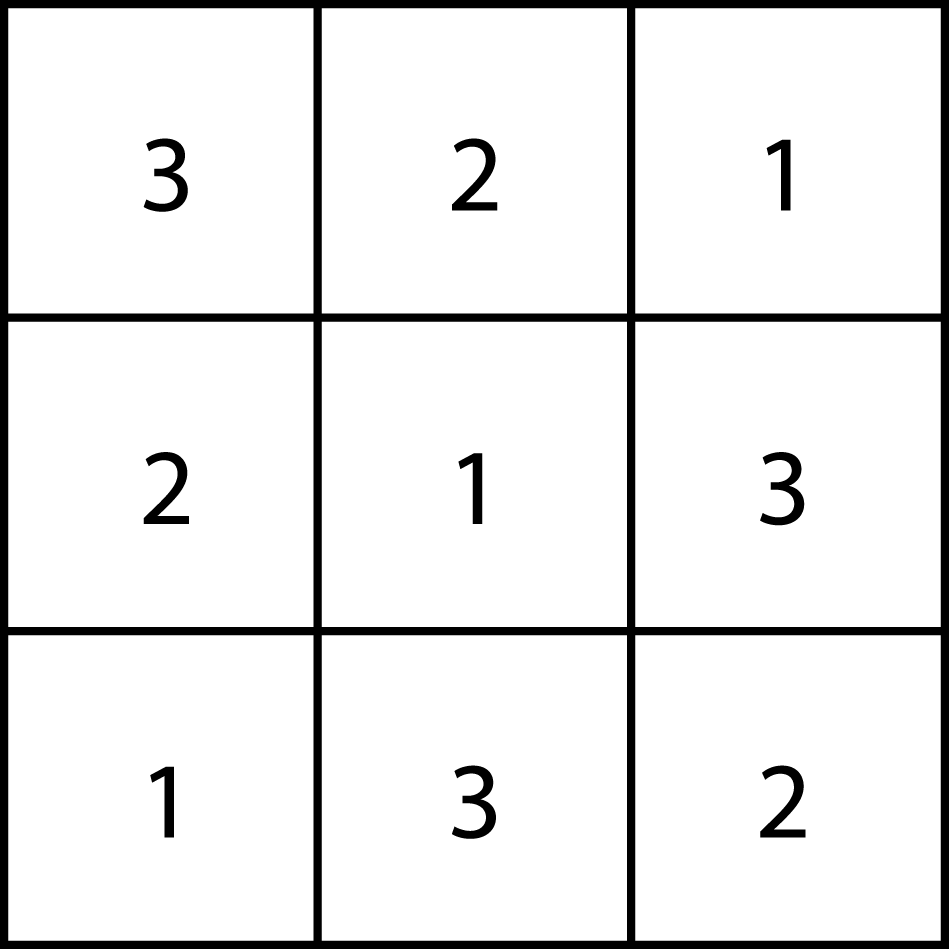}~~~~~~~~~~\includegraphics[height=2.5cm]{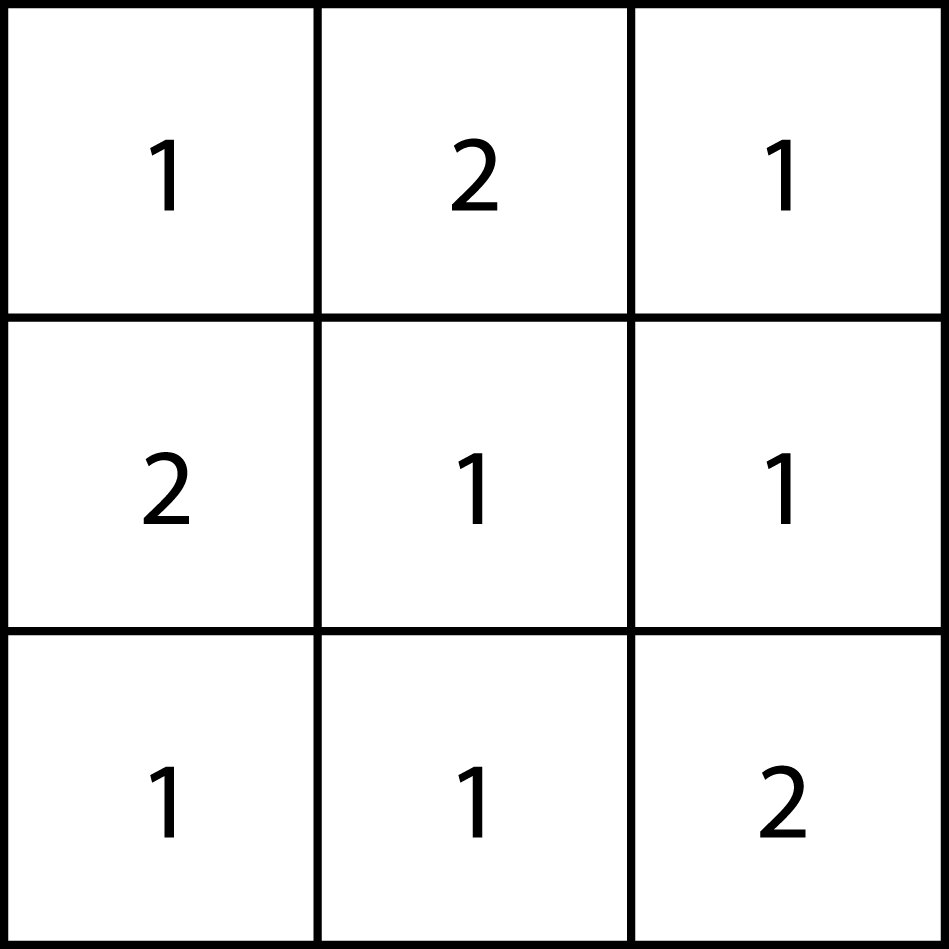}
\par\end{centering}

\protect\caption{\label{fig:another Latin square and a free one with multiset}Different
latinizations of the same board}
\end{figure}

\end{example}

\begin{example}
\noindent The triangular and hexagonal objects in Figs. \ref{fig:Latin triangle example}
and \ref{fig:Latin hexagon example} are unifit Latin boards with
respective sets $\left\{ 1,2,3,...,12\right\} $ and $\left\{ 1,2,3,...,18\right\} $.
They are examples of what will be called \emph{Latin triangles} and
\emph{Latin hexagons} later. As indicated in Sect. \ref{par:Points-and-Cells-in-boards},
the points here are the centroids of the cells. The superimposed graphical
structure \textendash the grids of squares and triangles\textendash{}
is just a way to rendering the board points and asterisms apparent.
\end{example}
\noindent In what follows we will usually drop the adjective ``unifit''.

\pagebreak{}

\subsection{Latinization Techniques}

\noindent Latinizations in Fig. \ref{fig:another Latin square and a free one with multiset}
were easily found by inspection, but latinizing a generic board is
not trivial. Several methods are available \cite{Latin Puzzles paper};
the author used \emph{constraint programming} techniques \cite{Principles of Constraint Programming Apt}
to generate the examples that follow.
\begin{example}
\noindent Fig. \ref{fig:LSandCriticalSet-1} (right) shows a Latin
square of order seven.
\end{example}

\begin{example}
\noindent The next figures show three different latinizations of the
board in Fig. \ref{fig:sudoku design pts-1}: Fig. \ref{fig:Sudoku}
(right), labeled with the set $\left\{ 1,2,3,4,5,6,7,8,9\right\} $;
Fig. \ref{fig:Sudoku Ripeto} (right), latinized with the multiset
$\left\{ 1,1,1,2,2,2,3,3,3\right\} $ and Fig. \ref{fig:Custom Sudoku}
(right) with the multiset \{N, E, O, S, U, D, O, K, U\}. 
\begin{figure}
\begin{centering}
\includegraphics[height=4.5cm]{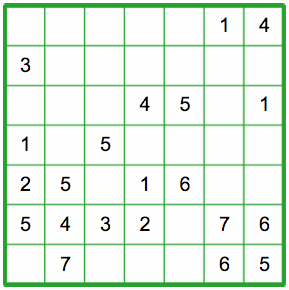}~~~~~~\includegraphics[height=4.5cm]{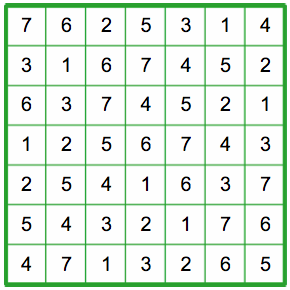}
\par\end{centering}

\protect\caption{\label{fig:LSandCriticalSet-1}Latin square puzzle and solution}
\end{figure}
\begin{figure}
\centering{}\includegraphics[height=4.5cm]{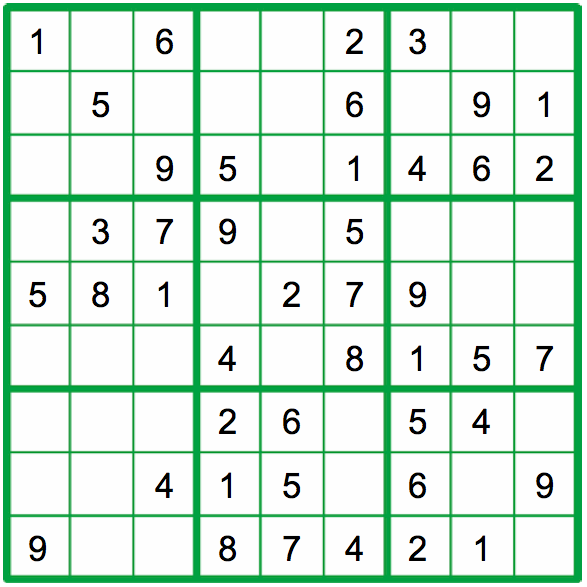}~~~~~~\includegraphics[height=4.5cm]{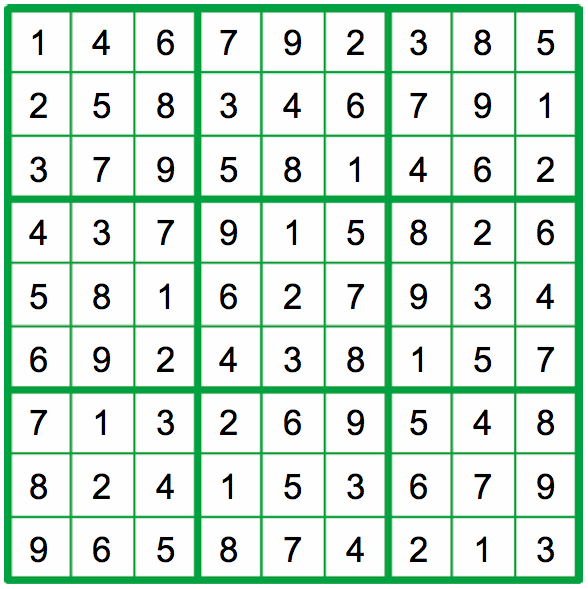}\protect\caption{\label{fig:Sudoku}\emph{Sudoku} puzzle and solution}
\end{figure}
\begin{figure}
\begin{centering}
\includegraphics[height=4.5cm]{Sudoku_Repeto_Clues}~~~~~~\includegraphics[height=4.5cm]{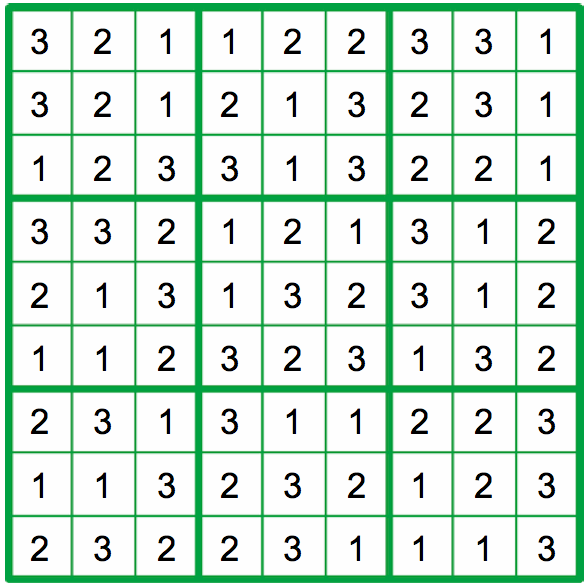}
\par\end{centering}

\protect\caption{\label{fig:Sudoku Ripeto}\emph{Sudoku Ripeto} puzzle and solution}
\end{figure}
\begin{figure}
\begin{centering}
\includegraphics[height=4.5cm]{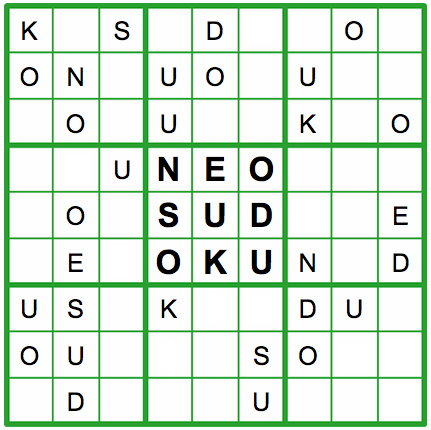}~~~~~~\includegraphics[height=4.5cm]{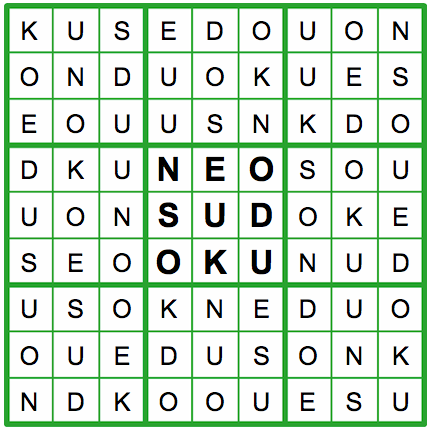}
\par\end{centering}

\protect\caption{\label{fig:Custom Sudoku}\emph{Custom Sudoku} puzzle and solution}
\end{figure}

\end{example}

\subsubsection{Latinizing the Fano Plane}
\begin{figure}
\begin{centering}
\includegraphics[height=3.4cm]{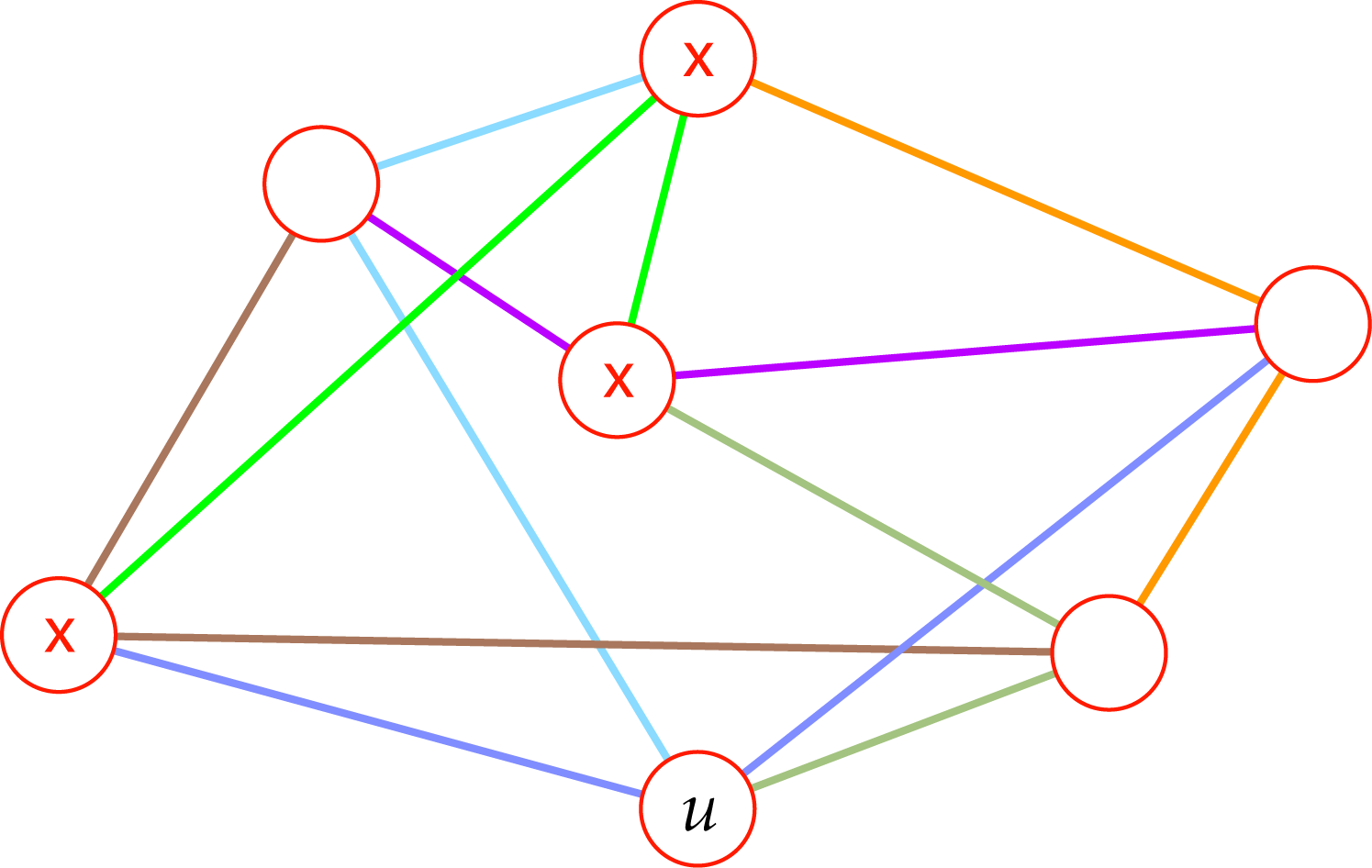}
\par\end{centering}

\protect\caption{\label{fig:FanoPlaneUnsymmetricUnnumberedWoneLabel}Attempt to latinize
the Fano Plane}
\end{figure}

\begin{defn}
In a \emph{trivial latinization} all labels in the board are equal.\end{defn}
\begin{thm}
The only latinization with a $3$-multiset admitted by the Fano Plane
is the trivial one.\end{thm}
\begin{proof}
As the asterisms in the Fano Plane have size three, there are only
two types of candidate $3$-multisets for latinization: one with at
least one unique label \textendash this one includes the case in which
all labels are different\textendash{} and another with all labels
equal. Let \emph{u} be the unique label in a multiset of the first
type. Fig. \ref{fig:FanoPlaneUnsymmetricUnnumberedWoneLabel} shows
the Fano Plane with one of its points labeled with this label. The
``light green'' asterism should have also a point labeled \emph{u},
but that point cannot be its leftmost one because it is being shared
with the dark blue line, already holding the label. The same applies
to the two other points in the light green line, due to the light
blue and dark green lines. As any point in the Fano Plane has the
same connectivity pattern, the argument is valid for any point, so
the only possibility is to use a multiset with all labels equal, i.e.:
one that produces the trivial latinization.
\end{proof}

\subsection{Latin Puzzles\label{sub:Latin-Puzzles}}

Latin squares are redundant structures: they contain more information
than the strictly necessary to describe them, as the next example
shows.
\begin{example}
If we remove all labels in the top row of the Latin square in Fig.
\ref{fig:latinSquareOrder3} the result is completable to \emph{just}
the original square. 
\end{example}
\noindent Squares with missing labels having the property of being
completable to a single Latin square are called \emph{uniquely completable
partial Latin squares} \cite{Latin squares and their applications}
and \emph{Latin square puzzles} in \cite{Latin Puzzles paper}. 
\begin{example}
Fig. \ref{fig:LSandCriticalSet-1} shows a Latin square (right) and
an associated Latin square puzzle (left).
\end{example}
\noindent The corresponding objects for Latin boards are called\emph{
Latin puzzles,} and are defined likewise \cite{Latin Puzzles paper}.
Methods to find Latin puzzles are intimately related to those used
to find Latin boards.
\begin{example}
The following examples were found by the author using constraint programming
techniques \cite{Latin Puzzles paper}: Fig. \ref{fig:Sudoku} (left)
shows a \emph{Sudoku} puzzle \cite{Taking Sudoku Seriously}. Fig.
\ref{fig:Sudoku Ripeto} (left) is an example of \emph{Sudoku Ripeto}
puzzle \cite{Latin Puzzles site}. Fig. \ref{fig:Custom Sudoku}\emph{
}(left) shows a\emph{ Custom Sudoku} puzzle \cite{Latin Puzzles site}.
More Latin puzzles may be found in \cite{Latin Puzzles paper,Latin Puzzles site}.
\end{example}

\subsection{Equivalence Classes\label{sub:4.1-Equivalence-classes-Latin boards}}

\noindent Latin boards may be easily derived from a given one in some
cases. We illustrate this with examples involving Latin squares, \emph{Sudokus}
and Latin boards in general.

\subsubsection{Equivalence Classes of Latin Squares\label{sub:Equivalence-classes-of-Latin squares}}

It is easily seen that, if we permute the labels, the rows or the
columns of a Latin square we obtain another Latin square. Both squares
are said to be \textit{isotopic}. \textit{Isotopism} is an equivalence
relation that partitions the set of Latin squares of a particular
order into \textit{isotopy classes} \cite{Latin squares and their applications}.

On the other hand, a Latin square of order \textit{n} with labels
\{1,\ldots{},\textit{n}\} may be represented with an \textit{n}\textsuperscript{2}-set
of $3$-tuples (\textit{i}, \textit{j}, \textit{k}), one for each
cell, where \textit{i} is the index of the row, \textit{j} that of
the column and \textit{k} the label in the cell. This set, called
the \textit{orthogonal array representation} of the square \cite{Latin squares and their applications}
makes apparent that any pair of components in the tuple may act as
row and column indices in yet another Latin square, with the remaining
one becoming the corresponding cell content. As there are six possible
ways to do this, from one Latin square we can obtain a maximum of
five others\footnote{Or less, since some permutations may result in the same Latin square.}
called \textit{conjugates }of the original one.

Both operations can be combined: two Latin squares are said to be
\textit{paratopic} if one of them is isotopic to a conjugate of the
other. This is an equivalence relation that creates \textit{paratopy
classes}, each one containing up to six isotopy classes\footnote{For the number of paratopy classes of Latin squares as a function
of the order see \cite{Handbook of Combinatorial Designs,Number of Paratopy classes of Latin Squares}.}.

\subsubsection{\noindent Equivalence Classes of \emph{Sudokus}}

\noindent Similarly, from a \emph{Sudoku} (see Fig. \ref{fig:Sudoku}
right) we can easily obtain new ones by any combination of the following
operations:
\begin{itemize}
\item a permutation of the labels 
\item a permutation of rows within a 3-row block
\item a permutation of the 3 row blocks 
\item a permutation of columns within a 3-column block
\item a permutation of the 3 column blocks 
\item a reflection across the vertical or horizontal symmetry axis
\item a $90\text{º}$, $180\text{º}$ or $270\text{º}$ rotation around
the center of the board
\end{itemize}
We can likewise say that two \emph{Sudokus} are \emph{equivalent}
if we can transform one into the other using a combination of the
mentioned operations\footnote{In fact, not all of them are necessary to generate all equivalent
boards, just a subset of the operations.}. We can divide then the set of \emph{Sudoku} boards into several
equivalence classes\footnote{Enumerating the classes is a somewhat involved process, see \cite{Enumerating possible Sudoku grids}.}.

\subsubsection{Equivalence Classes of Latin Boards}

If we can transform a Latin board into another one with a set of operations
we can define a notion of equivalence among boards. Each set of operations
will generate different equivalence classes. The more elements in
the set, the less equivalence classes there will be and the more boards
per class will result.

\noindent Equivalence classes are important when enumerating Latin
boards. If several boards are found with a particular method \cite{Latin Puzzles paper},
they will be more valuable if they belong to different equivalence
classes, since if they are not, from one of them we can easily find
the others using the available operations. To avoid unnecessary calculations,
algorithms that enumerate or count Latin boards can eliminate the
need to look for boards in the same class of those already found \cite{Handbook of Combinatorial Designs,Latin Puzzles paper},
or eliminate dead ends in the search for that matter. This needs as
a key element the set of operations that generates the classes, that
should be as large as possible for the reasons just explained.

As we will see in Sect. \ref{sub:Equivalence-Classes-of-Latin-Polytopes},
when a Latin board is symmetric and there are no other equivalence
operations available, those linked to its symmetry group can be used
to generate equivalence classes.

\subsection{Woven Boards\label{sec:Woven-boards}}

Given a Latin board, we can naturally extend the underlying board
by adding to it an asterism for each set of points with the same label.
\begin{defn}
\label{def:warp class}Let $B=\left(\left\{ P,C\right\} ,L,M\right)$
be a Latin board. Let \textit{$C\lyxmathsym{\textquoteright}=\left\{ a_{l},l\in L\right\} $}
be a constellation each of whose asterisms contains all points with
the same label. Then the board $W=\left(P,C\cup C\lyxmathsym{\textquoteright}\right)$
is the \emph{woven board} of \emph{$B$}.
\end{defn}
\noindent Latin boards with different permutations of the labels have
then the same woven board. In this context, the asterisms in $C$
are called\emph{ weft asterisms} and $C$ proper, the \emph{weft constellation}.
The asterisms in $C\text{\textquoteright}$ are the\emph{ warp asterisms,}
while $C\text{\textquoteright}$ proper is the \emph{warp constellation}.
Like in a loom, in a woven board the weft constellation is interwoven
with the warp one.
\begin{example}
Fig. \ref{fig:Latin board and its woven board} (left) shows a latinization
of the board in Fig. \ref{fig:3x3points and vert and hor asterisms}
with $\left\{ 1,2,3\right\} $. On the right, its woven board. Solid
colored segments indicate weft asterisms; warp asterisms are represented
by dashed segments. Woven boards obtained from Latin squares are combinatorial
designs known as $3$\emph{-nets} \cite{Handbook of Combinatorial Designs}.
\begin{figure}
\begin{centering}
\includegraphics[height=2.5cm]{latinSquareOrder3}~~~~~~~~~~\includegraphics[height=2.5cm]{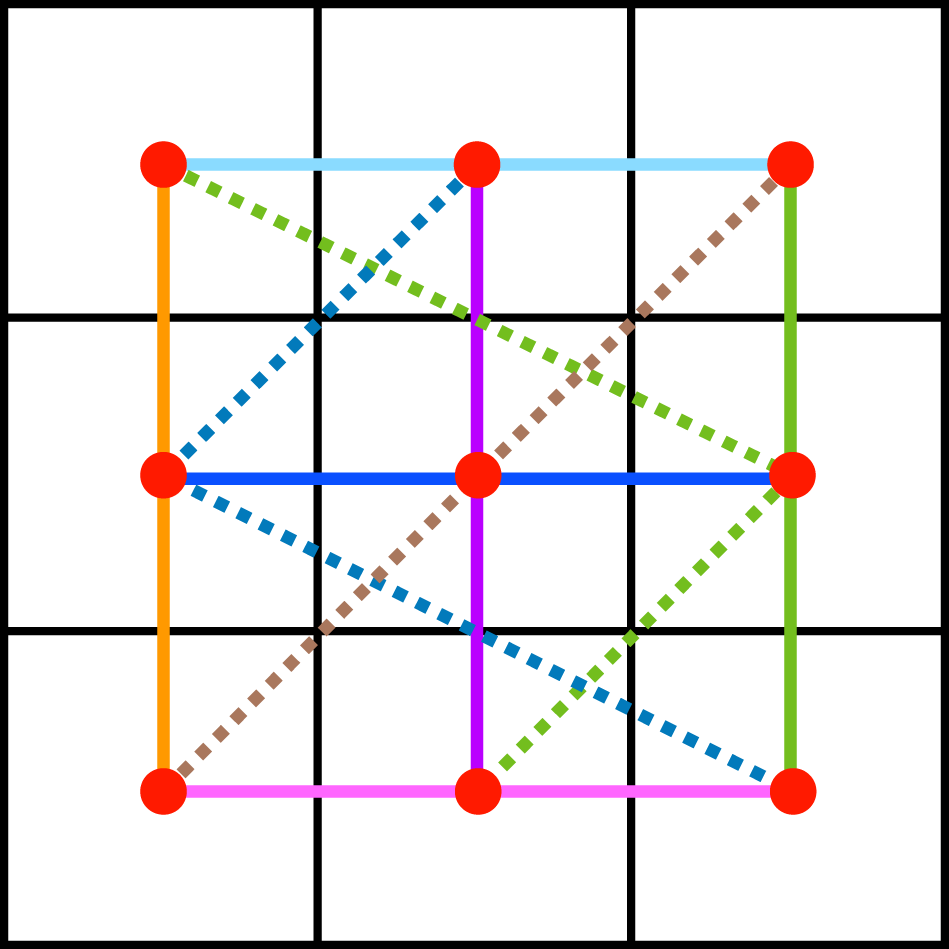}
\par\end{centering}

\protect\caption{\label{fig:Latin board and its woven board}Latin board and woven
board}
\end{figure}
\end{example}
\begin{lem}
The warp constellation\emph{ }of a woven board is a parallel class.\end{lem}
\begin{proof}
As all points in the originating Latin board are labeled, the warp
asterisms partition them, so they constitute a parallel class.\end{proof}
\begin{lem}
\label{lem:size of warp asterism}Let $B=(P,C)$ be a unifit Latin
board. If the weft constellation in its woven board contains a $k$-uniform
parallel class, then the size of each warp asterism is $\frac{|P|}{k}\text{·}c$,
where c is the count, in the multiset of labels, of the label the
asterism was derived from.\end{lem}
\begin{proof}
As the parallel class partitions $P$, the number of weft asterisms
in the class is $\frac{|P|}{k}$. As the Latin board is unifit, a
label with count $c$ in the multiset will be present in $c$ points
on each weft asterism, so the warp asterism corresponding to this
label has $\frac{|P|}{k}\text{·}c$ points. \end{proof}
\begin{thm}
If a unifit Latin board with \emph{$|P|$} points has a k-uniform
parallel class, and each label in the multiset has count $\frac{k^{2}}{|P|}$,
then its woven board is also k-uniform. \end{thm}
\begin{proof}
By Lemma \ref{lem:size of warp asterism} the size of each warp asterism
is $\frac{|P|}{k}\text{·}c=\frac{|P|}{k}\text{·\ensuremath{\frac{k^{2}}{|P|}}=\emph{k}}.$
As weft asterisms are also $k$-uniform, the woven board is \emph{also
k}-uniform.
\end{proof}
\noindent The repetition pattern of labels in the multiset serves
then to fine-tune the combinatorial properties of the resulting woven
board, like the number and size of the warp asterisms and its SIN
(set of intersection numbers, see Def. \ref{def:SIN}).

\subsection{Sequences of Latin Boards}

\subsubsection{Sequences of Latin Squares}

\begin{figure}
\begin{centering}
\includegraphics[height=2.5cm]{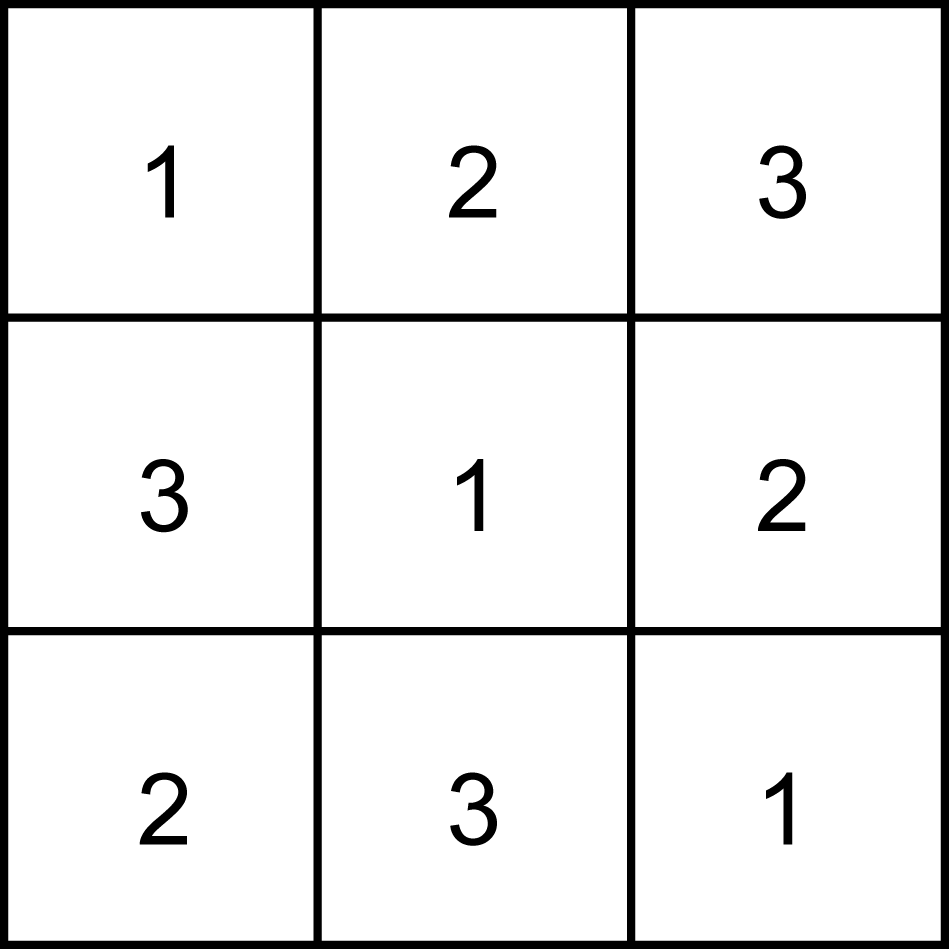}
\par\end{centering}

\protect\caption{\label{fig:latinize-Weave sequence Latin Squares}Latinized woven
board}
\end{figure}

\begin{figure}
\begin{centering}
~~~\includegraphics[height=3.9cm]{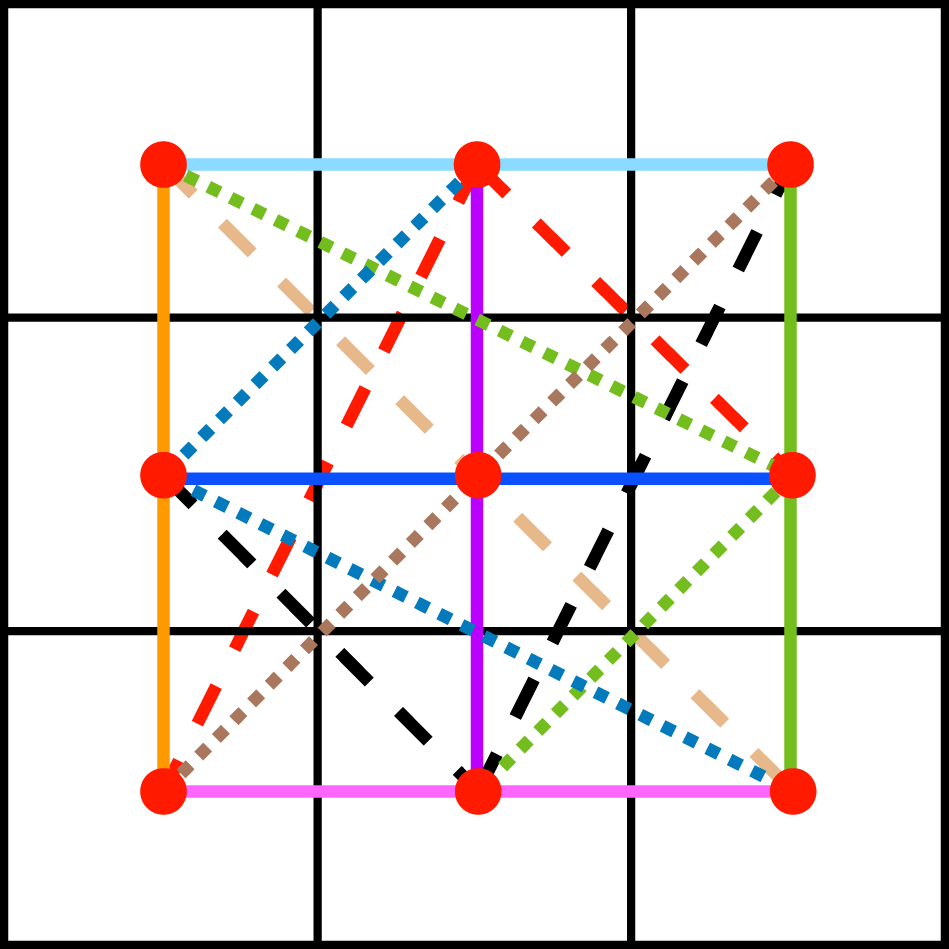}
\par\end{centering}

\protect\caption{\label{fig:second woven board}Second woven board}

\end{figure}

If we latinize the woven board in Fig. \ref{fig:Latin board and its woven board}
(right) with $\left\{ 1,2,3\right\} $ we can obtain the Latin Square
in Fig. \ref{fig:latinize-Weave sequence Latin Squares}. Weaving
in turn this second Latin square we obtain the woven board in Fig.
\ref{fig:second woven board}, where the three new warp asterisms
are represented with even more spaced dashed lines. This board features
four parallel classes, each with three asterisms of three points each.
If we identify the points as in Fig. \ref{fig:numberedPointsB1},
the classes are:\renewcommand{\arraystretch}{1.3}

\[
\begin{array}{c}
\left\{ \left\{ 1,2,3\right\} ,\left\{ 4,5,6\right\} ,\left\{ 7,8,9\right\} \right\} \\
\left\{ \left\{ 1,4,7\right\} ,\left\{ 2,5,8\right\} ,\left\{ 3,6,9\right\} \right\} \\
\left\{ \left\{ 1,6,8\right\} ,\left\{ 2,4,9\right\} ,\left\{ 3,5,7\right\} \right\} \\
\left\{ \left\{ 1,5,9\right\} ,\left\{ 2,6,7\right\} ,\left\{ 3,4,8\right\} \right\} 
\end{array}
\]
\renewcommand{\arraystretch}{1.5}

\noindent These classes have an interesting property: any two asterisms,
each from a different class, intersect at exactly one point. The two
intervening Latin squares in the sequence constitute a set of \emph{mutually
orthogonal Latin squares} (MOLS) \cite{Handbook of Combinatorial Designs}.
Can we extend the just constructed \emph{Latin sequence}? It is known
that the maximum size of a set of MOLS of order $n$ is $n-1$ \cite{Handbook of Combinatorial Designs},
so no further extensions of the sequence are possible.

\noindent 

\subsubsection{Sequences of Free Latin Squares}

\begin{figure}
\begin{centering}
\includegraphics[height=4.5cm]{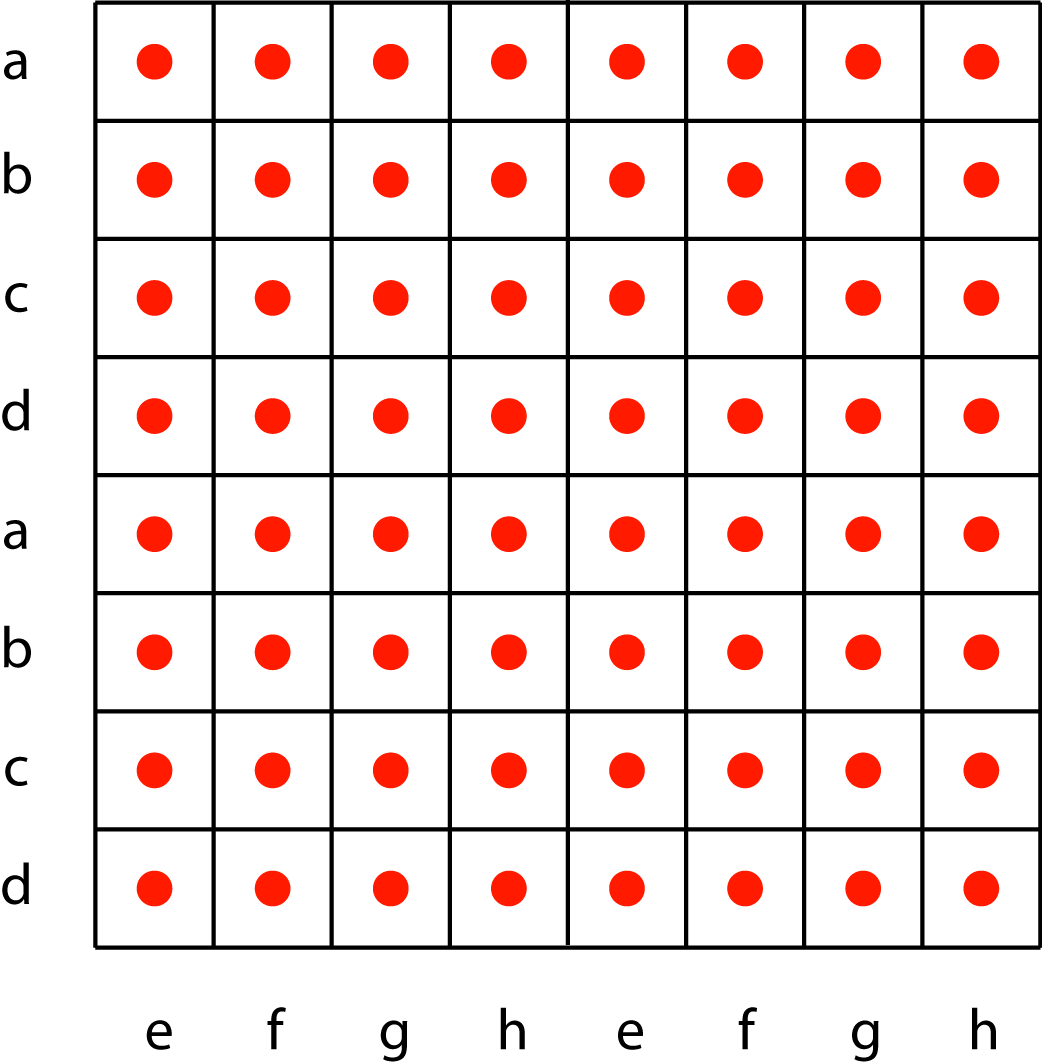}
\par\end{centering}

\protect\caption{\label{fig:board Free Latin Square}Board with SIN $\left\{ 0,4\right\} $}
\end{figure}

The same process can be applied to more general Latin boards. In the
next example we construct a Latin sequence that illustrates the generalization
of orthogonality to \textit{uniform intersection} among parallel classes.
Let \textit{$P$} be the set of points in Fig. \ref{fig:board Free Latin Square}.
Let\textquoteright s make an asterism out of points indexed by the
same letter: each asterism has then the points on either two rows
or two columns. The result is a board with SIN $\left\{ 0,4\right\} $
and parallel classes: $\left\{ a,b,c,d\right\} $ and $\left\{ e,f,g,h\right\} $.

\begin{figure}
\begin{centering}
\includegraphics[height=3.7cm]{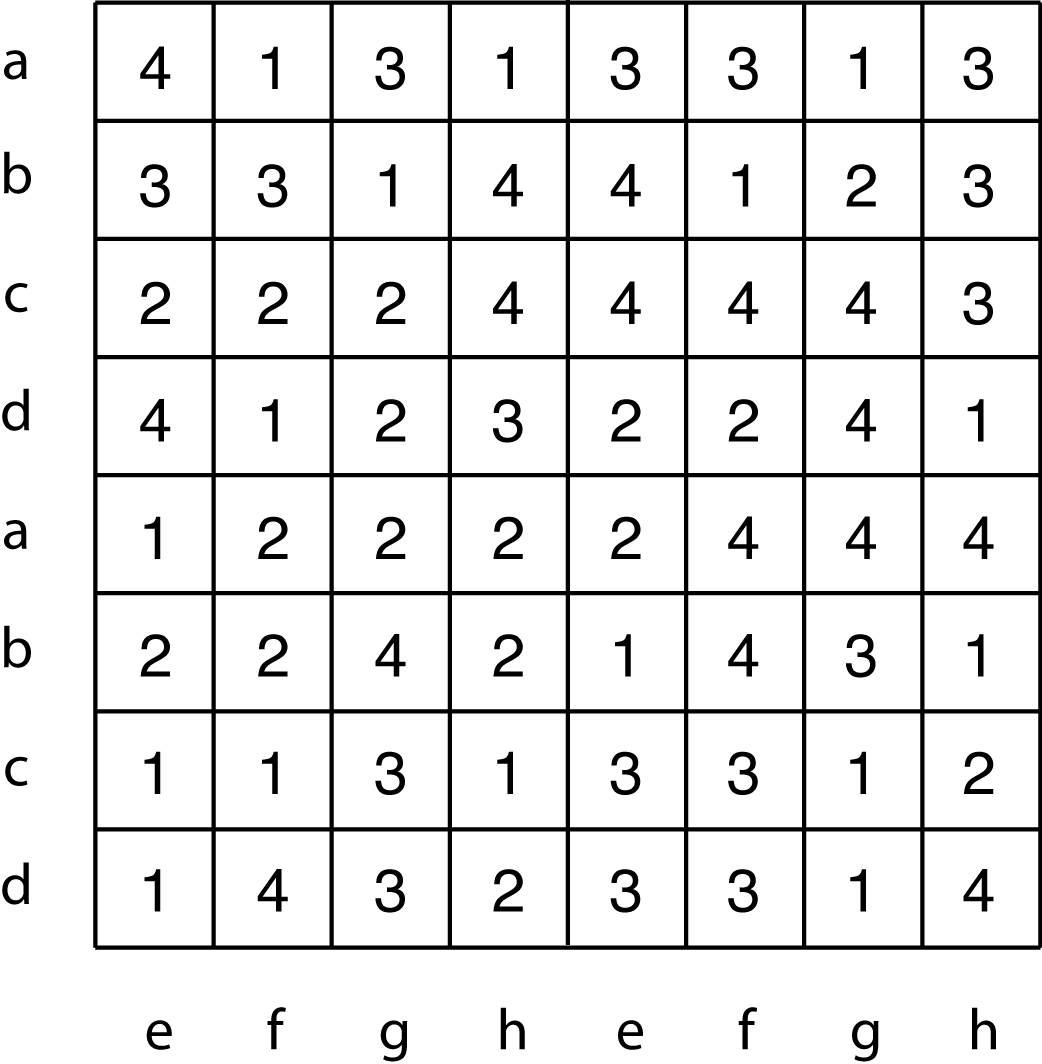}~~~~\includegraphics[height=3.7cm]{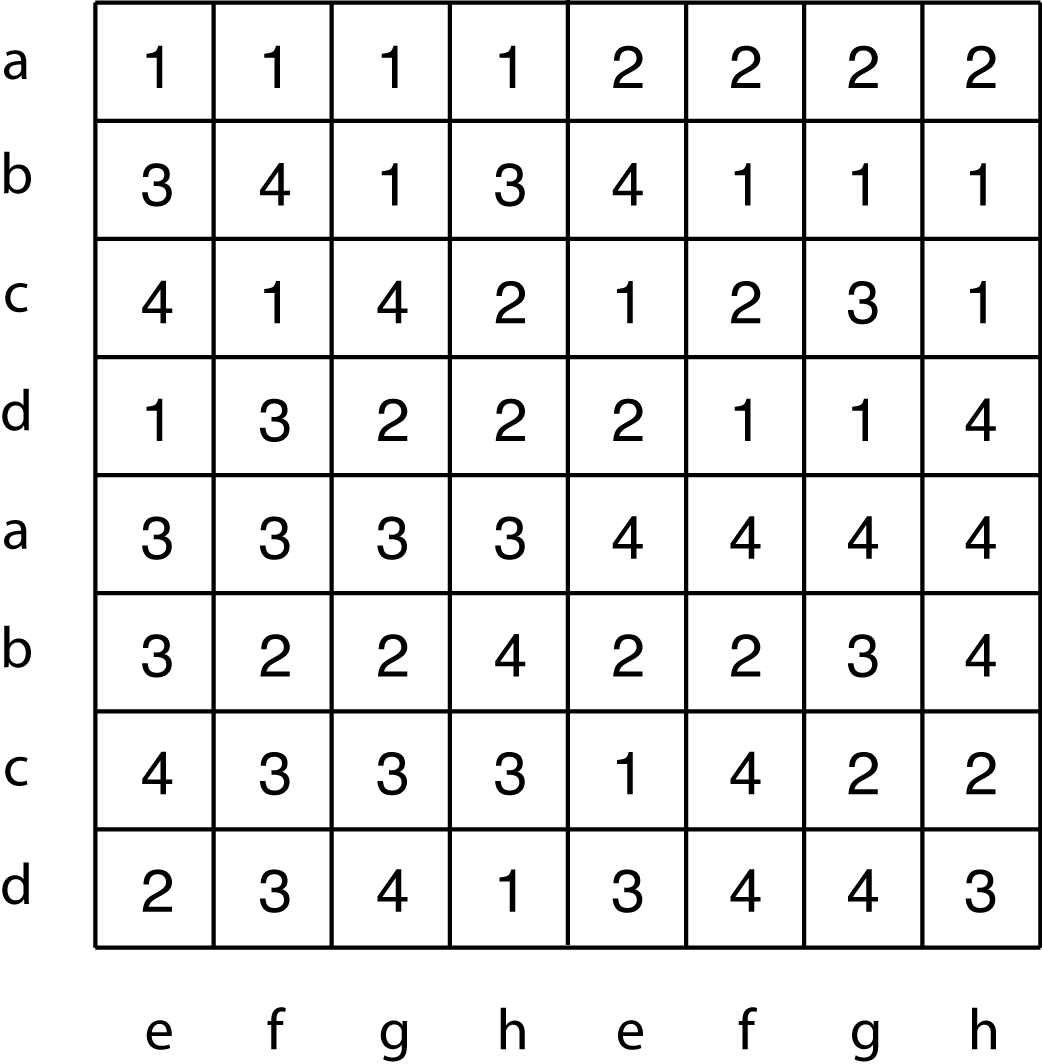}~~~~\includegraphics[height=3.7cm]{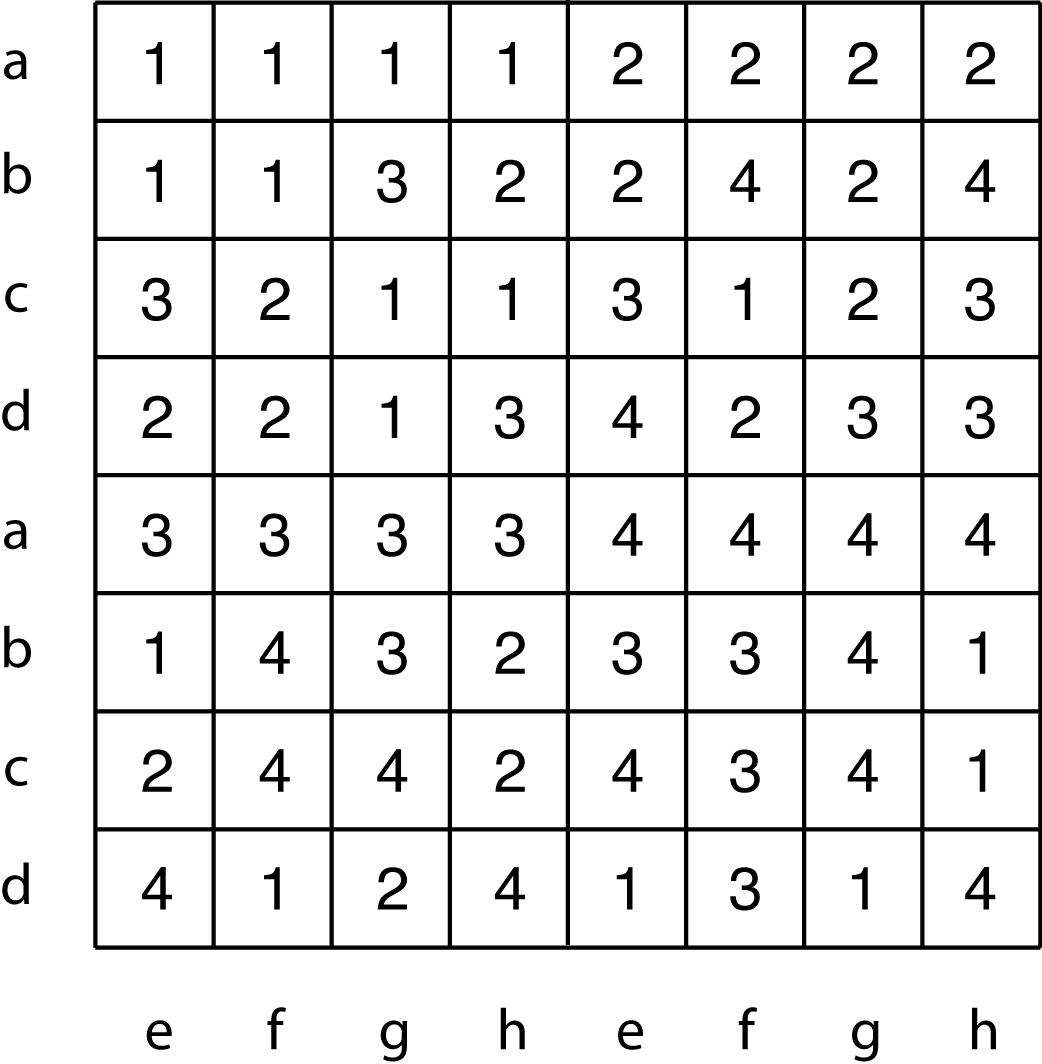}
\par\end{centering}

\protect\caption{\label{fig:Latin sequence of Free Latin squares}Uniform intersection
among asterisms}
\end{figure}

Using the multiset $\left\{ 1,1,1,1,2,2,2,2,3,3,3,3,4,4,4,4\right\} $
we have obtained the three Latin boards\footnote{As we will see in Sect. \ref{sub:Free-Latin-Squares}, these Latin
boards are called \emph{free Latin squares}.} shown in Fig. \ref{fig:Latin sequence of Free Latin squares}. The
derived woven boards have also SIN $\left\{ 0,4\right\} $. Let\textquoteright s
consider the following parallel classes, each with four asterisms:
\begin{itemize}
\item \textit{$H$}: $\left\{ a,b,c,d\right\} $
\item \textit{$V$}: $\left\{ e,f,g,h\right\} $
\item \textit{$S_{1}$}: containing one asterism per set of points with
the same symbol in the first Latin board in the figure
\item \textit{$S_{2}$}: \textit{idem} in the second one
\item \textit{$S_{3}$}: \textit{idem} in the third one
\end{itemize}
\pagebreak{}It is easily verified by inspection that any two asterisms,
each from a different class and not necessarily from the same board,
intersect at exactly four points: the three Latin boards constitute
then a Latin sequence that features uniform intersection among parallel
classes.

\section{Symmetric Boards\label{sec:Symmetric-Boards}}

In this section we study boards with geometric symmetry and how this
translates into combinatorial and algebraic symmetries.

\subsection{Symmetry Groups\label{sub:Symmetry-Groups}}
\begin{defn}
\label{def:Full-symmetry-groups}The (\textit{full}) \textit{symmetry
group} \textit{T }of a geometric object embedded in a space \textit{S}
is the group of all transformations\textit{ }under which the object
remains invariant, with composition of transformations as the group
operation.
\end{defn}
\noindent \textit{\emph{If $S$ is a metric space, $T$ is a subgroup
of the }}\textit{isometry group}\textit{\emph{ of $S$ \cite{Space Groups for Scientists and Engineers}.}}
\begin{example}
\label{exa:dihedral group D4}Consider the following eight  $\mathbb{R}{}^{2}$
transformations of the square:\end{example}
\begin{itemize}
\item 0º, 90º, 180º and 270º counter-clockwise rotations
\item 2 reflections across the two diagonals
\item 2 reflections, one across the vertical symmetry axis and another across
the horizontal one
\end{itemize}
Each of them leaves the square invariant. Together with the operation
of composition of transformations, they form the dihedral group \textit{$D_{4}$}
\cite{Regular polytopes} \textendash a subgroup of \emph{Euclidean
group} in two dimensions\textendash{} the algebraic structure that
captures the symmetry of the square. 
\begin{example}
The symmetry group of the equilateral triangle is the dihedral group
$D_{3}$, that comprises three counter-clockwise rotations ($0\text{º}$,
$120\text{º}$ and $240\text{º}$) and three reflections across the
axes that go through the vertices to the mid-point of the edges each.
\end{example}

\subsubsection{Symmetry Groups from Regular Polytopes\label{sub:Symmetry-Groups-from-finite-regular-polytopes}}

\noindent Natural candidates for \textit{T} in Def. \ref{def:Full-symmetry-groups}
are the symmetry groups of regular polytopes (see Sect. \ref{sub:Sources-of-Points}).
Examples in \textit{\emph{$\mathbb{R}^{2}$}} are the symmetry groups
of the regular polygons (the dihedral groups $D_{n}$). In \textit{\emph{$\mathbb{R}^{3}$}}
we have the symmetry groups of the finite real convex regular polyhedra
(see Table \ref{tab:symmetryGroupsPlatonicSolids-1}). For symmetry
groups in higher dimensions see \cite{Space Groups for Scientists and Engineers,Regular polytopes}.

\subsection{Actions of Symmetry Groups}

The effect of a symmetry group on a finite set of points is an instance
of the more general notion of \emph{action} of a generic group on
a generic finite set:
\begin{defn}
\label{def:group action}If \textit{T} is a group and \textit{P} is
a finite set, a \textit{group action on P} is a group homomorphism
\textendash called the \textit{action homomorphism}\textendash{} from
\textit{T} to \emph{Sym}(\textit{P}) that satisfies two axioms:\end{defn}
\begin{enumerate}
\item the identity element of \textit{T} is mapped to the identity permutation
of \textit{P} 
\item a composition of two elements of \textit{T} is mapped to the composition
of the corresponding permutations
\end{enumerate}
We say in this case that the group \emph{acts} on the set, and also
that the set is \emph{acted upon} by the group. There is a natural
action of symmetry groups on finite sets of geometric points, as the
next example shows.
\begin{example}
\label{exa:D3 action}Let $P$ be the set of vertices of an equilateral
triangle. If we map the identity transformation in $D_{3}$ to the
identity permutation in $Sym\left(P\right)$, and the other transformations
in $D_{3}$ to the permutation of the vertices that result after the
transformation is applied to the triangle, then it is easy to prove
that this map is a group homomorphism between $D_{3}$ and $Sym\left(P\right)$
that fulfills the axioms in Def. \ref{def:group action}, i.e.: the
map is a group action of $D_{3}$ on $P$. We say that $D_{3}$ acts
on the vertices of the triangle. By a similar argument we can prove
that $D_{3}$ also acts on the sides of the triangle.
\end{example}

\subsection{Definition of Symmetric Board\label{sub:Definition-of-Symmetric-Boards}}

\noindent A symmetry group that acts on the points of a board may
also act on its asterisms. We define:
\begin{defn}
A \textit{symmetric} \textit{board} is one whose asterisms are acted
upon by a symmetry group.
\end{defn}
\noindent In a symmetric board the elements of the group take asterisms
to asterisms, i.e.: it preserves the points, the constellation and
hence the board itself.
\begin{example}
\noindent The asterisms in the board in Fig. \ref{fig:3x3points and vert and hor asterisms}
are acted upon by \textit{$D_{4}$}, so the board is symmetric.
\end{example}

\subsubsection{Class-symmetric Boards\label{sub:Class-symmetric-Boards}}
\begin{defn}
The action of a group \textit{$T$} on a set \textit{$P$} is\textbf{
}\textit{transitive} if \textit{$P$} is non-empty, and if for any
\textit{$p_{1}$}, \textit{$p_{2}$} $\in$ \textit{$P$} there exists
a \textit{$t$} $\in$ \textit{$T$} such that\textit{\emph{ the action
of $t$ on}}\textit{ $p_{1}$} is $p_{2}$.
\end{defn}
\noindent So any element in \textit{$P$} can be reached from any
other if the action is transitive.
\begin{example}
\noindent \textit{$D_{4}$} does not act transitively on the asterisms
of the board in Fig. \ref{fig:3x3points and vert and hor asterisms},
as there is no element in \textit{$D_{4}$} that brings the top horizontal
asterism to the middle horizontal one. For the same reasons\textit{,
$D_{3}$} does not act transitively on the Fano Plane asterisms (see
Fig. \ref{fig:FanoPlaneUnsymmetricAndSymmetric} right).
\end{example}
\noindent If the symmetric board is resolvable (see Def. \ref{def:resolvable board})
then additional symmetry may be present at another level in the board:
\begin{defn}
Let\textit{ $B$} be a resolvable symmetric board with symmetry group
\textit{$T$}. If \textit{$T$} acts transitively on the parallel
classes, then \textit{$B$} is a \textit{class-symmetric }\textit{\emph{board}}.
\end{defn}
\noindent In a class-symmetric board the elements of the group not
only take classes to classes: there is always a sequence of group
elements that take any class in the board to any other too. A symmetric
board may not be class-symmetric for several reasons:
\begin{itemize}
\item \noindent there are no parallel classes (the case of the Fano Plane)
\item \noindent there are parallel classes but the group does not act upon
them
\item \noindent there are parallel classes acted upon by the group, but
the action is not transitive (the case of the board in Fig. \ref{fig:sudoku design pts-1})\end{itemize}
\begin{example}
\textit{$D_{4}$} acts on the set of parallel classes of the board
in Fig. \ref{fig:3x3points and vert and hor asterisms} (it takes
the set of rows to the set of columns and vice versa), so the board
is class-symmetric. This shows that a board may be class-symmetric
even if the group does not act transitively \textit{\emph{on its asterisms}}\emph{. }
\end{example}

\subsection{Permutations linked to Actions of Symmetry Groups}

The permutations linked to the action of a symmetry group form themselves
a group:
\begin{lem}
Let B = (P, C) be a symmetric board and T the symmetry group acting
on the asterisms in C. Then the permutations of P linked to the action
form a group with composition of permutations as the group operation.\end{lem}
\begin{proof}
Let \textit{h} be the action homomorphism between\textit{ T }and \emph{Sym}(\textit{P}),
with $P$ being a set of points with symmetry $T$. Since \textit{T}
acts on \textit{C} then, as each asterism is a subset of the points
and the union of all asterisms is $P$, \textit{T} also acts on \textit{P.}
Since a group homomorphism preserves the subgroups, and as \textit{T}
is an improper subgroup of itself, then \textit{h}(\textit{T}) is
a subgroup of \emph{Sym}(\textit{P}), hence a group.
\end{proof}

\subsection{Board Automorphisms from Geometric Symmetry\label{sec:B1 automorphisms from symmetry}}

Board automorphisms are point permutations that preserve the board
(see Sect. \ref{sub:Board-Automorphisms}). We see next that some
automorphisms come from the geometric symmetry of the board, i.e.:
how geometric symmetry entails algebraic geometry.
\begin{lem}
\label{lem:AT(B) included in Aut(B)}Let B = (P, C) be a symmetric
board and let T be the symmetry group acting on the asterisms in C.
Let A\textsubscript{\textit{T}}(B) be the group of permutations of
P linked to the action. Then A\textsubscript{\textit{T}}(B) ${\subseteq}$
Aut(B).\end{lem}
\begin{proof}
As \textit{T} acts on the asterisms in \textit{C,} then the elements
in their image \emph{A\textsubscript{T}}(\textit{B)} under the action
homomorphism are isomorphisms of \textit{B}, hence \emph{A\textsubscript{T}}(\textit{B})\textit{
}${\subseteq}$\textit{ }\emph{Aut}(\textit{B}).\end{proof}
\begin{example}
Consider the eight transformations in $D_{4}$ (see Example \ref{exa:dihedral group D4})
and the board in Fig. \ref{fig:3x3points and vert and hor asterisms}.
Each transformation is linked to a point permutation that leaves the
asterisms invariant (i.e.: to automorphisms). For example a 90º counter-clockwise
rotation induces the permutation (see Fig. \ref{fig:numberedPointsB1}):
\[
\left(\begin{array}{ccccccccc}
1 & 2 & 3 & 4 & 5 & 6 & 7 & 8 & 9\\
2 & 4 & 1 & 8 & 5 & 2 & 9 & 6 & 3
\end{array}\right)
\]

\end{example}
\noindent which is an automorphism. As is the permutation linked to
a reflection across the vertical axis:

\[
\left(\begin{array}{ccccccccc}
1 & 2 & 3 & 4 & 5 & 6 & 7 & 8 & 9\\
3 & 2 & 1 & 6 & 5 & 4 & 9 & 8 & 7
\end{array}\right)
\]
\smallskip{}

\noindent Other automorphisms of the board have no correspondence
in \textit{$D_{4}$. }\textit{\emph{A}}n example\textit{\emph{ is
the a}}utomorphism linked to permutation \ref{eq:a column permutation}:
it swaps the two leftmost columns in the resulting Latin squares,
something that no element of $D_{4}$ can perform.

\subsection{Group Hierarchy\label{sub:Group-Hierarchy}}

The relation among symmetry groups, automorphism groups and symmetric
groups (see Sect. \ref{sub:Board-Automorphisms}) is given by Lemma
\ref{lem:AT(B) included in Aut(B)} and Relation \eqref{eq:Aut(B) included in Sym(P)}
in Sect. \ref{sub:Board-Automorphisms}\label{Relation three groups }:
\[
A_{T}(B)\subseteq Aut(B)\subseteq Sym(P)
\]

\noindent This shows that if we want interesting boards with large
\emph{Aut}(\textit{B}), a symmetric board whose asterisms are acted
upon by a symmetry group \textit{T} provides at least \emph{A\textsubscript{T}}(\textit{B})
upfront. 
\begin{example}
Let \textit{B} = (\textit{P}, \textit{C}) be a board with \textit{$P=\left\{ v_{i},i=1,2,3\right\} $}
being the set of vertices of an equilateral triangle, with \textit{$C=\left\{ \left\{ v_{i}\right\} ,i=1,2,3\right\} $}
(i.e.: each asterism has just a single point). As the dihedral group
\textit{$D_{3}$} acts on the asterisms in \textit{$C$,} \textit{$B$}
is a symmetric board. We have in this case $A_{D_{3}}(B)=Aut(B)=Sym(P)$,
all with order $6.$
\end{example}

\begin{example}
Let \textit{$B=(P,C)$} be a board with \textit{$P=\left\{ v_{i},i=1,2,3,4\right\} $}
being the set of vertices of the square, with \textit{$C=\left\{ \left\{ v_{i}\right\} ,i=1,2,3,4\right\} $}.
As \textit{$D_{4}$} acts on the asterisms in \textit{$C$,} \textit{$B$}
is a symmetric board. We now have \emph{$A_{D_{4}}(B)=Aut(B)$} and
\emph{$Aut(B)\subset Sym(P)$}. \emph{$Sym(P)$} has order $4!=24$,
while $A_{D_{4}}(B)$ and \emph{$Aut(B)$} have both order $8$.
\end{example}

\subsection{Symmetric Boards from Sources of Symmetry\label{sub:Sources of Symmetry}}

A way to find symmetric boards\textit{ }is to start with a \emph{source
of symmetry} (or \textit{source} for short): a geometric object with
a symmetry group \textit{$T$} (see Sect. \ref{sub:Symmetry-Groups-from-finite-regular-polytopes}).
Board points can then be picked in the source and asterisms formed
in such a way that \textit{T} acts on them. If we want the board to
be class-symmetric, classes can be chosen so that $T$ acts transitively
on them. The resulting board \textit{$B$} is assured to have at least
\emph{$A_{T}(B)$} as an automorphism group following Lemma \eqref{lem:AT(B) included in Aut(B)}.

Sources are usually drawn along with the boards to render the symmetry
apparent. If the source is intricate or has many dimensions it may
be given analytically instead of graphically.

We describe next how to obtain different sources, and give example
boards based on them. To highlight the effect of geometric symmetry
on combinatorial symmetry (see Sect. \ref{sub:Combinatorial-Properties}),
the examples include properties such as the pattern of intersection
of the asterisms (SIN), existence of parallel classes, orthogonality,
resolvability, etc.

\subsubsection{Sources from Regular Polytopes \label{sub:Sources-of-Points}}

\noindent A \emph{polytope} is a geometrical figure bounded by portions
of lines, planes and hyperplanes in a particular space \cite{Regular polytopes}.
\emph{Polygons} and \emph{polyhedra} are polytopes in $\mathbb{R}^{2}$
and $\mathbb{R}^{3}$ respectively. \emph{Regular polytopes} are polytopes
with special properties of symmetry that make them ideal sources.
Among the finite real regular polytopes we have\footnote{For a formal definition and a comprehensive list of regular real polytopes
\textendash including \emph{apeirotopes}, the infinite ones\textendash{}
see \cite{Regular polytopes}. For the \emph{complex regular polytopes}
see \cite{Regular Complex Polytopes}; for \emph{abstract polytopes}
-a purely combinatorial generalization- see \cite{Abstract Regular Polytopes}.}:
\begin{itemize}
\item \noindent in $\mathbb{R}$ the closed intervals
\item \noindent in $\mathbb{R}^{2}$ the regular polygons (convex) and the
regular star-polygons (non-convex)
\item \noindent in $\mathbb{R}^{3}$ the five platonic polyhedra (convex)
and the four Kepler-Poinsot polyhedra (non-convex)
\item \noindent in $\mathbb{R}^{4}$ the six convex regular $4$-polytopes
and the ten non-convex Schläfli-Hess 4-polytopes
\item \noindent in $\mathbb{R}^{n}$ (\textit{$n>4$}) the simplexes, cubes
and cross-polytopes\end{itemize}
\begin{example}
\label{exa:helios design points}The source in Fig. \ref{fig:emptyBrdHelios}
(left) has the symmetry of the regular dodecahedron (dihedral symmetry
\textit{$D_{12}$}). Let\textquoteright s take as board points the
intersections between segments (thirty six intersections in all) and
let\textquoteright s define a $6$-uniform asterism out of each $6$-set
of points lying on the same segment (twelve asterisms in all). The
board is shown on the right. The asterisms are acted upon by \textit{$D_{12}$}
and hence we have a symmetric board. The board is not class-symmetric
because there are no parallel classes here.
\begin{figure}
\begin{centering}
\includegraphics[height=4.5cm]{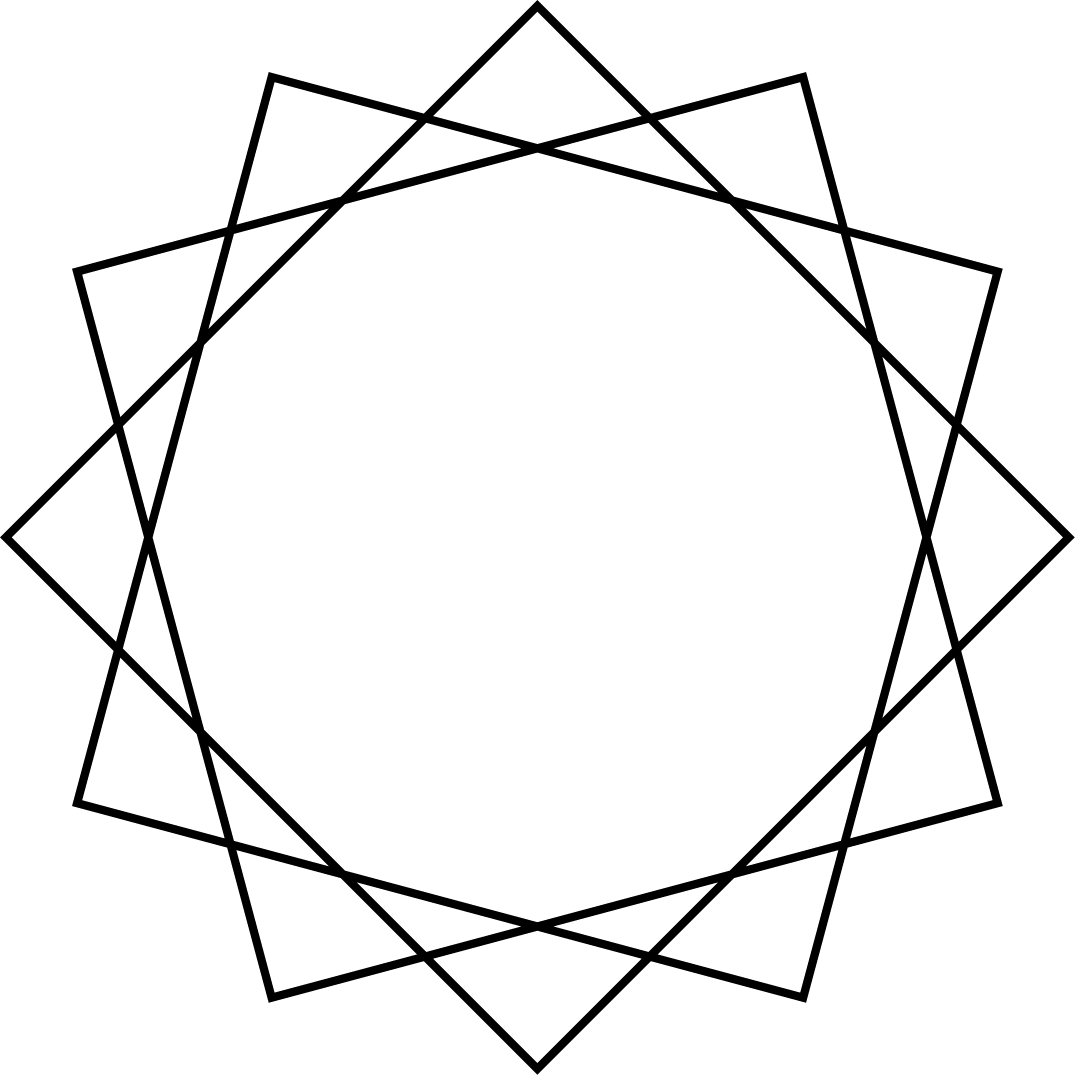}~~~~~\includegraphics[height=4.5cm]{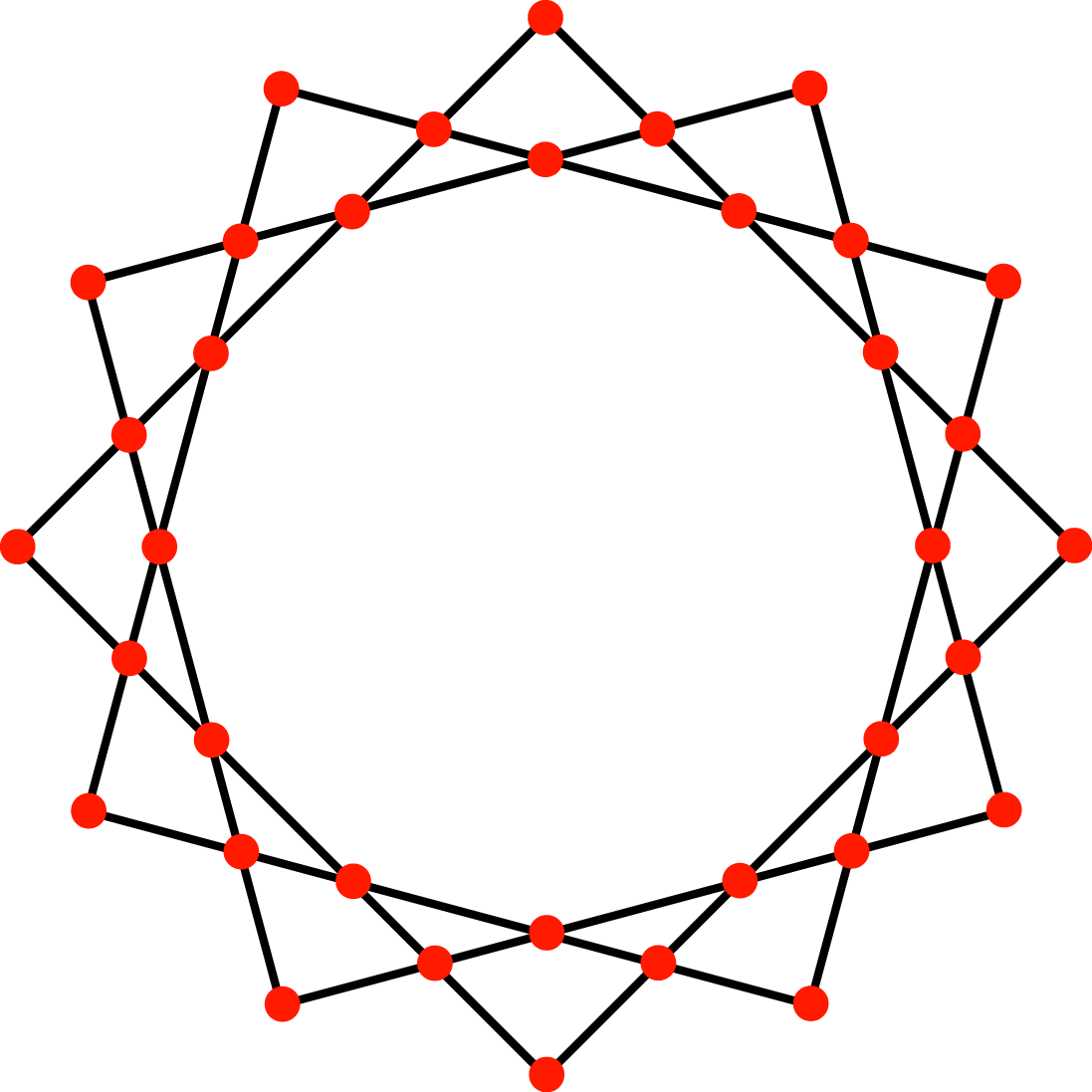}
\par\end{centering}

\protect\caption{\label{fig:emptyBrdHelios}Source with $D_{12}$ symmetry and symmetric
board}
\end{figure}

\end{example}

\subsubsection{Sources from Projections}

A way to obtain sources in a space with a given dimension is to properly
project onto it sources from higher dimensions. We can obtain bidimensional
sources with dihedral symmetry as proper projections of the regular
polyhedra \cite{Regular polytopes}.

\subsubsection{Sources from Kaleidoscopic Sets\label{sub:Sources-from-Actions-Symmetry-groups}}

Sources can also be built from scratch in any geometric space using
symmetry groups (see Sect. \ref{sub:Symmetry-Groups-from-finite-regular-polytopes}),
group actions (see Def. \ref{def:group action}) and ideas from \cite{Regular polytopes}.
\begin{defn}
Let \textit{$T$} be a group acting on a set of points \textit{$P$}.
The \textit{orbit} of a point \textit{$p$} is the set of points in
\textit{$P$} to which \textit{$p$} is moved by the elements of \textit{$T$}.
\end{defn}

\begin{defn}
Let \textit{$T$} be a group acting on a space \textit{$S$}. A \textit{fundamental
region} is a subset of \textit{$S$} containing exactly one point
from the orbit of every point in \textit{$S$}.
\end{defn}

\begin{defn}
Let \textit{$T$} be a symmetry group acting on the points of a space
\textit{$S$}. Let \textit{$F$} be a set of points in a fundamental
region of \textit{$S$}. Then \textit{$K=\left\{ orbit\,of\,p,p\in F\right\} $}
is the \textit{kaleidoscopic set }\textit{\emph{of}}\textit{ $F$
}\textit{\emph{by}}\textit{ $T.$}
\end{defn}
\noindent A kaleidoscopic set so constructed is evidently a source.
The elements in the definition have clear counterparts in a kaleidoscope:
points in \textit{$F$} are the colored beads inside the tube, while
the fundamental region is the chamber in which they are confined;
the set of images of a bead is the orbit of a point; \textit{$S$}
is the plane in which the beads move; the kaleidoscopic set is the
image seen for a certain rotation of the tube, while \textit{$T$}
generators are the inside mirrors.
\begin{example}
Let \textit{$S$} be \textbf{$\mathbb{R}^{2}$} and \textit{$T$}
be \textit{$D_{3}$}. Fig. \ref{fig:smallerTriangKaleidoSetSmall}
(left) shows the fundamental regions, six in this case. \textit{$F$}
is in the middle and \textit{$K$} on the right. Now let\textquoteright s
choose as board points its fifteen vertices and make an asterism out
of each pair of vertices on the same segment (thirty asterisms in
all). The board is shown on Fig. \ref{fig:Board from kaleidoscopic set}.
As \textit{$D_{3}$} acts on the asterisms we have a symmetric board.
\end{example}

\begin{example}
\noindent Let \textit{$S$} be \textbf{$\mathbb{R}^{2}$} and \textit{$T$}
be\textit{ $D_{4}$}. The eight fundamental regions are shown in Fig.
\ref{fig:biregSqreAsKaleido} (left), \textit{$F$} is in the center
and \textit{$K$ }on the right. If we choose the centers of the squares
as board points and make asterisms out of every row and every column
of points we obtain the board Latin squares are based on.
\begin{figure}
\begin{centering}
\includegraphics[height=3.1cm]{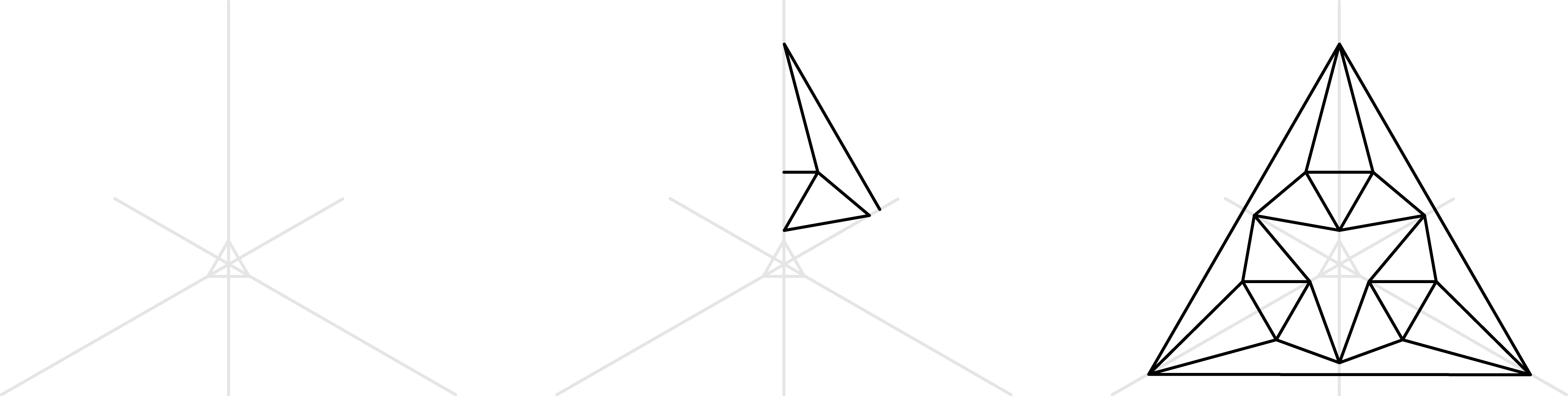}
\par\end{centering}

\protect\caption{\label{fig:smallerTriangKaleidoSetSmall}Kaleidoscopic set by $D_{3}$}
\end{figure}
\begin{figure}
\begin{centering}
\includegraphics[height=2.7cm]{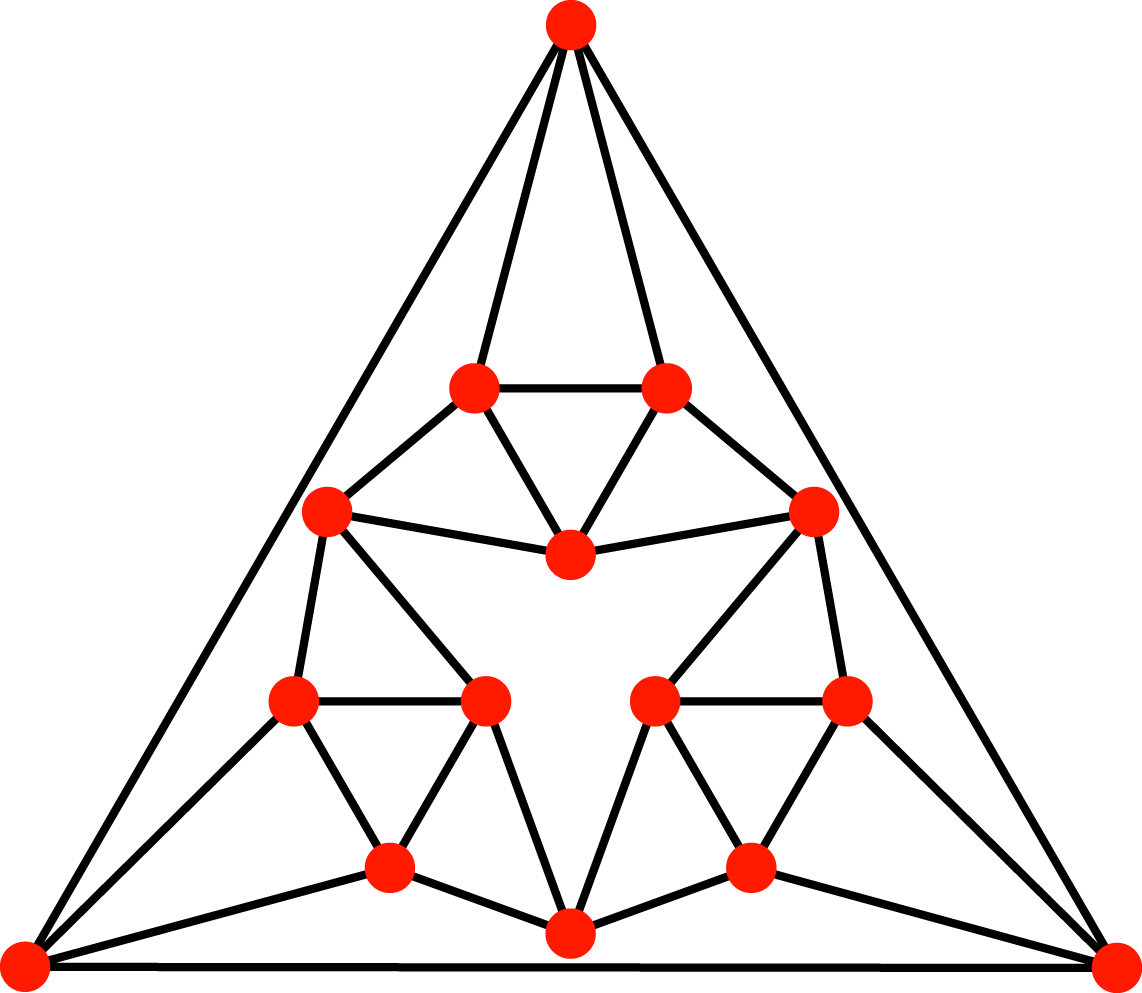}
\par\end{centering}

\protect\caption{\label{fig:Board from kaleidoscopic set}$D_{3}$ symmetric board}
\end{figure}
\begin{figure}
\begin{centering}
\includegraphics[height=3cm]{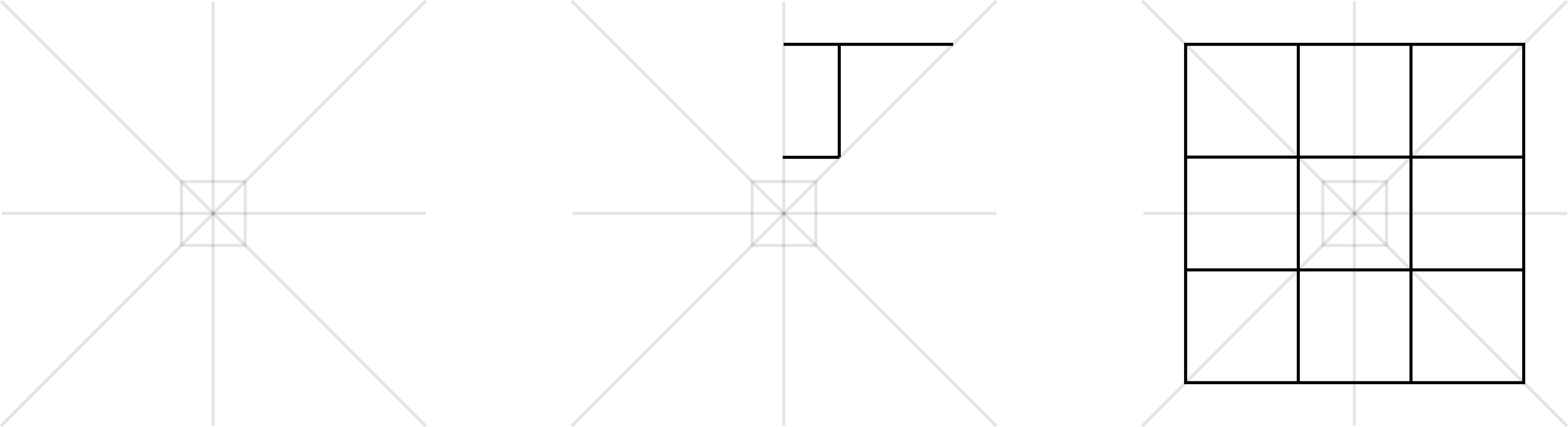}
\par\end{centering}

\protect\caption{\label{fig:biregSqreAsKaleido}Kaleidoscopic set by $D_{4}$}
\end{figure}

\end{example}

\subsubsection{Enriching Sources}

Sources with rich structure give more options to define board points
and asterisms. We can obtain a richer source from a given one with
a combination of techniques that add structure without breaking the
symmetry, like the ones described next.

\paragraph{Composition}

We can combine several sources with the same symmetry to obtain a
richer one. \emph{Compound polytopes} \cite{Regular polytopes} \textendash sets
of regular polytopes with a common center\textendash{} are an example.

\paragraph{Truncation}

A source properly truncated may gain new vertices, edges, faces, etc.
depending on the dimension of the source. See \cite{Regular polytopes}
for truncated versions of the regular polytopes.

\paragraph{Kaleidoscopic Sets}

If we have a source with symmetry group $T$ we can enrich it while
preserving its symmetry by adding a kaleidoscopic set generated by
$T$ (see Sect. \ref{sub:Sources-from-Actions-Symmetry-groups}).
We can add in this way more than one point per edge, face or vertex.

\paragraph{Extra Regular Tilings. Biregular Polytopes\label{sub:Biregular-polytopes}}

We define\emph{ biregular polytopes} as regular polytopes tiled with
other regular polytopes. In two dimensions for example we have \emph{biregular
polygons}: regular polygons tiled with regular polygons. In the context
of sources biregular polytopes thus enrich regular polytopes.

The number of tiling polygon sides on each side of the tiled polygon
is the \textit{order} of the biregular polygon. The tiling polygons
are called \textit{faces} or \emph{cells}. Since there are only three
regular tessellations in \textbf{$\mathbb{R}^{2}$} \cite{Regular polytopes}
it is easy to prove that in this space there are only three families
of biregular polygons (see Fig. \ref{fig:Examples of biregular triangle, square and hexagon}):
\begin{itemize}
\item \noindent equilateral triangles tiled with equilateral triangles (\textit{biregular
triangles})
\item \noindent squares tiled with squares (\textit{biregular squares})
\item \noindent regular hexagons tiled with equilateral triangles (\textit{biregular}
\textit{hexagons})
\end{itemize}
\noindent As biregular polygons feature dihedral symmetry, they are
eligible sources for both symmetric and class-symmetric boards (see
Sect. \ref{sub:Definition-of-Symmetric-Boards}). They also feature
structure and scalability while keeping the symmetry, so entire families
of symmetric boards can be defined in a very compact way. Table \ref{tab:biregular polys params}
shows relevant parameters of each family. We can enrich candidate
sources by transforming parts of them that are regular polytopes:
\begin{itemize}
\item polygons into biregular polygons (see Sects. \ref{exa:Boards-from-Birregular-tri},
\ref{exa:Boards-from-Birregular-squares}, \ref{exa:Boards-from-Birregular-hex}
and \ref{exa:Boards-from-Polyhedra with bireg faces})
\item polyhedra into biregular polytopes
\item \textit{$n$}-polytopes into biregular polytopes
\begin{figure}
\begin{centering}
\includegraphics[height=2.4cm]{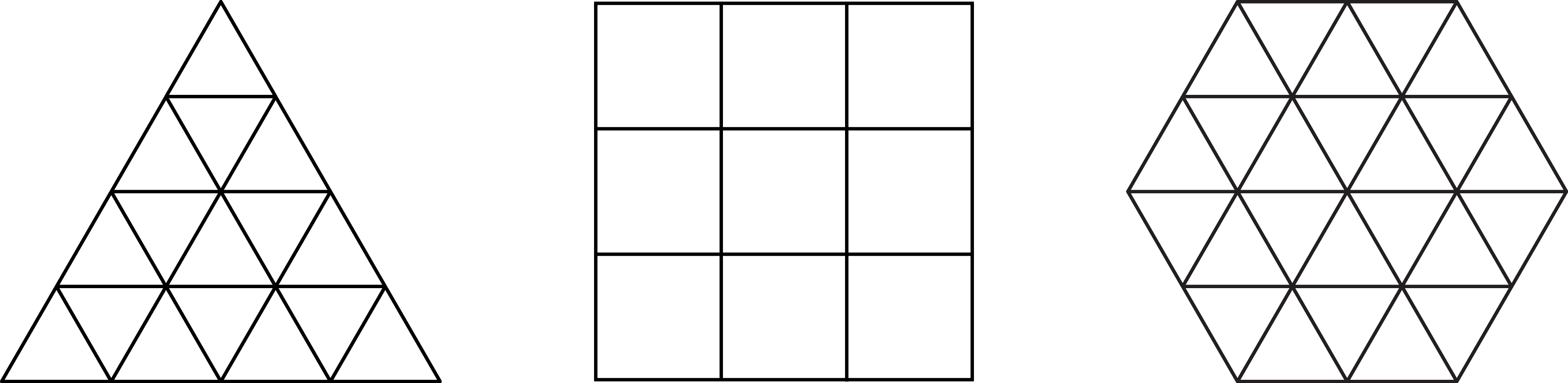}
\par\end{centering}

\protect\caption{\label{fig:Examples of biregular triangle, square and hexagon}Biregular
polygons}
\end{figure}
\begin{table}
\begin{centering}
\begin{tabular}{>{\centering}p{1.1cm}>{\centering}p{0.6cm}>{\centering}p{1.1cm}>{\centering}p{1.1cm}>{\centering}p{1.1cm}>{\centering}p{1.6cm}c>{\centering}p{0.8cm}}
\textbf{\scriptsize{}biregular polygon} & \textbf{\scriptsize{}order} & \textbf{\scriptsize{}tiling polygon} & \textbf{\scriptsize{}tiled polygon} & \textbf{\scriptsize{}symmetry group} & \textbf{\scriptsize{}\#vertices} & \textbf{\scriptsize{}\#edges} & \textbf{\scriptsize{}\#faces}\tabularnewline
\hline 
\noalign{\vskip0.1cm}
\textbf{\tiny{}triangle} & \emph{\footnotesize{}n} & {\scriptsize{}triangle} & {\scriptsize{}triangle} & \textit{\footnotesize{}D}{\footnotesize{}\textsubscript{4}} & {\footnotesize{}$\frac{(n+1)\text{·}(n+2)}{2}$} & {\footnotesize{}$\frac{3\text{·}n\text{·}(n+1)}{2}$} & {\footnotesize{}$n^{2}$}\tabularnewline
\noalign{\vskip0.1cm}
\textbf{\tiny{}square} & \emph{\footnotesize{}n} & {\scriptsize{}square} & {\scriptsize{}square} & \textit{\footnotesize{}D}{\footnotesize{}\textsubscript{4}} & $(n+1)^{2}$ & $2n(n+1)$ & $n^{2}$\tabularnewline
\noalign{\vskip0.1cm}
\textbf{\tiny{}hexagon} & \emph{\footnotesize{}n} & {\scriptsize{}triangle} & {\scriptsize{}hexagon} & \textit{\footnotesize{}D}{\footnotesize{}\textsubscript{6}} & {\footnotesize{}$3n^{2}+3n+1$} & $3n(3n+1)$ & $6n^{2}$\tabularnewline
\end{tabular}
\par\end{centering}

\protect\caption{\label{tab:biregular polys params}Parameters of biregular polygons}

\end{table}

\end{itemize}

\subparagraph{Biregular Triangles\label{exa:Boards-from-Birregular-tri}}
\begin{example}
\label{exa:latin triangle faces}Fig. \ref{fig:raw la tri order 6 faces}
shows a biregular triangle of order six where face barycenters (red
dots) have been chosen as board points. An asterism here is the set
of points pointed to by the same letter. Each asterism has twelve
points and belongs to one of three parallel classes: $\left\{ a,b,c\right\} $,
$\left\{ d,e,f\right\} $ or $\left\{ g,h,i\right\} $. We have then
a uniform resolution. Furthermore, any two asterisms from different
classes intersect in four points. To see this, take the smaller, three-row
upper triangle, and turn it clockwise around the center of the larger
triangle's right side until both halves of this side coincide. We
see that every asterism in the first class intersects in four points
any other  in the second. The symmetry of the board allows us to conclude
that this applies to any two asterisms from different classes, so
the SIN of the board is $\left\{ 0,4\right\} $. It is simple to verify
that \textit{$D_{3}$} acts transitively on the classes, so we have
a class-symmetric board (see Sect. \ref{sub:Class-symmetric-Boards}).
In fact, provided that the order \textit{n} of the biregular triangle
is even, we have a whole family of class-symmetric boards likewise
built, each with \textit{n}\textsuperscript{2} points and three uniform
parallel classes, each with \textit{$\frac{n}{2}$} asterisms of size
$2n$\textit{ }each.
\begin{figure}
\begin{centering}
\includegraphics[height=5cm]{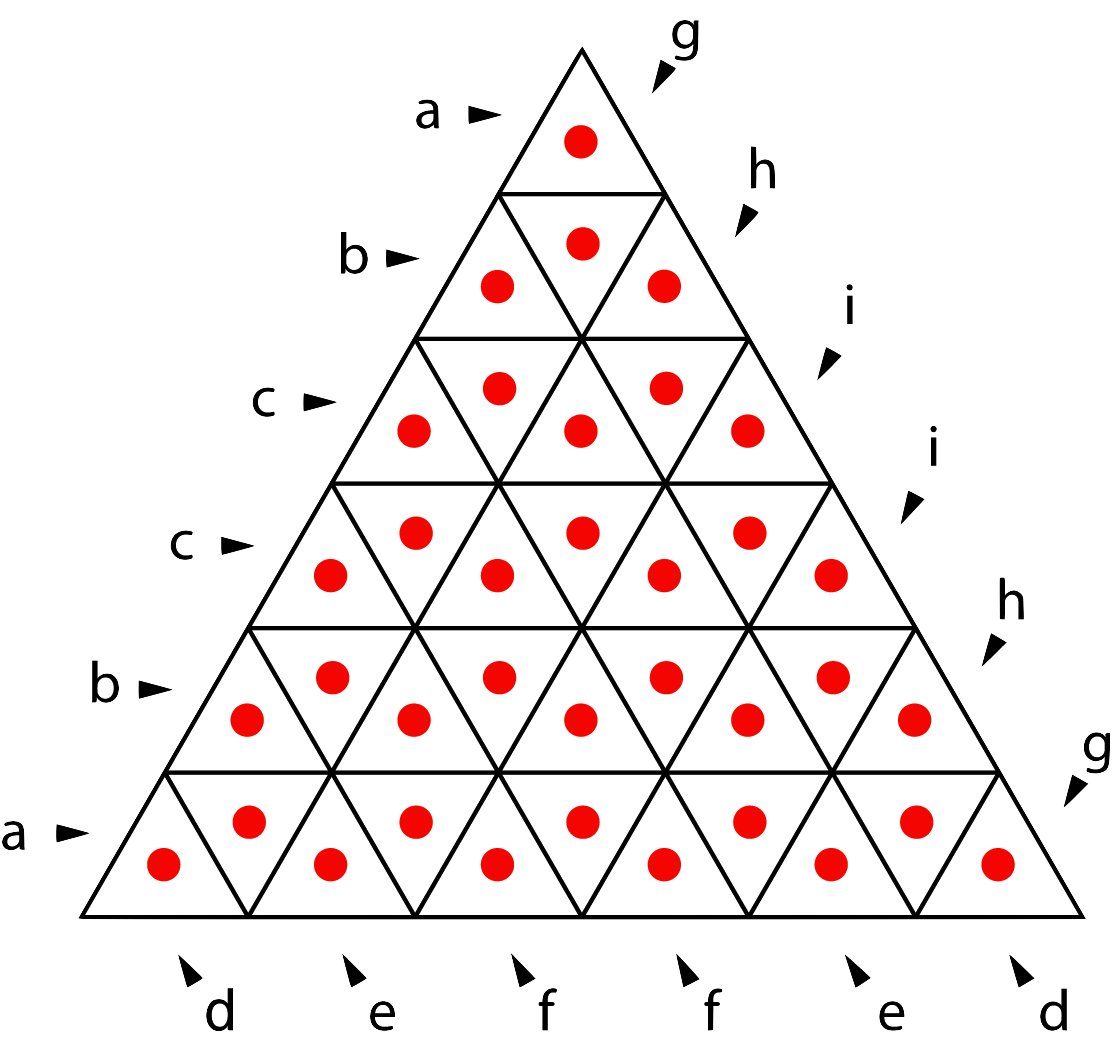}
\par\end{centering}

\protect\caption{\label{fig:raw la tri order 6 faces}Triangular board with points
at the barycenters}
\end{figure}

\end{example}

\begin{example}
Fig. \ref{fig:raw la tri order 7 vtcesSmall} shows a biregular triangle
of order seven with vertices as board points. An asterism here is
the set of points pointed to by the same letter. Each has nine points
and belongs to one of three parallel classes: $\left\{ a,b,c,d\right\} $,
$\left\{ e,f,g,h\right\} $ or $\left\{ i,j,k,l\right\} $. We have
then a board with uniform resolution. It is simple to verify that
\textit{$D_{3}$} acts transitively on the classes. In fact, provided
that the order \textit{$n$} of the biregular triangle is odd, we
have a whole family of class-symmetric boards, each with $\frac{(n+1)(n+2)}{2}$
points and three uniform parallel classes, each with $\frac{(n+1)}{2}$
asterisms of size \textit{$n+2$} each.
\begin{figure}
\begin{centering}
\includegraphics[height=5cm]{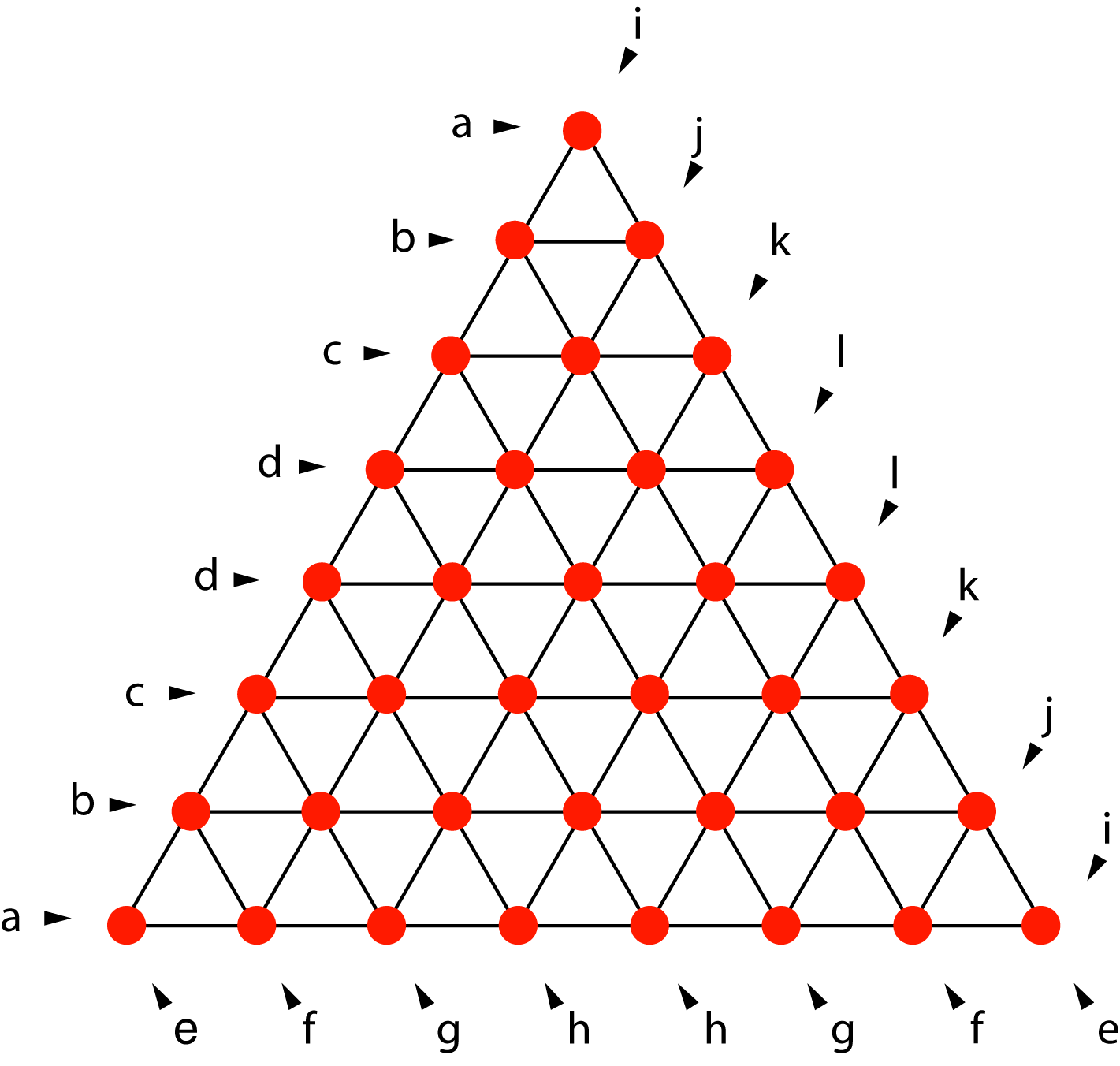}
\par\end{centering}

\protect\caption{\label{fig:raw la tri order 7 vtcesSmall}Triangular board with points
at the vertices}
\end{figure}

\end{example}

\begin{example}
Fig. \ref{fig:raw lat tri order 4 edges 3x10P} shows a biregular
triangle of order five on which we have chosen edge centers as board
points. As before, an asterism is specified as the set of points pointed
to by the same letter. Each asterism has ten points and belongs to
one of three parallel classes: $\left\{ a,b,c\right\} $, $\left\{ d,e,f\right\} $
or $\left\{ g,h,i\right\} $. We have again a uniform resolution.
It is simple to verify that \textit{$D_{3}$} acts transitively on
the classes, so this is a class-symmetric board.
\begin{figure}
\begin{centering}
\includegraphics[height=5.8cm]{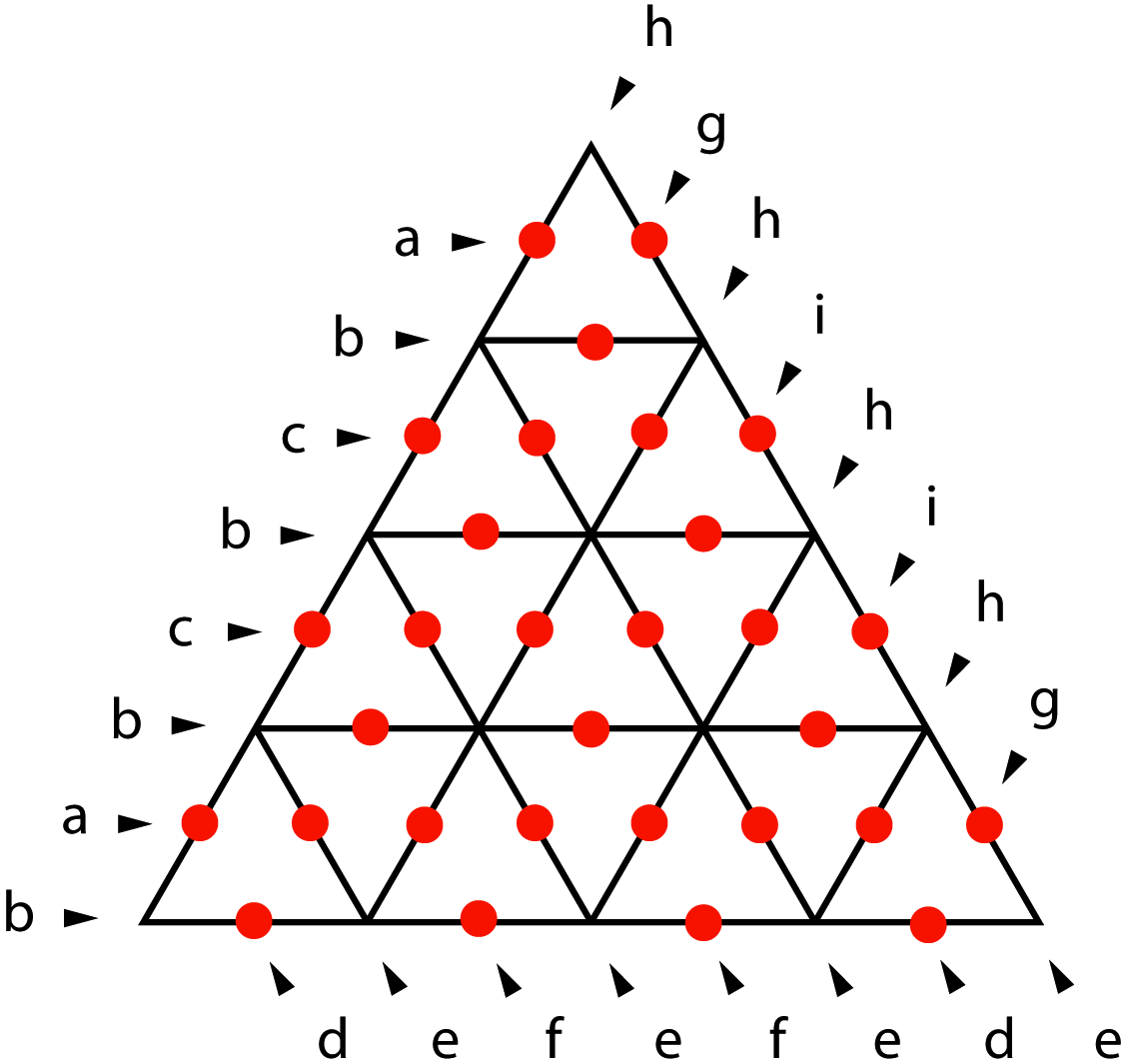}
\par\end{centering}

\protect\caption{\label{fig:raw lat tri order 4 edges 3x10P}Triangular board with
points at the edge midpoints}
\end{figure}

\end{example}

\subparagraph{Biregular Squares\label{exa:Boards-from-Birregular-squares}}
\begin{example}
\label{exa:free LS board}Fig. \ref{fig:board Free Latin Square}
shows a biregular square of order eight (see Sect. \ref{sub:Biregular-polytopes})
with face centers as board points. A board asterism is made up of
sets of points pointed to by the same letter. Each asterism has sixteen
points and belongs to parallel class $\left\{ a,b,c,d\right\} $ or
$\left\{ e,f,g,h\right\} $. We have then a $16$-uniform resolution
in which every pair of asterisms not in the same class share four
points, so the board has SIN is $\left\{ 0,4\right\} $. It is straightforward
to verify that \textit{$D_{4}$} acts transitively on the classes,
so the board is class-symmetric. We can build symmetric boards on
biregular squares of any order \textit{$n$} \textit{\emph{with}}\textit{
$n^{2}$} points. If the order is prime we can only have symmetric
boards that will originate conventional Latin squares. For composite
orders we have several options to group rows \textendash and columns\textendash{}
into board asterisms. If \textit{$n$} is even and asterisms come
in two pieces as in Example \ref{exa:free LS board}, the board will
have two uniform parallel classes, each with \textit{$\frac{n}{2}$}
asterisms of size $2n$ each. The pairing chosen will determine if
the board is class-symmetric, just symmetric or non-symmetric.
\end{example}

\begin{example}
\label{exa:Tartan board-1}Fig. \ref{fig:tartan board-1} shows the
board in Fig. \ref{fig:board Free Latin Square} with an extra parallel
class whose asterisms are the four $4\times4$ highlighted sub-squares.
Its SIN is also $\left\{ 0,4\right\} $. The board is evidently symmetric
but not class-symmetric, as there is no way to take the extra class
to one of the others with the elements of $D_{4}$.
\begin{figure}
\centering{}\includegraphics[height=4.5cm]{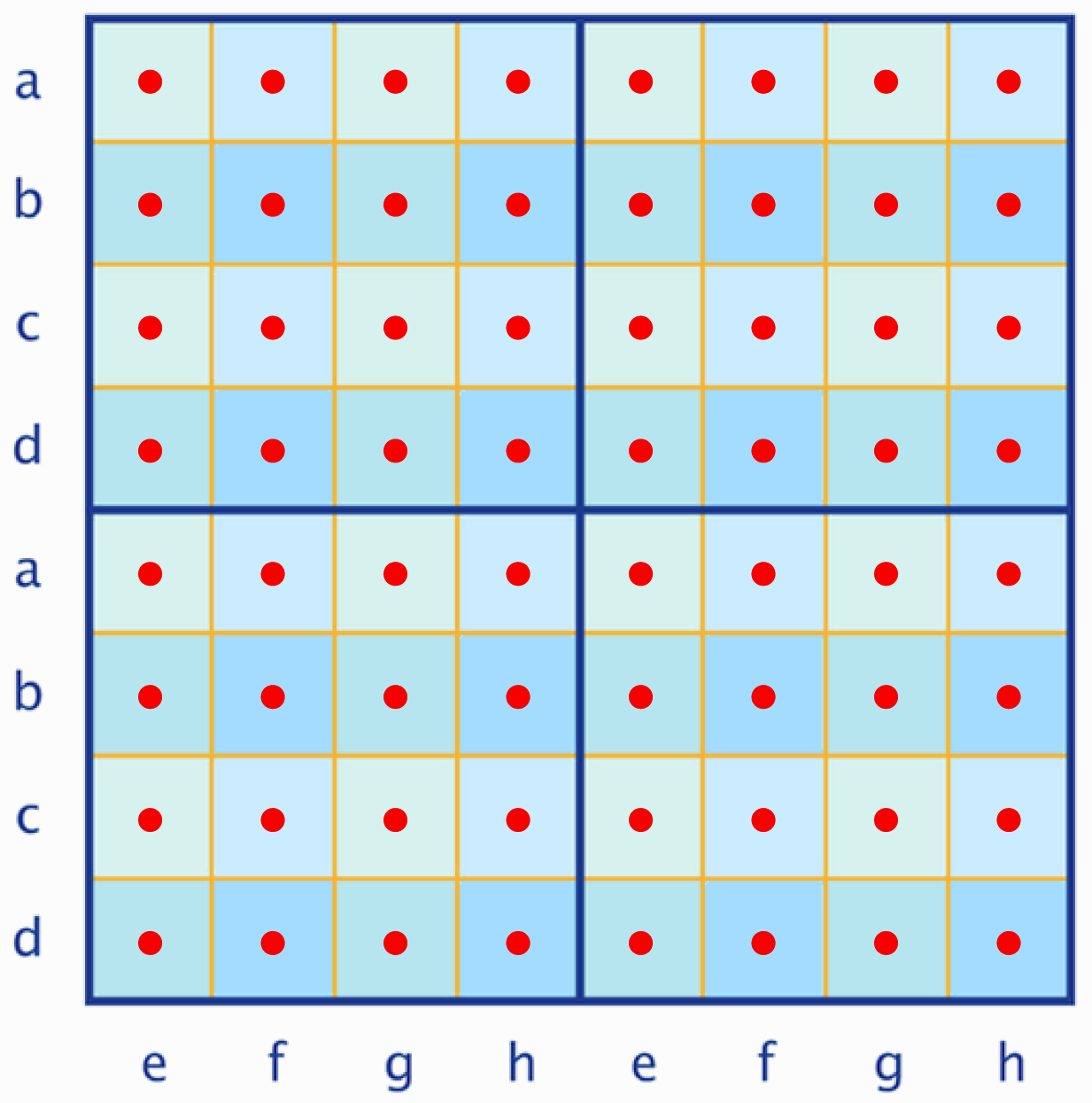}\protect\caption{\label{fig:tartan board-1}Square board with SIN $\left\{ 0,4\right\} $}
\end{figure}

\end{example}

\begin{example}
Fig. \ref{fig:punctSquareBoard} shows an arrangement of points with
\textit{$D_{4}$} symmetry. They come from the vertices and tile centers
\textendash with one missing\textendash{} of a biregular square of
order four. They have been organized into asterisms, each one formed
by points pointed to by the same letter. The asterisms have eight
points each and are acted upon by \textit{D}\textsubscript{4 }\textendash asterisms
$b$ and $h$ for example are swapped by a vertical reflection. This
is thus a 8-uniform \textit{$D_{4}$} board that is also class-symmetric.
Its SIN is $\left\{ 0,2\right\} $, as any two asterisms either intersect
at two points or do not intersect at all.
\begin{figure}
\begin{centering}
\includegraphics[height=5cm]{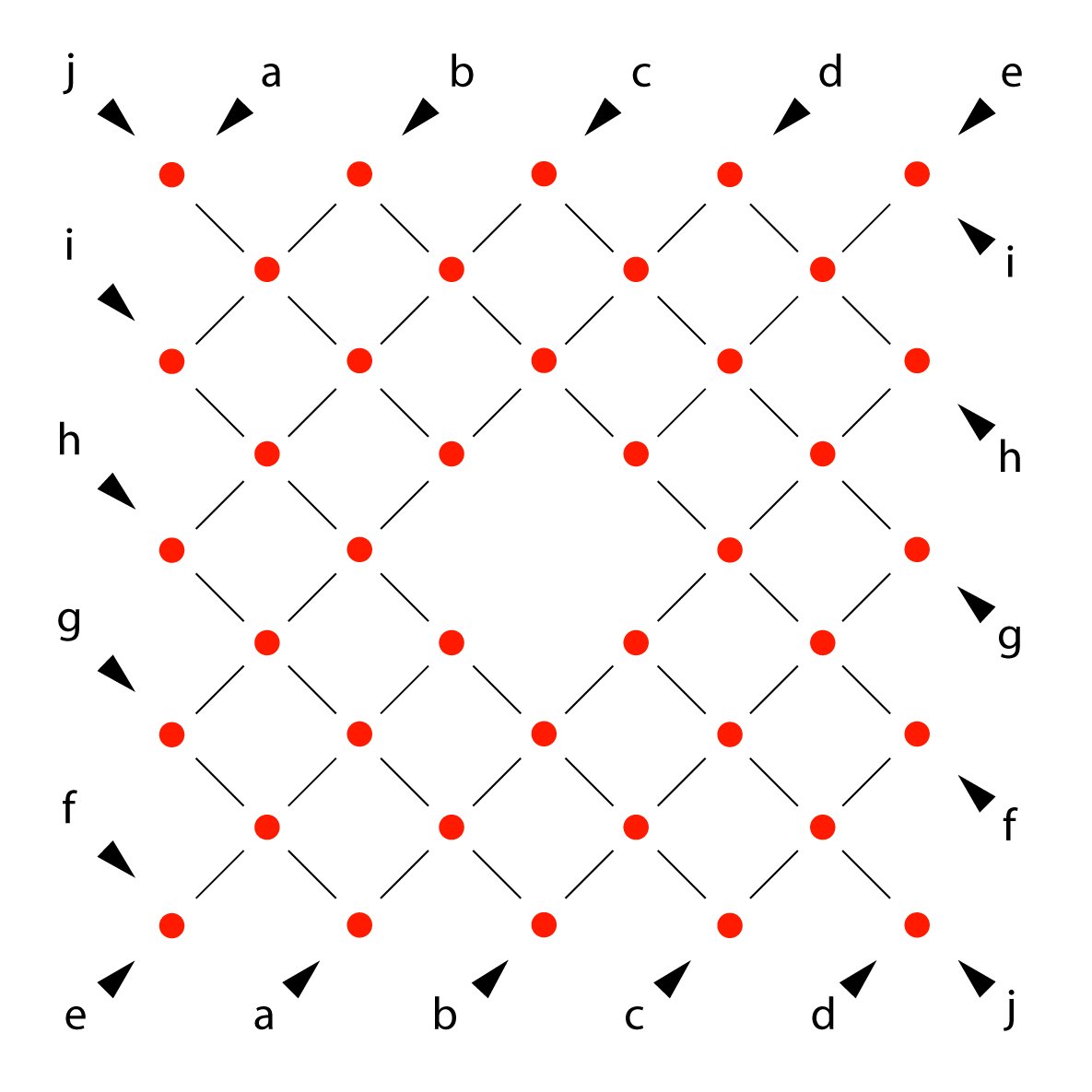}
\par\end{centering}

\protect\caption{\label{fig:punctSquareBoard}Square board with SIN $\left\{ 0,2\right\} $}
\end{figure}

\end{example}

\subparagraph{Biregular Hexagons\label{exa:Boards-from-Birregular-hex}}
\begin{example}
Fig. \ref{fig:hexFacesSymmetricBoard} shows a biregular hexagon of
order three on which we have chosen face centers as board points.
An asterism is made up of the eighteen points pointed to by the same
letter, and belongs to one of three parallel classes: $\left\{ a,b,c\right\} $,
$\left\{ d,e,f\right\} $ or $\left\{ g,h,i\right\} $. We have then
a board with a $18$-uniform resolution. It is easy to verify that
\textit{$D_{6}$} acts on the asterisms and transitively on the classes,
so the board is class-symmetric\footnote{\textit{$D_{3}$} also acts on the asterisms here but we choose to
focus on the largest symmetry group.}. It is also easy to see that its SIN is $\left\{ 0,6,8\right\} $.
We can have symmetric boards like this one for any order \textit{n}
of the biregular hexagon ($6n^{2}$ points) with three uniform parallel
classes, each with \textit{n} asterisms of size $6n$ each.
\begin{figure}
\begin{centering}
\includegraphics[height=5cm]{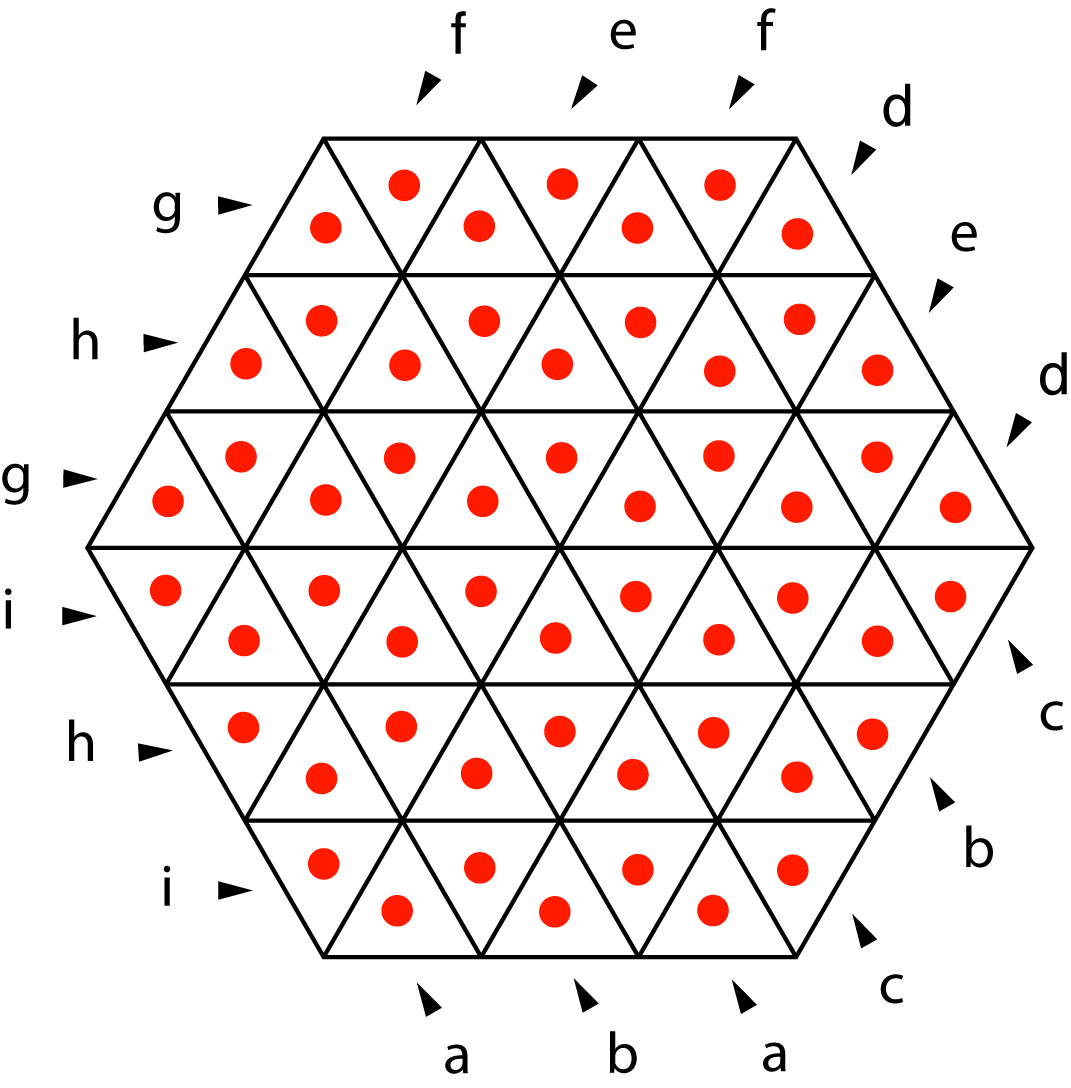}
\par\end{centering}

\protect\caption{\label{fig:hexFacesSymmetricBoard}Hexagonal board with points at
the barycenters}
\end{figure}

\end{example}

\begin{example}
Fig. \ref{fig:hexSymBoardEdgesOnSide} shows a biregular hexagon of
order three on which we have chosen edge centers as board points.
An asterism is made up of sets of points pointed to by the same letter.
Each asterism has fourteen points and belongs to one of three parallel
classes: $\left\{ a,b,c\right\} $, $\left\{ d,e,f\right\} $ or $\left\{ g,h,i\right\} $.
The board has a $14$-uniform resolution. It is simple to verify that
\textit{$D_{6}$} acts transitively on the classes, so the board is
class-symmetric.
\begin{figure}
\begin{centering}
\includegraphics[height=5.5cm]{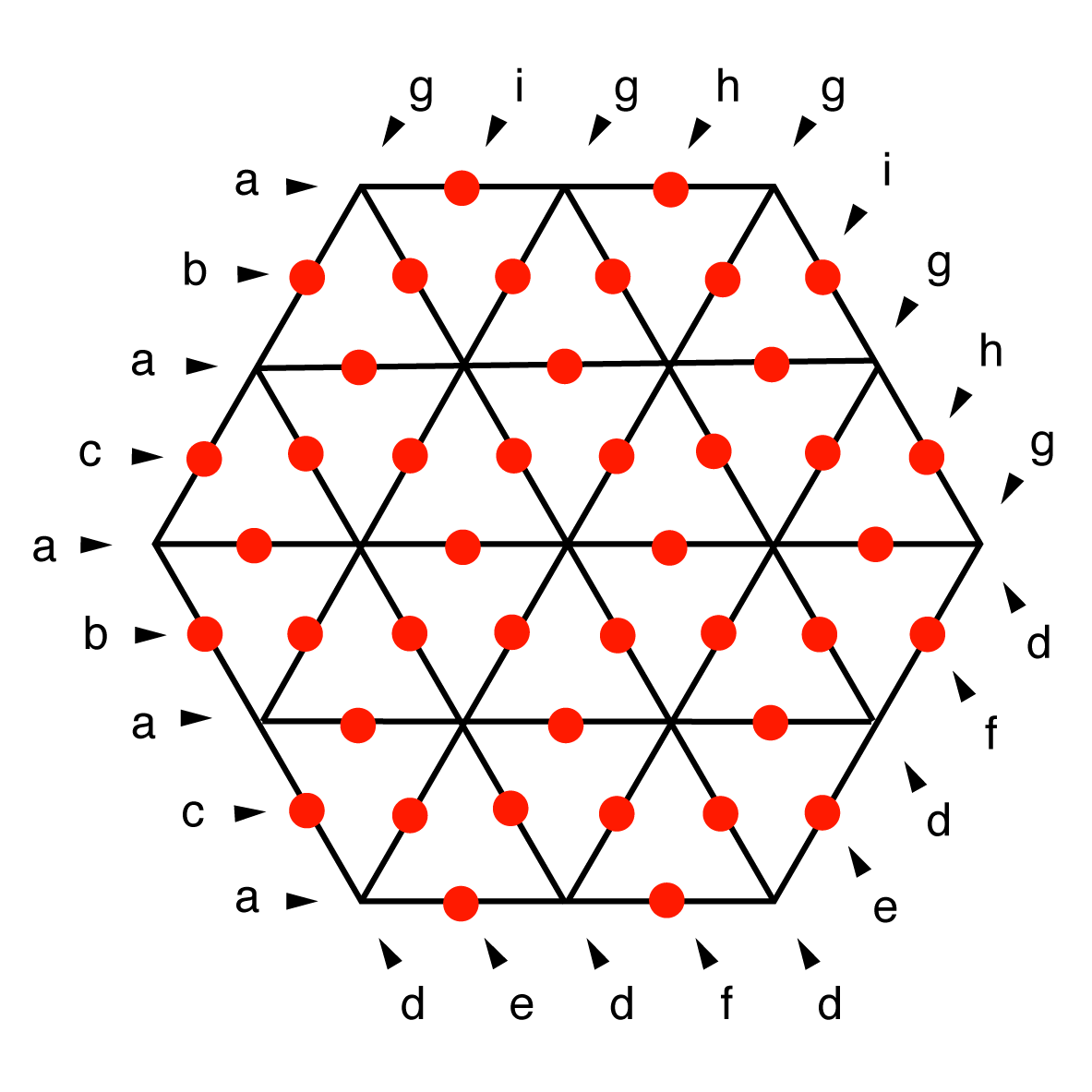}
\par\end{centering}

\protect\caption{\label{fig:hexSymBoardEdgesOnSide}Hexagonal board with points at
the edge midpoints}
\end{figure}

\end{example}

\subparagraph{Regular Polyhedra with Biregular Faces\label{exa:Boards-from-Polyhedra with bireg faces}}

\noindent There are five finite real convex regular polyhedra, and
three polyhedral full symmetry groups (see Table \ref{tab:symmetryGroupsPlatonicSolids-1}).
Each of the following examples consists of a polyhedral source of
symmetry (see sect. \ref{sub:Sources of Symmetry}) and a board based
on it. The corresponding pictures are meant to be folded along the
inner black lines and glued along the outer ones. With the exception
of those in the example dodecahedral board, every asterism in the
boards is made up of points enclosed by two lines of the same color.
The asterism becomes a closed band \textendash so to speak\textendash{}
when the board is folded.
\begin{table}
\begin{centering}
\begin{tabular}{ccc}
\textbf{\scriptsize{}group} & \textbf{\scriptsize{}regular polyhedron} & \textbf{\scriptsize{}order}\tabularnewline
\hline 
\noalign{\vskip0.1cm}
{\footnotesize{}tetrahedral} & {\footnotesize{}tetrahedron} & {\footnotesize{}24}\tabularnewline
\noalign{\vskip0.1cm}
{\footnotesize{}octahedral} & {\footnotesize{}cube, octahedron} & {\footnotesize{}48}\tabularnewline
\noalign{\vskip0.1cm}
{\footnotesize{}icosahedral} & {\footnotesize{}dodecahedron, icosahedron} & {\footnotesize{}120}\tabularnewline
\end{tabular}
\par\end{centering}

\protect\caption{\label{tab:symmetryGroupsPlatonicSolids-1} Polyhedral symmetry groups
of the platonic solids }
\end{table}

\begin{example}
\textbf{Tetrahedral Board}. Fig. \ref{fig:biregularTetrahedron} (left)
shows a tetrahedral source. On the right we have selected face centroids
(triangle orthocenters) as board points. As indicated above, asterisms
are made up of points contained between lines of the same color. The
board has SIN $\left\{ 0,4\right\} $\footnote{Because when unfolded, the board is also a Latin triangle like the
one in Fig. \ref{fig:latin triangle faces Small}, which has SIN $\left\{ 0,4\right\} $.} and is resolvable. It is also class-symmetric, since there are three
parallel classes that are acted upon transitively by the tetrahedral
group.
\begin{figure}
\begin{centering}
\includegraphics[height=5cm]{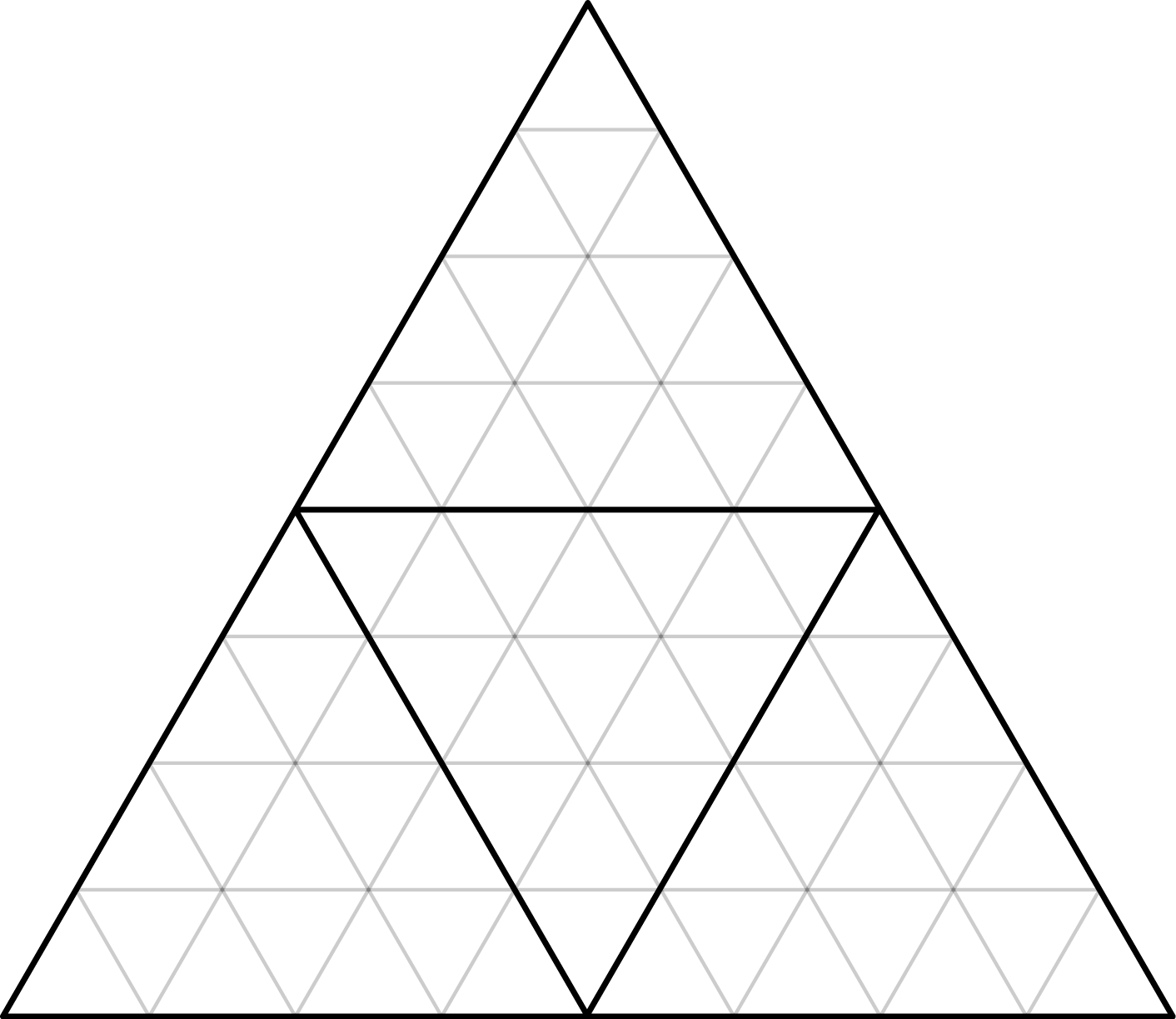}~~~~~~~~~~\includegraphics[height=5cm]{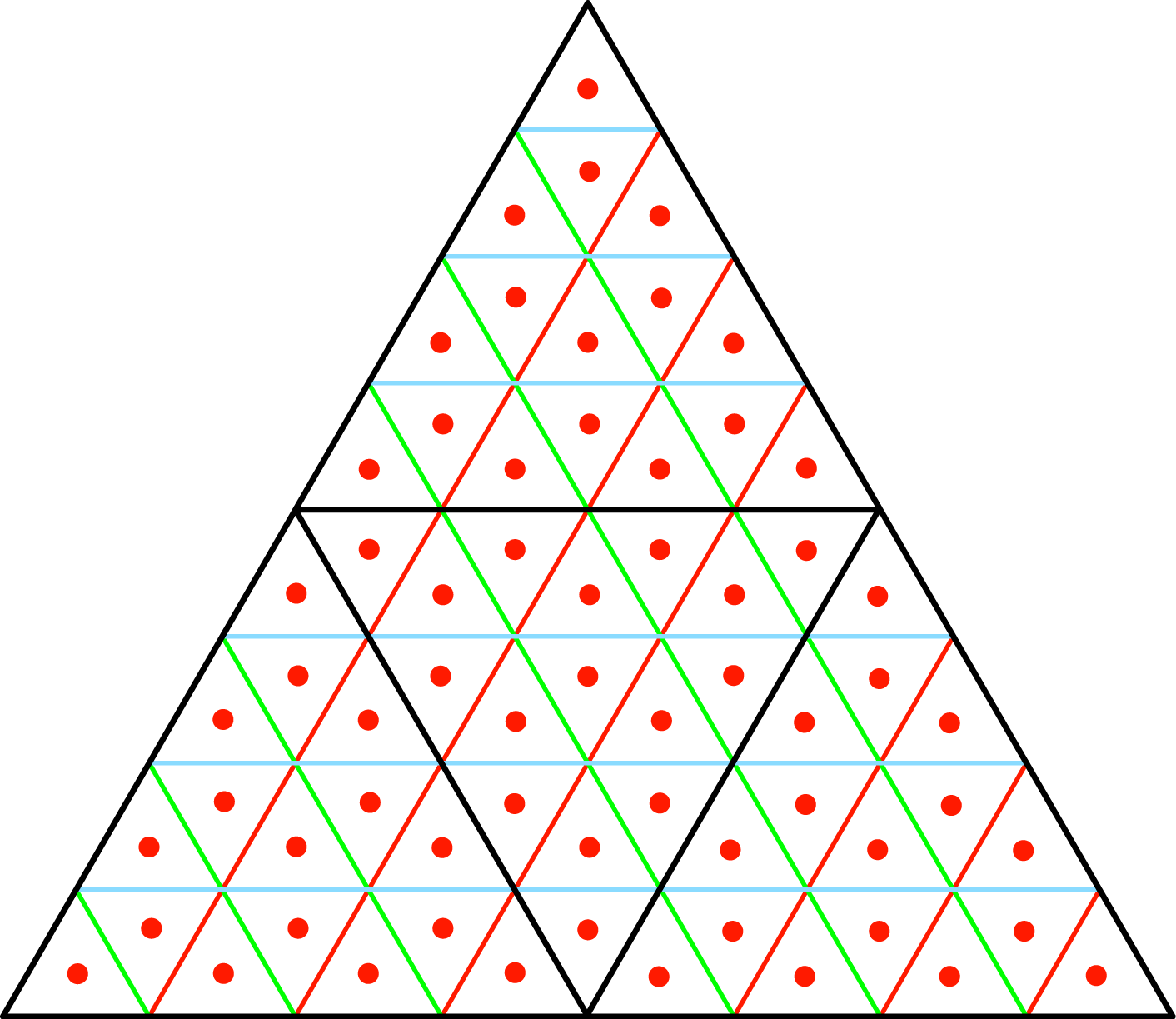}
\par\end{centering}

\protect\caption{\label{fig:biregularTetrahedron}Tetrahedral source and board}
\end{figure}

\end{example}

\begin{example}
\textbf{Cubical} \textbf{Board}. The picture in Fig. \ref{fig:biregularCube}
is a cubical source. In Fig. \ref{fig:latin cube symmetric board}
we have chosen the face centroids as board points. The asterisms,
defined as before as the points between lines of the same color, are
acted upon by the octahedral group. The board has SIN $\left\{ 0,8\right\} $
and no parallel classes, so it is not resolvable.
\begin{figure}
\begin{centering}
\includegraphics[height=8cm]{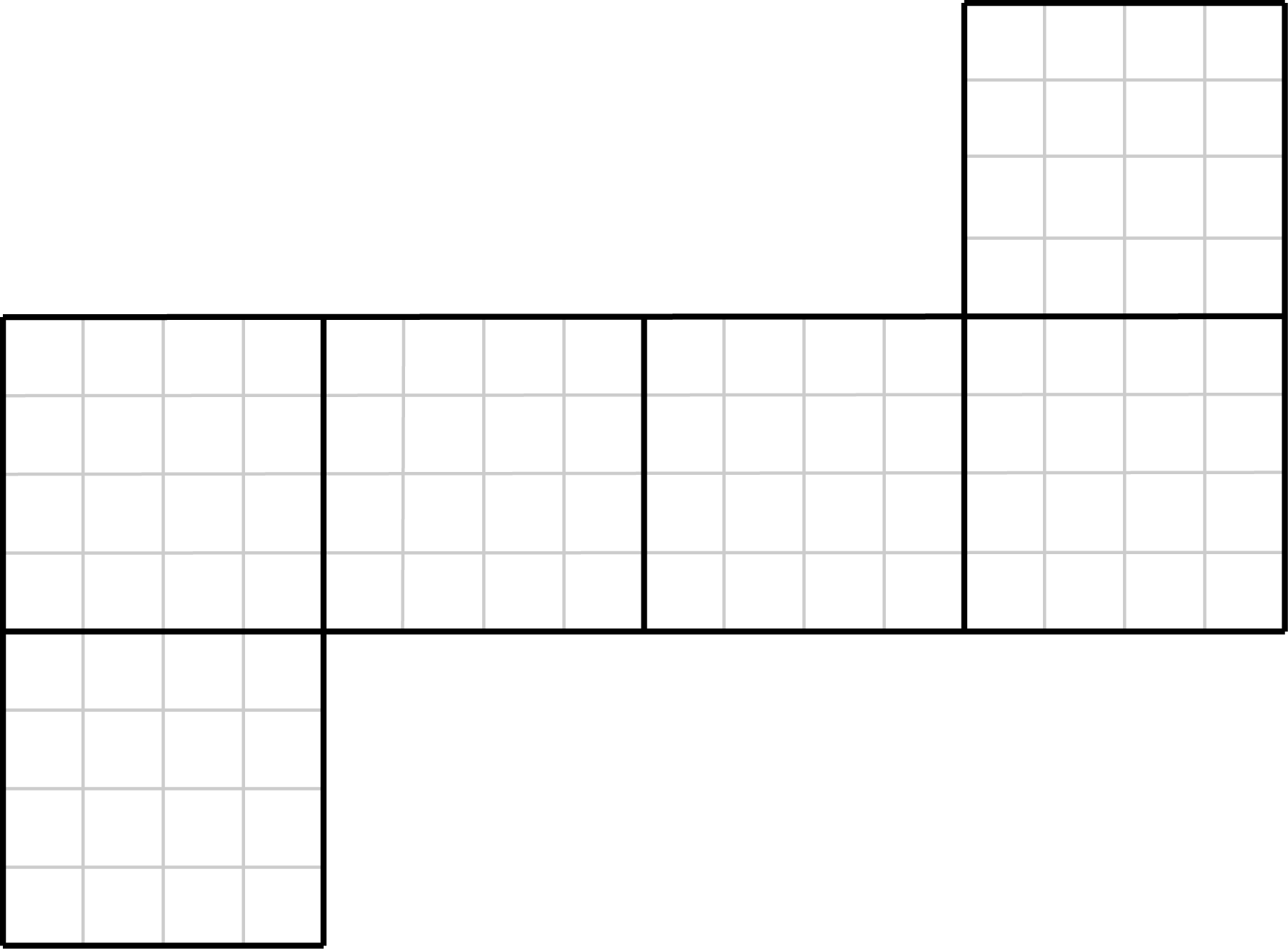}
\par\end{centering}

\protect\caption{\label{fig:biregularCube}Cubical source}
\end{figure}
\begin{figure}
\begin{centering}
\includegraphics[height=8cm]{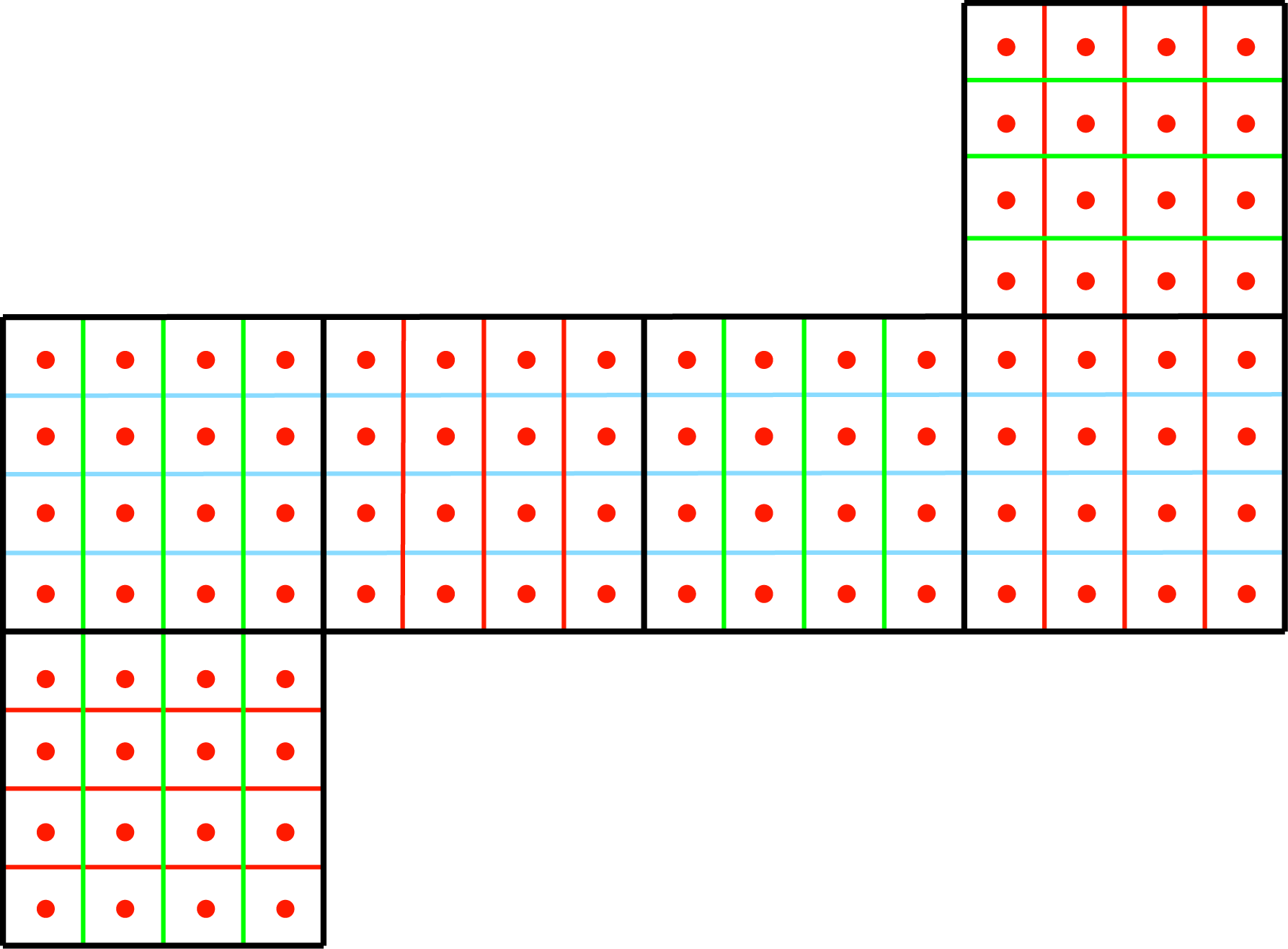}
\par\end{centering}

\protect\caption{\label{fig:latin cube symmetric board}Cubical board}
\end{figure}

\end{example}

\begin{example}
\textbf{Octahedral} \textbf{Board}. The picture in Fig. \ref{fig:biregularOctahedron}
is an octahedral source. In Fig. \ref{fig:latin octahedron symmetric board}
we have chosen face orthocenters as board points and asterisms as
indicated before. The board is symmetric since the asterisms are acted
upon by the octahedral group. The board has SIN $\left\{ 0,4\right\} $
and no parallel classes, and hence it is not resolvable.
\begin{figure}
\begin{centering}
\includegraphics[height=8cm]{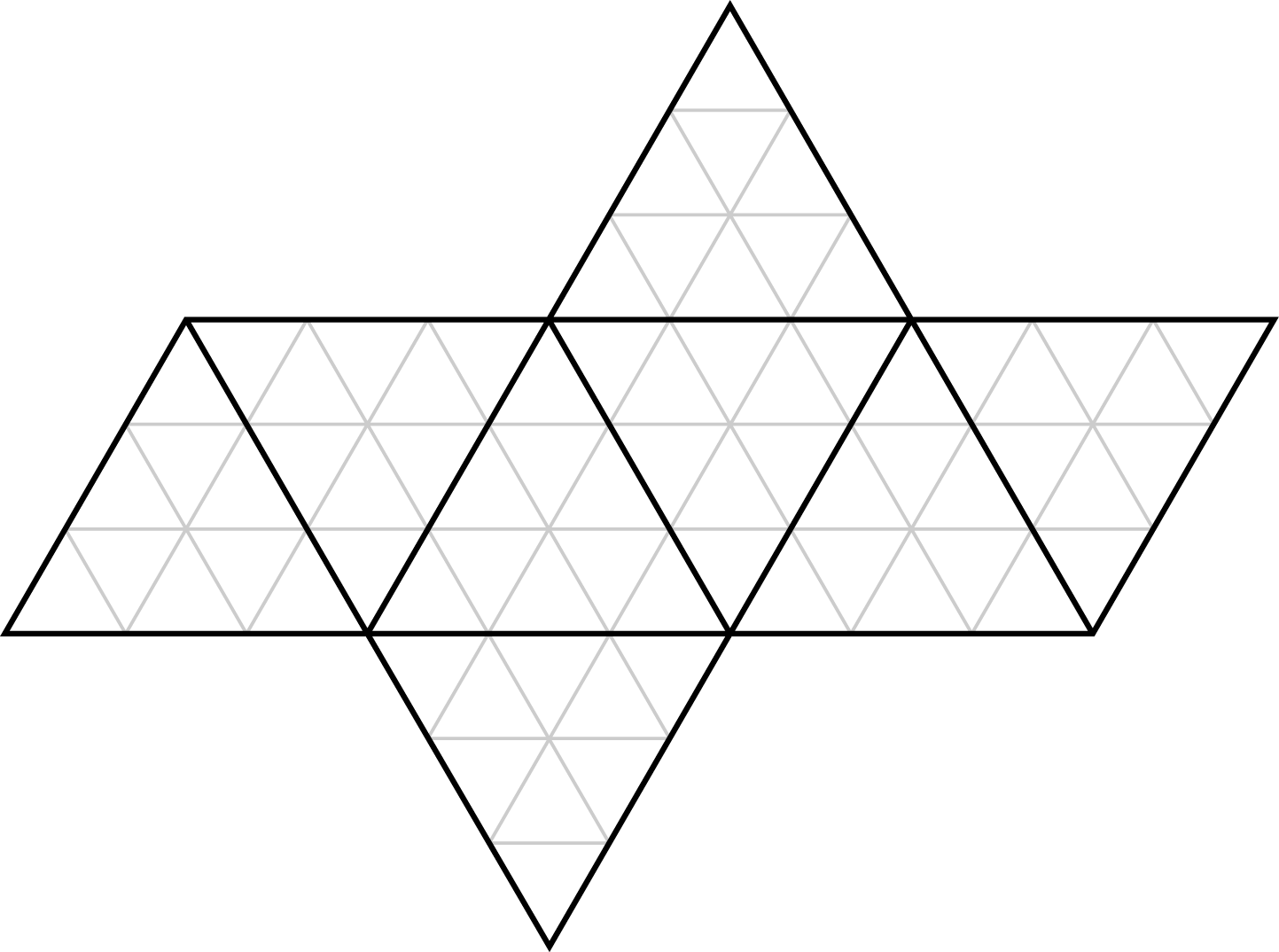}
\par\end{centering}

\protect\caption{\label{fig:biregularOctahedron}Octahedral source}
\end{figure}
\begin{figure}
\begin{centering}
\includegraphics[height=8cm]{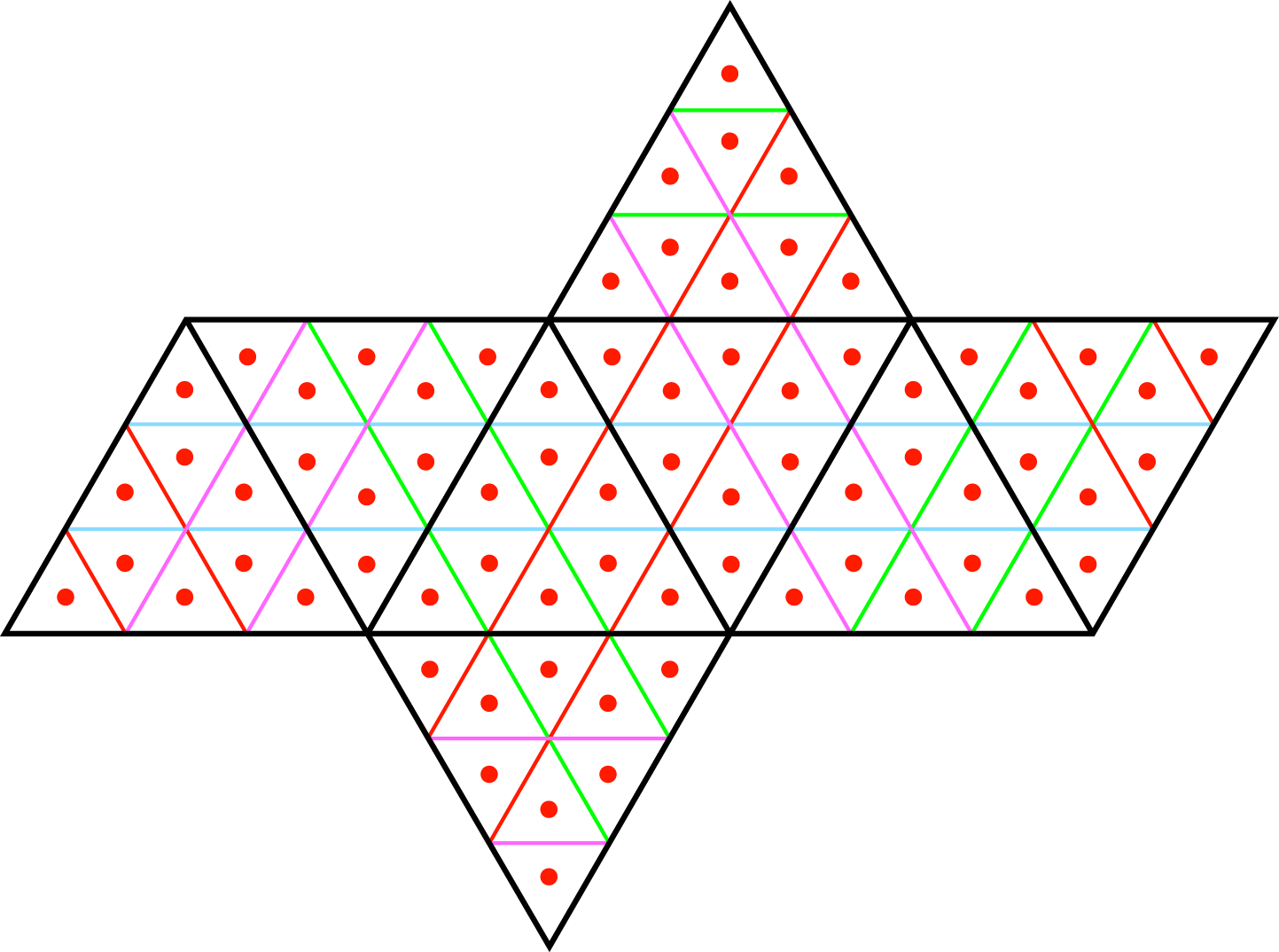}
\par\end{centering}

\protect\caption{\label{fig:latin octahedron symmetric board}Octahedral board}
\end{figure}

\end{example}

\begin{example}
\textbf{Icosahedral} \textbf{Board}. The picture in Fig. \ref{fig:biregularIcosahedron}
is an icosahedral source. In Fig. \ref{fig:latin icosahedron symmetric board}
we have chosen face orthocenters as board points. The asterisms, defined
as above, are acted upon by the icosahedral group, so the board is
symmetric. The board has SIN $\left\{ 0,4\right\} $ and no parallel
classes.
\begin{figure}
\begin{centering}
\includegraphics[height=5.8cm]{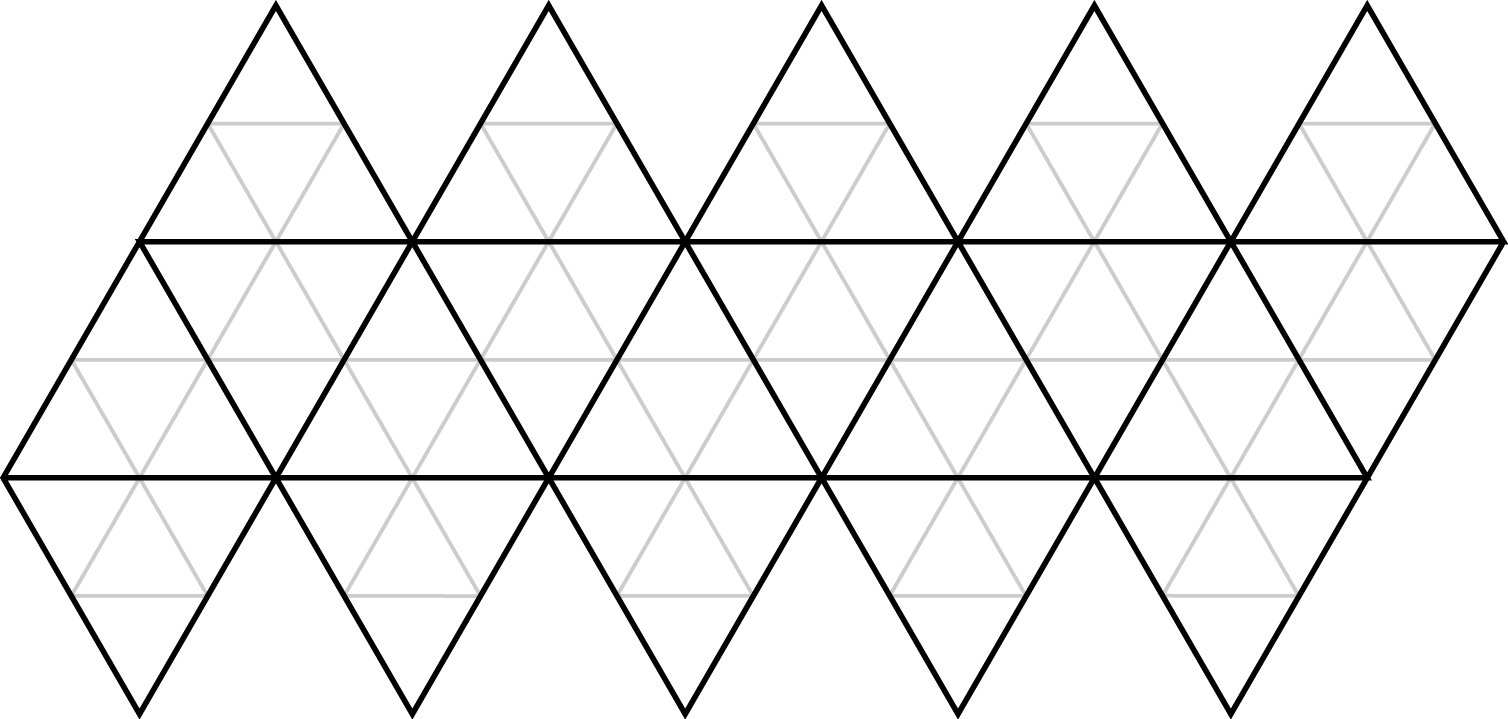}
\par\end{centering}

\protect\caption{\label{fig:biregularIcosahedron}Icosahedral source}
\end{figure}
\begin{figure}
\begin{centering}
\bigskip{}
\includegraphics[height=5.8cm]{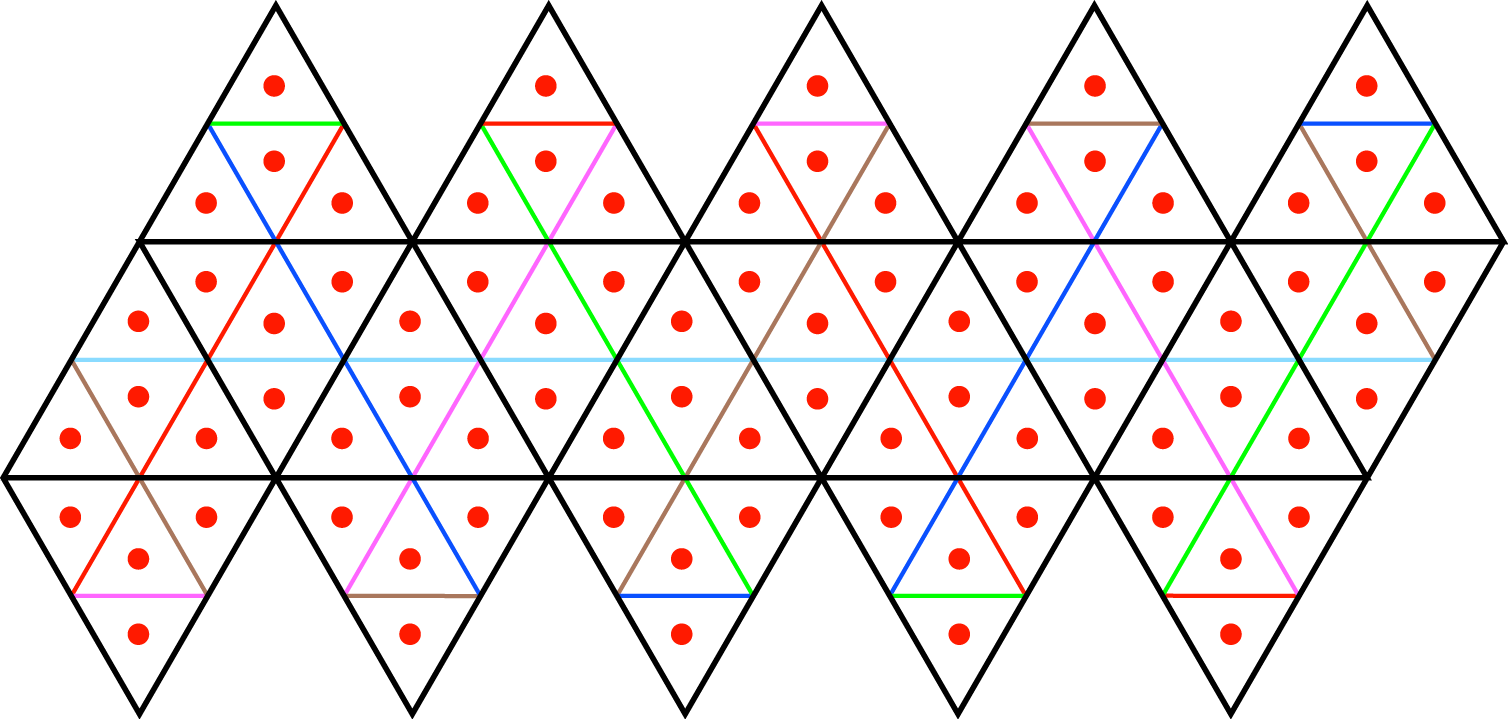}
\par\end{centering}

\protect\caption{\label{fig:latin icosahedron symmetric board}Icosahedral board}
\bigskip{}
\end{figure}

\end{example}

\begin{example}
\textbf{Dodecahedral Board}. The picture in Fig. \ref{fig:unfoldedDodecahedron}
is a dodecahedral source. Its pentagonal faces are trivial biregular
pentagons. In Fig. \ref{fig:latinBaseDodecahedronEdgesAndLines} we
have chosen the edge centers as board points. Each asterism is made
up of points lying on a drawn line of a particular color. Once folded,
each asterism becomes a closed polyline on the dodecahedron surface.
The board is symmetric because the asterisms so defined are acted
upon by the icosahedral group. It has SIN $\left\{ 0,2\right\} $
and no parallel classes, and so it is not resolvable.
\begin{figure}
\begin{centering}
\includegraphics[height=6cm]{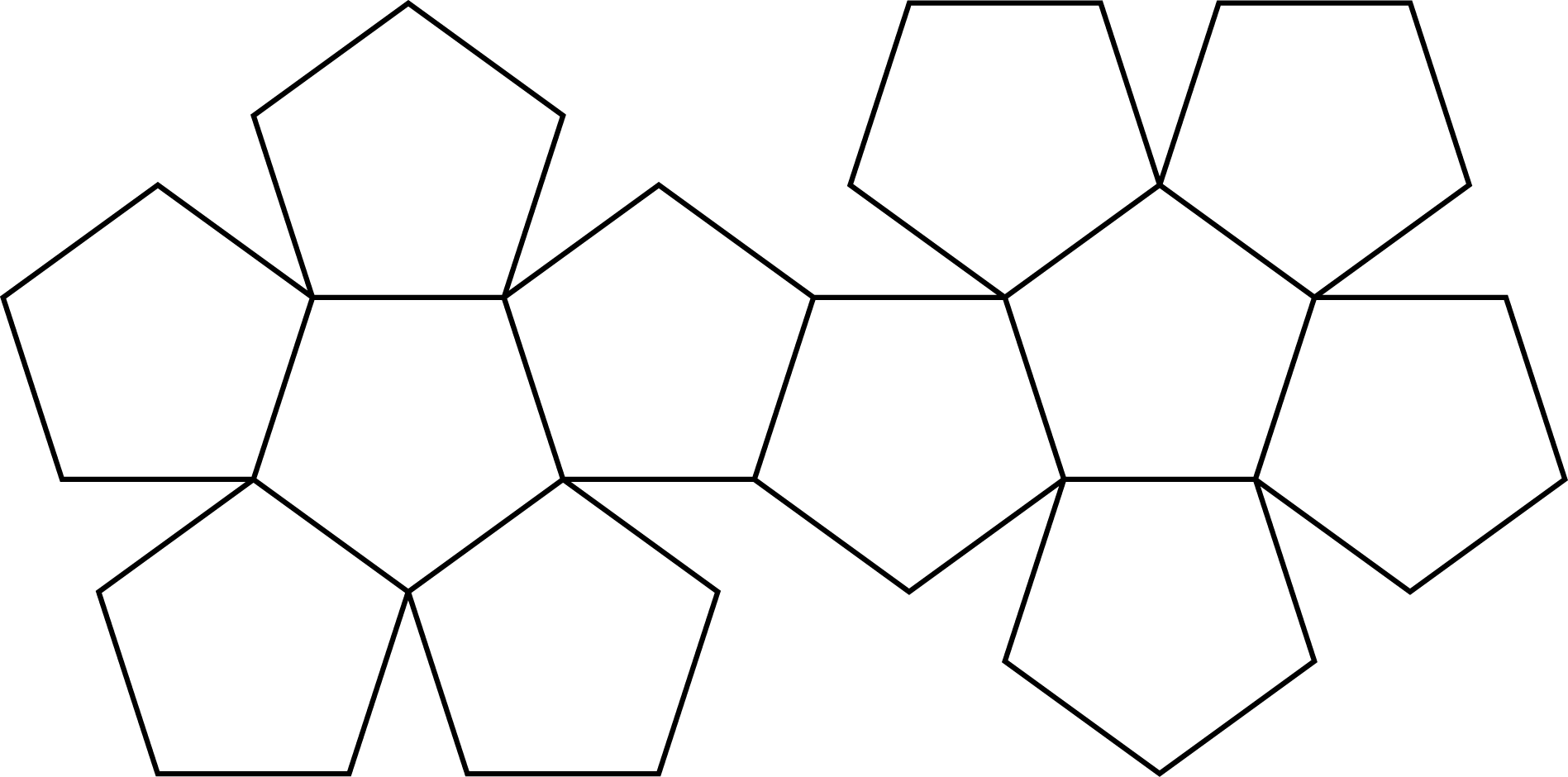}
\par\end{centering}

\protect\caption{\label{fig:unfoldedDodecahedron}Dodecahedral source}
\end{figure}
\begin{figure}
\begin{centering}
\includegraphics[height=6cm]{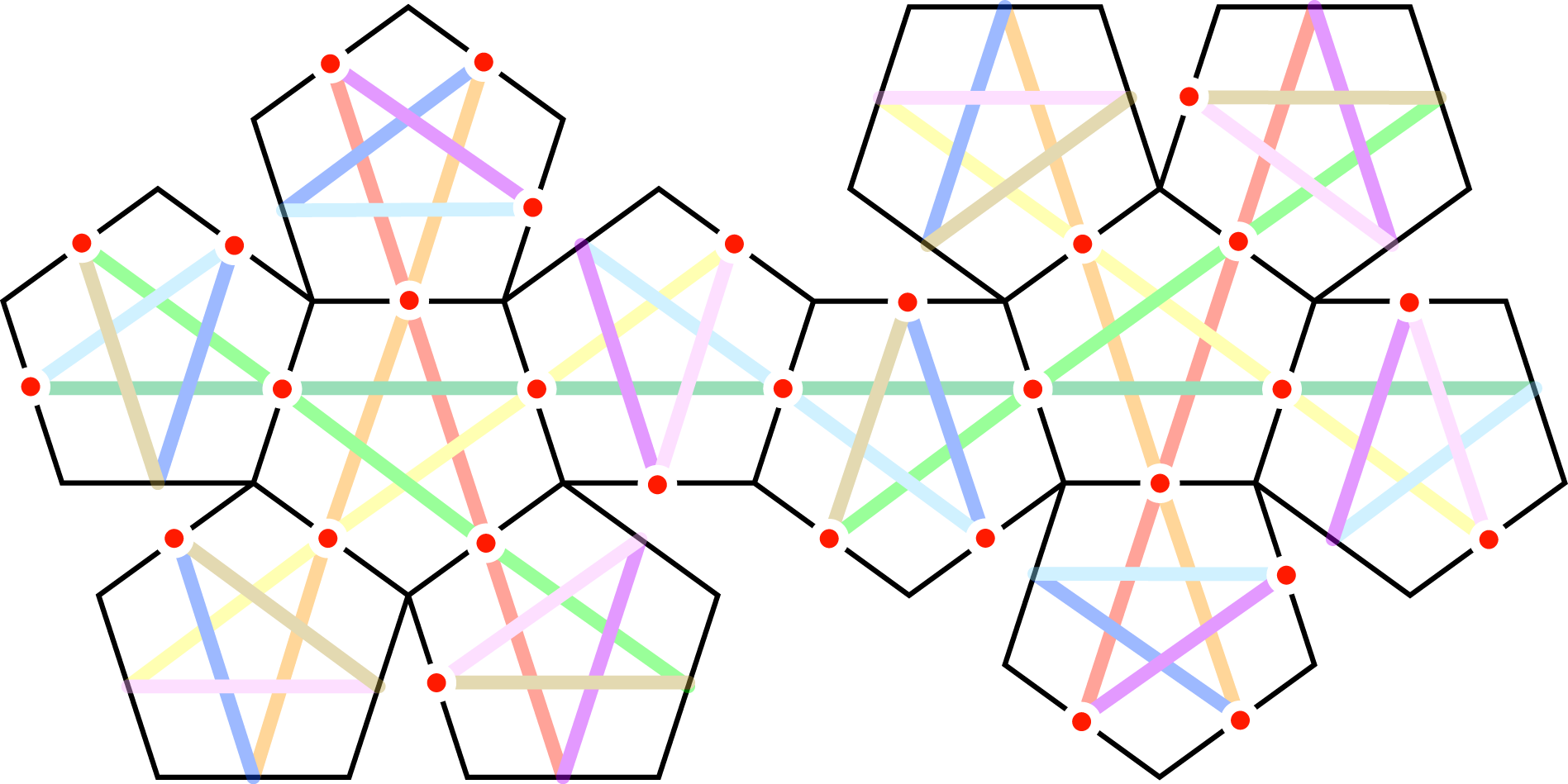}
\par\end{centering}

\protect\caption{\label{fig:latinBaseDodecahedronEdgesAndLines}Dodecahedral board}
\end{figure}

\end{example}

\paragraph{Extra General Tilings}

We can also enrich sources by tiling some of its elements (edges,
faces, hyperplanes) with tiles that are not necessarily regular.
\begin{example}
\noindent The picture in Fig. \ref{fig:latinBaseDodecahedronWtrianglesInFaces}
is meant to be folded along the inner black lines and glued along
the outer ones. The resulting dodecahedron, with faces tiled with
non-equilateral triangles, is a source with icosahedral symmetry (see
Table \ref{tab:symmetryGroupsPlatonicSolids-1}).
\begin{figure}
\begin{centering}
\includegraphics[height=6cm]{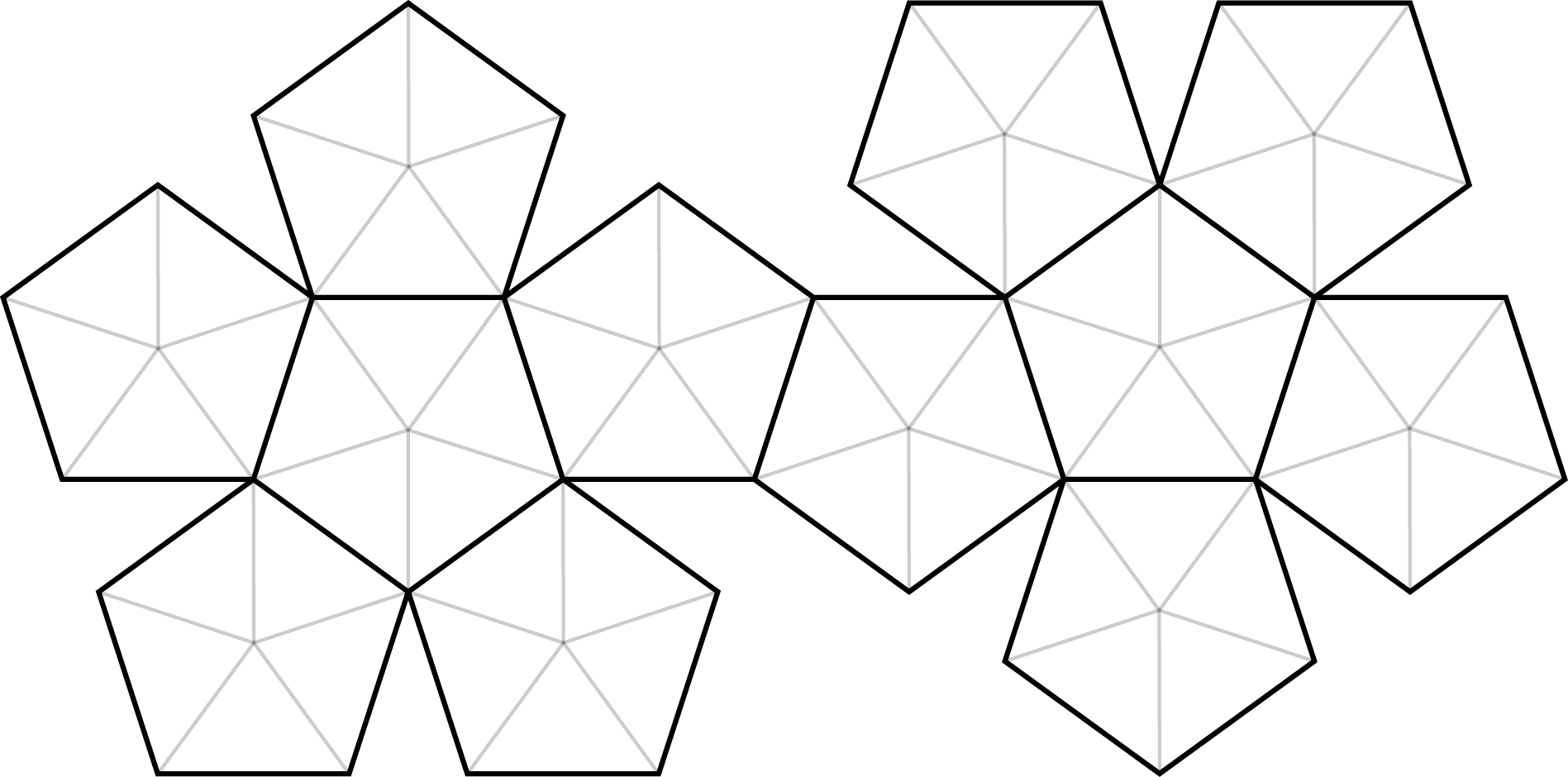}
\par\end{centering}

\protect\caption{\label{fig:latinBaseDodecahedronWtrianglesInFaces}Source with dodecahedral
symmetry}
\end{figure}
\begin{figure}
\begin{centering}
\includegraphics[height=6cm]{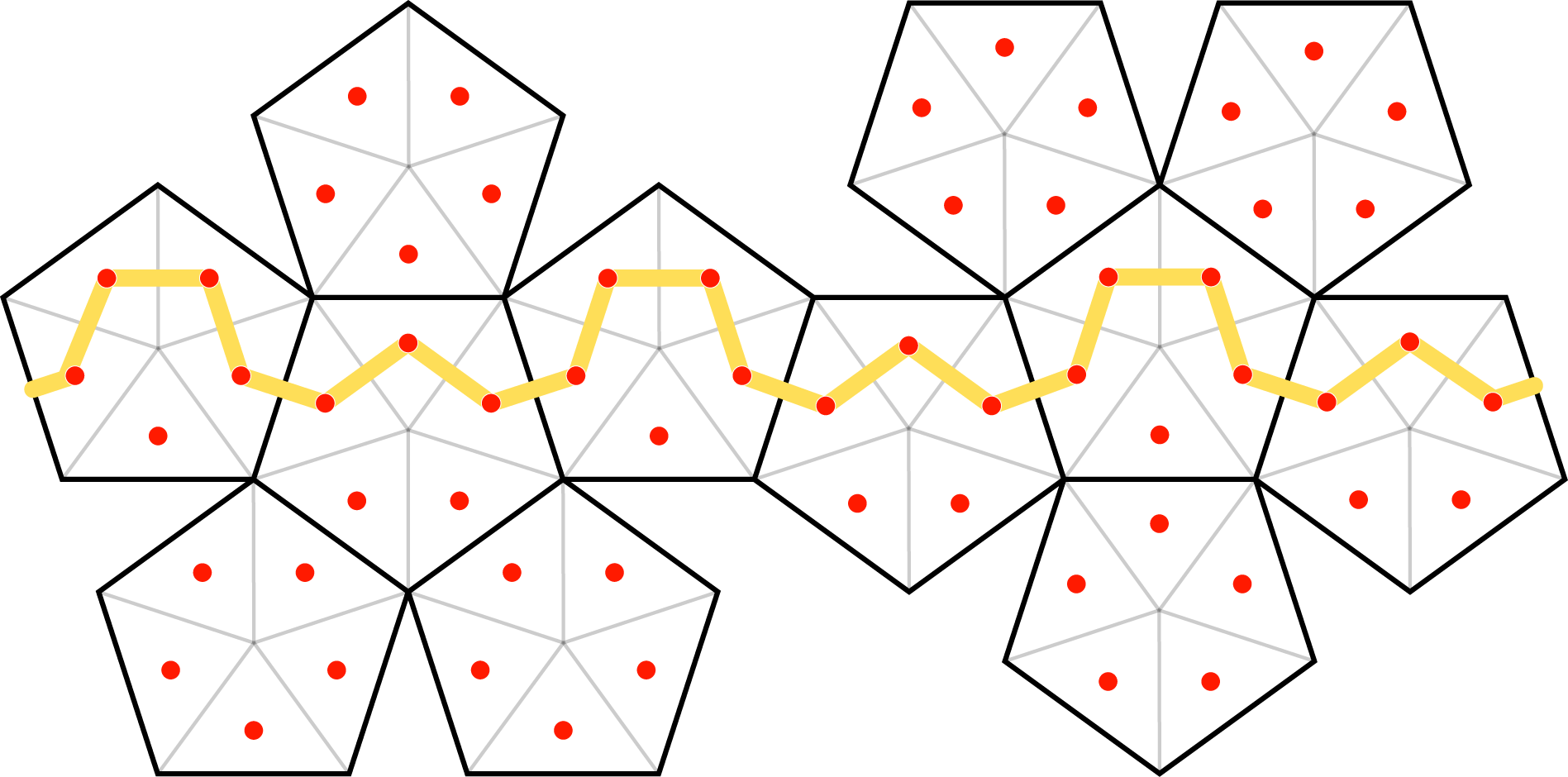}
\par\end{centering}

\protect\caption{\label{fig:designPtsDodecahedronWtrianglesInFacesAndLineHint}Points
and asterism on a dodecahedral board}
\end{figure}
\begin{figure}
\begin{centering}
\includegraphics[angle=90,origin=c,height=15cm]{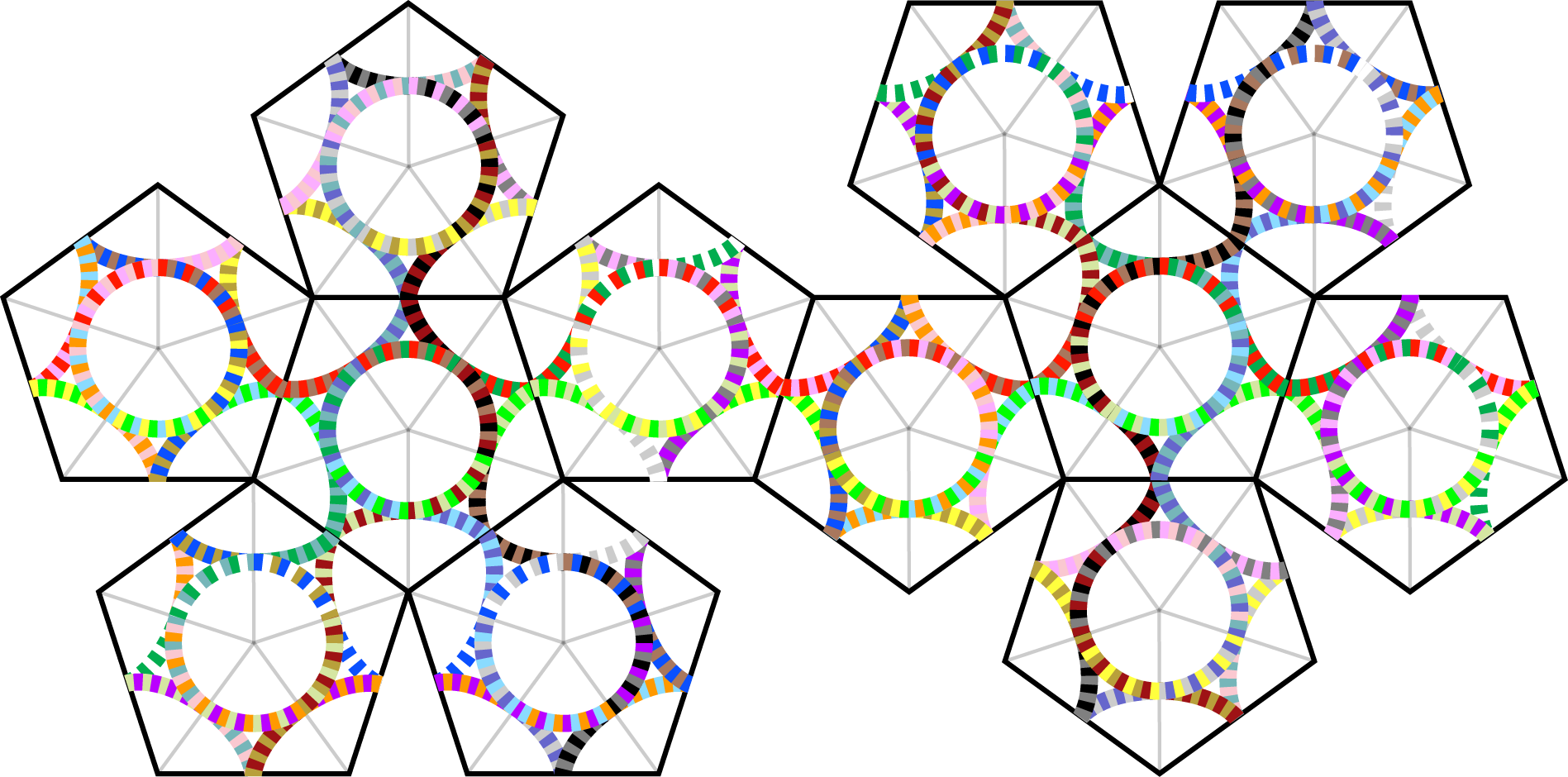}
\par\end{centering}

\protect\caption{\label{fig:easyDodecahedronFaces}Points and asterisms on a dodecahedral
board}
\end{figure}

\noindent In Fig. \ref{fig:designPtsDodecahedronWtrianglesInFacesAndLineHint}
we have a board whose points are the orthocenters of the triangles
inside the pentagons. Each of its asterisms is made up of points lying
on lines like the one shown in yellow. The asterisms, each indicated
by a colored dashed line, are shown in Fig. \ref{fig:easyDodecahedronFaces}.
Once the board folded, the portions of each asterism become one on
the polyhedron surface. The board has no parallel classes and thus
it is not resolvable, but it is symmetric because the asterisms are
acted upon by the icosahedral group.
\end{example}

\subsection{Symmetric Boards from Asymmetric Boards}

In Sect. \ref{sub:Sources of Symmetry} we have obtained symmetric
boards for sources of symmetry. We discuss next how to obtain symmetric
boards from asymmetric ones.

\subsubsection{Board Symmetrization Problem (BSP)}
\begin{defn}
\emph{Given an asymmetric board, find symmetric boards isomorphic
to it}. We call this problem the \emph{Board Symmetrization Problem}
(BSP).\end{defn}
\begin{example}
Fig. \ref{fig:FanoPlaneUnsymmetricAndSymmetric} (left) shows the
Fano Plane as an asymmetric board. We would like to find boards like
the one in the center, which is both isomorphic to the first \textendash this
is easily verified by inspection\textendash{} and symmetric \textendash $D_{3}$
acts on the asterisms\footnote{The boards derived from this one by the action of $D_{3}$ are also
solutions to the problem.}. 
\begin{figure}
\begin{centering}
\includegraphics[height=2.5cm]{FanoPlaneUnsymmetricUnnumbered}~~~~~\includegraphics[height=2.5cm]{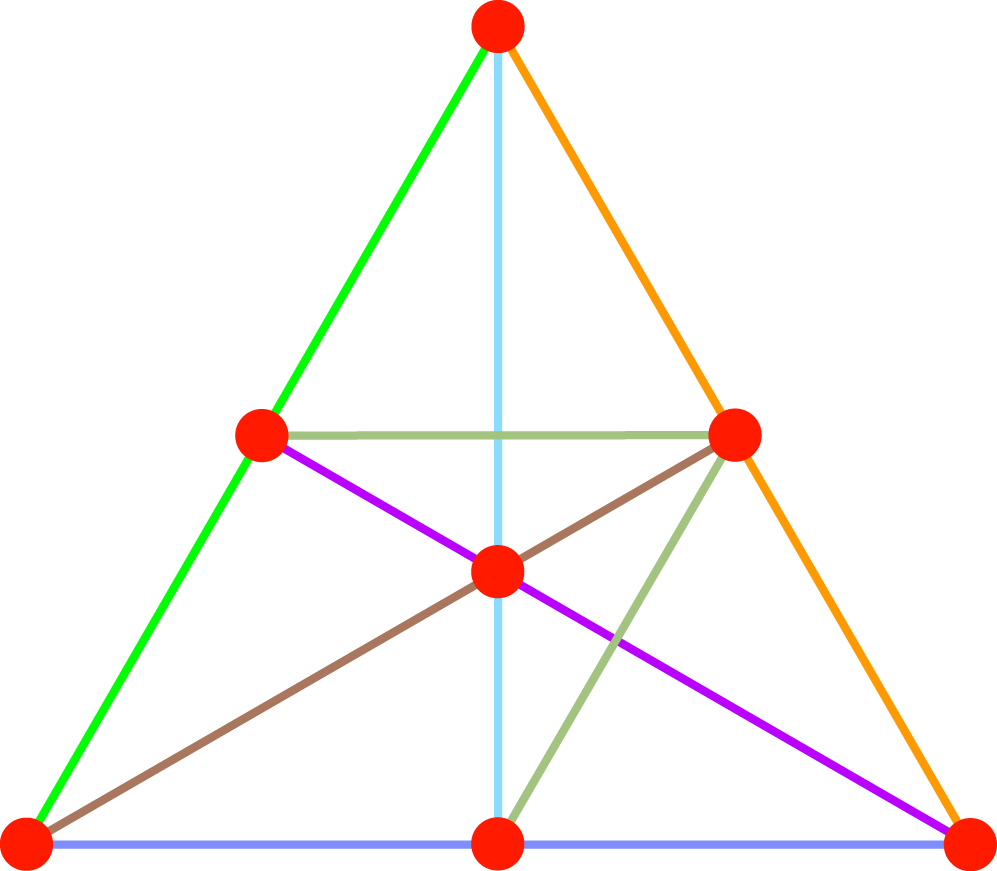}~~~~~~~~~\includegraphics[height=2.5cm]{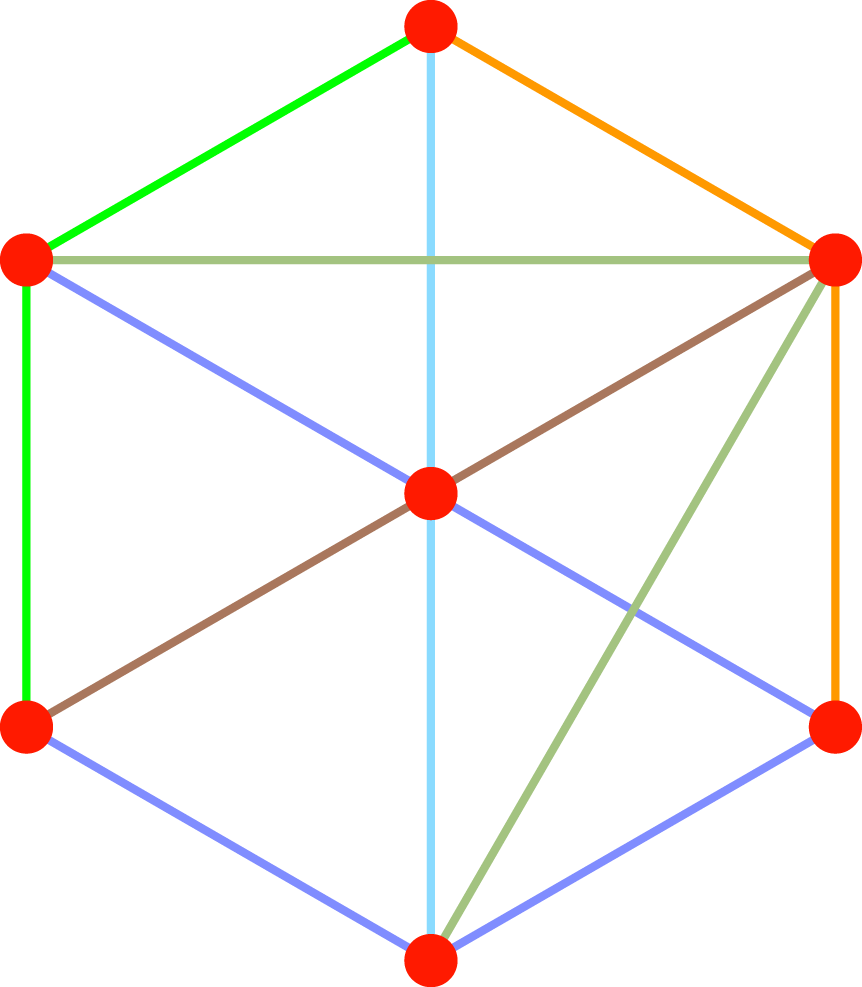}
\par\end{centering}

\protect\caption{\label{fig:FanoPlaneUnsymmetricAndSymmetric}The Fano Plane in three
different boards}
\end{figure}

\end{example}

\subsubsection{Symmetric Boards with unknown Automorphism Group\label{sub:Asymmetric-Boards-with-unknown automorphism group}}

\noindent Let \textbf{$B=\left(P,C\right)$} be our asymmetric board
and let's suppose that $Aut\left(B\right)$ is unknown. We can approach
the BSP geometrically by first establishing a bijection between points
in \textit{$P$} and other points in $P'$, that may reside in the
same space or in a different one. We define an asterism as the set
of points in $P'$ that are linked by the bijection to points in an
asterism in $P$. If we group all resulting asterisms in constellation
$C'$ we obtain the isomorphic board $B'=\left(P',C'\right)$.

\noindent Now we consider a set $S_{SG}$ of candidate symmetry groups
(see Sect. \ref{sub:Symmetry-Groups}). For each of them we arrange
the points in $P'$ so that the group acts on them, then repeatedly
permute the points to see if the group could also act on the asterisms
in $C'$. When this is achieved the resulting board is a solution
to the BSP for $B$.
\begin{example}
\noindent By taking $D_{3}$ as a candidate symmetry group for the
Fano Plane in Fig. \ref{fig:FanoPlaneUnsymmetricAndSymmetric} (left)
and arranging its points as described we have obtained the solution
to the BSP shown in the center. We could try now $D_{6}$. Fig. \ref{fig:FanoPlaneUnsymmetricAndSymmetric}
(right) shows a failed attempt to rearrange the points in which $D_{6}$
acts on them but not on the asterisms (a $60\text{º}$ counter-clockwise
rotation for example does not take the red asterism to any other).
\end{example}
\noindent As we saw in sect. \ref{sec:B1 automorphisms from symmetry},
a symmetric board $B'$ isomorphic to another one $B$ provides a
subgroup of $Aut\left(B\right)$. So if $Aut\left(B\right)$ is unknown,
finding solutions to the BSP provides at least some previously unknown
automorphisms.

\subsubsection{Symmetric Boards with known Automorphism Group}

\begin{figure}
\begin{centering}
\includegraphics[height=3cm]{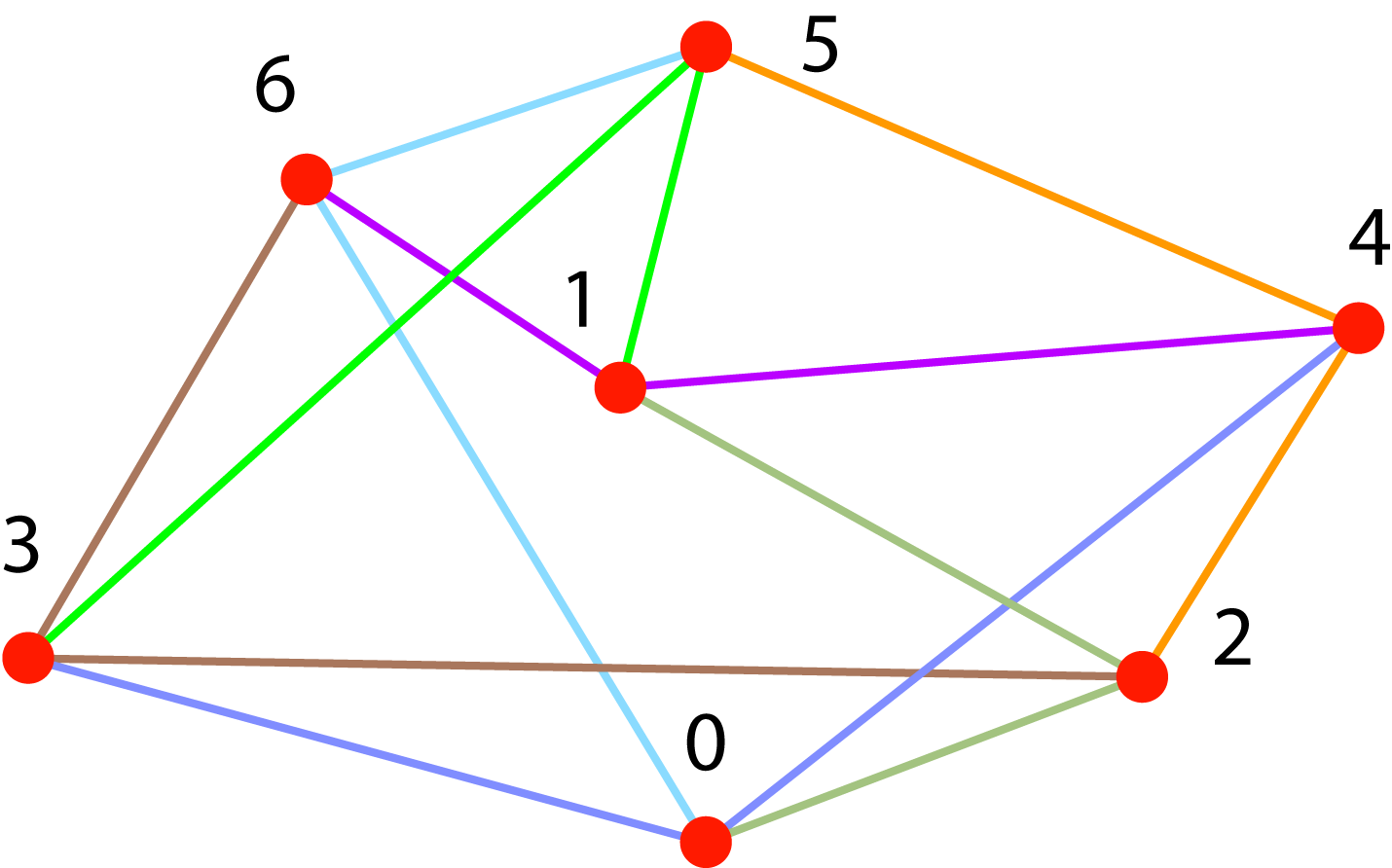}~~~~~~~~\includegraphics[height=3cm]{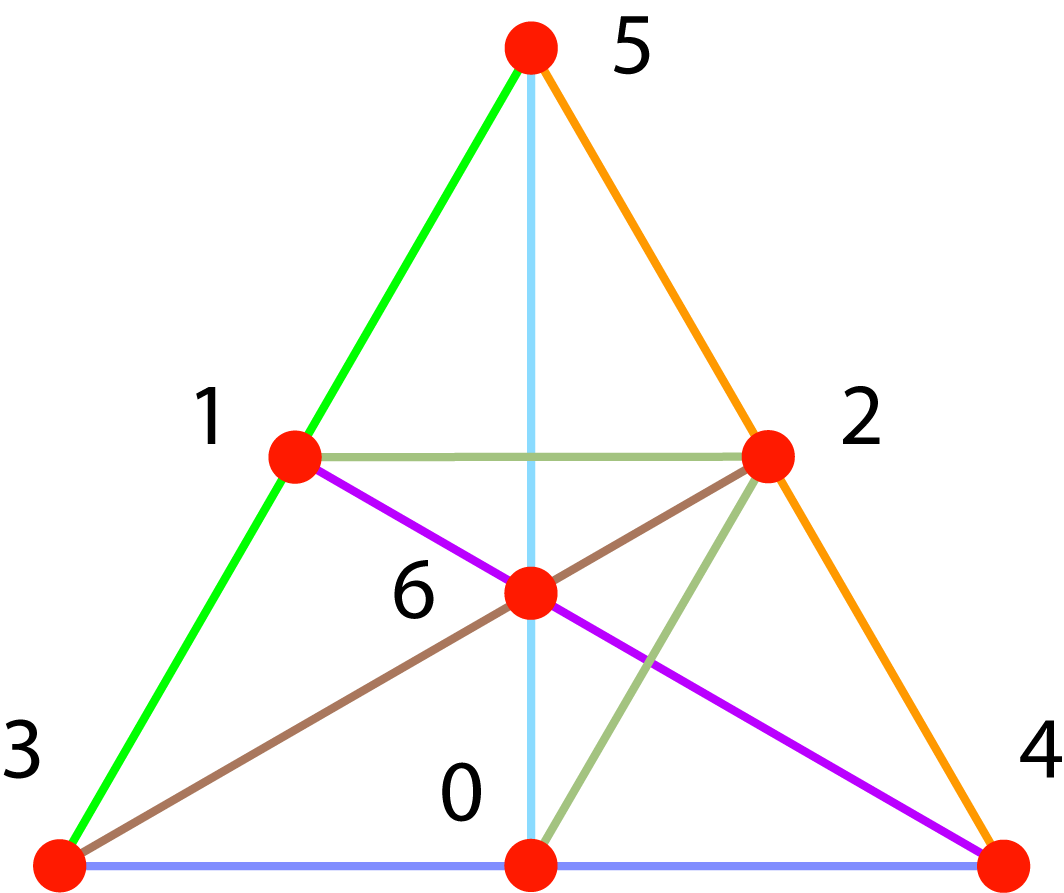}
\par\end{centering}

\protect\caption{\label{fig:Fano Plane asymmetrical and symmetrical numbered}Asymmetric
and symmetric boards for the Fano Plane}

\end{figure}

On the other hand, if \emph{$Aut(B)$} is known, we can further simplify
the BSP: Lemma \ref{lem:AT(B) included in Aut(B)} tells us that we
need to focus only in symmetry groups whose actions are linked to
permutation groups that are subgroups of \emph{$Aut(B)$}\footnote{This leaves us of course with the additional problem of finding the
subgroups of a group.}.
\begin{example}
The automorphism group of the Fano Plane $Aut\left(FP\right)$ has
$168$ elements \cite{Handbook of Combinatorial Designs}. If we number
the points as indicated in Fig. \ref{fig:Fano Plane asymmetrical and symmetrical numbered}
(left) we have:\renewcommand{\arraystretch}{1.5}
\[
\begin{array}{c}
P=\left\{ 0,1,2,3,4,5,6\right\} \\
C=\left\{ \left\{ 0,1,2\right\} ,\left\{ 0,3,4\right\} ,\left\{ 0,5,6\right\} ,\left\{ 1,3,5\right\} ,\left\{ 1,4,6\right\} ,\left\{ 2,3,6\right\} ,\left\{ 2,4,5\right\} \right\} 
\end{array}
\]
\renewcommand{\arraystretch}{1}
\end{example}
\noindent 

\noindent It is easy to verify, just by making the corresponding changes
in the asterisms in $C$, that any of the following permutations takes
any asterism in $FP=(P,C)$ to another one, i.e.: that it is an automorphism
of $FP$:\\
\smallskip{}

\noindent \footnotesize
\[
\left(\begin{array}{ccccccc}
0 & 1 & 2 & 3 & 4 & 5 & 6\\
0 & 1 & 2 & 3 & 4 & 5 & 6
\end{array}\right),\left(\begin{array}{ccccccc}
0 & 1 & 2 & 3 & 4 & 5 & 6\\
2 & 0 & 1 & 4 & 5 & 3 & 6
\end{array}\right),\left(\begin{array}{ccccccc}
0 & 1 & 2 & 3 & 4 & 5 & 6\\
1 & 2 & 0 & 5 & 3 & 4 & 6
\end{array}\right)
\]
\medskip{}

\noindent 
\[
\left(\begin{array}{ccccccc}
0 & 1 & 2 & 3 & 4 & 5 & 6\\
0 & 2 & 1 & 4 & 3 & 5 & 6
\end{array}\right),\left(\begin{array}{ccccccc}
0 & 1 & 2 & 3 & 4 & 5 & 6\\
1 & 0 & 2 & 3 & 5 & 4 & 6
\end{array}\right),\left(\begin{array}{ccccccc}
0 & 1 & 2 & 3 & 4 & 5 & 6\\
2 & 1 & 0 & 5 & 4 & 3 & 6
\end{array}\right)
\]
\normalsize\medskip{}

\noindent It is also easy to see that the six permutations form a
group, with composition of permutations as the group operation. We
notice how the group only leaves invariant point $6$. We wonder if
we could link this group to the action of a symmetry group. The obvious
candidate is $D_{3}$, because it has six elements (three rotations
and three reflections) and fixes a single point. 

We could try now to rearrange the points so that they exhibit $D_{3}$
symmetry, and permute them until $D_{3}$ acts on the new asterisms
as prescribed. Fig. \ref{fig:Fano Plane asymmetrical and symmetrical numbered}
(right) shows that this is possible. The action homomorphism here
maps the identity element in $D_{3}$ to the identity permutation
above; the $120\text{º}$ counter-clockwise rotation to the top center
permutation; the $240\text{º}$ counter-clockwise rotation to the
top right permutation; and each of the $3$ reflections to one of
the bottom permutations. For the Fano Plane we have then (see Sect.
\ref{sub:Group-Hierarchy}):

\[
A_{D_{3}}(FP)\subset Aut(FP)\subset Sym(P)
\]
\\
with respective orders 6, 168 and $7!=5040$. An example of an automorphism
not linked to $D_{3}$ is:\smallskip{}
\[
\left(\begin{array}{ccccccc}
0 & 1 & 2 & 3 & 4 & 5 & 6\\
2 & 5 & 4 & 3 & 6 & 1 & 0
\end{array}\right)
\]

\subsubsection{Sperner Families of Boards}

Let's suppose that $S_{SB}$ is the set of symmetric boards that are
solutions of the BSP for a given asymmetric board. If we put their
symmetry groups in a set $S_{SG}$ we can remove duplicate symmetry
groups. To further reduce redundancy in $S_{SB}$, we can retain instead
a subset $S'_{SB}$ of the former such that $S'_{SG}$, the set of
its symmetry groups, is the \emph{Sperner family} \cite{Combinatorics of Finite Sets}
of $S_{SG}$\footnote{Besides the identity, groups in the Sperner family may have other
common elements.}, i.e: a set with two properties:
\begin{itemize}
\item \textit{\emph{no group in }}$S'_{SG}$ \textit{\emph{is a proper subgroup
of another group in }}$S'_{SG}$ 
\item any group in $S_{SG}$ not in $S'_{SG}$ is a subgroup of at least
one member of $S'_{SG}$
\end{itemize}

\subsubsection{Computational Complexity of the BSP}

\noindent There are methods based on geometric symmetry to find graph
automorphisms\footnote{I.e.: permutations of the vertices that leave the edges invariant.}
\cite{Area requirement and symmetry display of planar upward drawings,Drawing graphs symmetrically in three dimensions},
a problem known to be in the class NP of computational complexity
\cite{Some NP-complete problems similar to graph isomorphism}. Since
graphs are a subset of hypergraphs, and since boards are geometric
hypergraphs, the BSP complexity class is also at least NP. For algebraic
methods to find automorphisms of hypergraphs see \cite{full automorphism groups of edge-transitive hypergraphs}\emph{.}

\section{Latin Polytopes\label{sub:Latin-Polytopes}}

In Sect. \ref{sec:Latin-Boards} we defined Latin boards and described
some of their properties. In this section we focus on a subset of
Latin boards with special properties of symmetry.

\subsection{Symmetric Latin Boards}
\begin{defn}
A Latin board is \emph{symmetric} under a geometric transformation
if it is transformed by this one into another Latin board.\end{defn}
\begin{example}
\noindent A Latin board is trivially symmetric under the identity
transformation. A Latin square is symmetric under a reflection across
its vertical symmetry axis.
\end{example}

\begin{example}
The board in Fig. \ref{fig:FanoPlaneUnsymmetricAndSymmetric} (center)
is symmetric under a $120\text{º}$ counter-clockwise rotation. The
board on its left is not symmetric under the same transformation.
\end{example}

\subsection{Definition of Latin Polytope}

The transformations under which a Latin board is symmetric may form
a group. When the group is the symmetry group of a finite regular
polytope we have the following objects:
\begin{defn}
A \emph{Latin finite regular polytope} is a Latin board with the symmetry
of a finite regular polytope.
\end{defn}
\noindent We will normally drop the adjectives and use simply ``Latin
polytopes'' when there is no risk of ambiguity. Latin squares for
example are symmetric under any of the symmetries of the square \textendash a
finite regular polytope with $D_{4}$ symmetry\textendash{} hence
Latin squares are Latin polytopes. The Latin boards in Figs. \ref{fig:Latin triangle example}
and \ref{fig:Latin hexagon example} are also Latin polytopes, since
they are symmetric under $D_{3}$ and $D_{6}$ respectively.

\subsection{Latin Polygons}
\begin{defn}
A \textit{Latin polygon }is a Latin polytope with the symmetry of
a finite regular polygon. 
\end{defn}
\noindent The symmetry groups of the regular polygons are the dihedral
groups $D_{n}$, with order 2\textit{n} (\textit{n} reflections and
\textit{n} rotations, see Sect. \ref{sub:Symmetry-Groups-from-finite-regular-polytopes}).
Some dihedral groups properly contain others, for example $D_{3}\subset D_{6}$.
When a Latin polygon has more than one symmetry group we will name
it after the polygon with the group of highest order.
\begin{example}
\noindent Conventional Latin squares have \textit{D}\textsubscript{4}
symmetry, so they are Latin polygons.\emph{ Sudokus} (see Fig. \ref{fig:Sudoku}
right) are also Latin polygons since they feature the same symmetry.
\end{example}

\subsubsection{About the Example Latin Polygons\label{sub:About-the-Example-Latin polygons}}

The next section defines some Latin polygons. The accompanying examples
are latinizations (see Sect. \ref{sub:Unifit-Latin-Boards}) of the
boards in Sect. \ref{sub:Biregular-polytopes}. The multiset used
in each latinization defines the balance and intersection pattern
in the resulting Latin polygon and its woven board (see Sect. \ref{sec:Woven-boards}).
The examples were created by the author using the technique described
in \cite{Latin Puzzles paper}.

\subsubsection{Latin Triangles}
\begin{defn}
A \textit{Latin triangle }is a Latin polytope with $D_{3}$ symmetry. \end{defn}
\begin{example}
\label{exa:If-we-latinize}If we latinize the board in Fig. \ref{fig:raw la tri order 6 faces}
with numbers $1$ to $12$ we obtain Latin triangles like the one
in Fig. \ref{fig:latin triangle faces Small} (left)\footnote{This example corresponds to the notion of Latin triangle defined in
\cite{Latin triangles}, which is a more restricted concept than the
one proposed here. See Fig. \ref{fig:hardLatinTetrahedron} for a
larger Latin triangle with points in the barycenters, like this one.}. To derive its woven board we make up a warp asterism out of the
points with the same label. Each warp asterism intersects in one point
each weft asterism. Warp asterisms constitute then an orthogonal class:
the resulting woven board has SIN $\left\{ 0,1,4\right\} $.
\end{example}

\begin{example}
Fig. \ref{fig:latin triangle faces Small} (right) shows a Latin triangle
with the board in Fig. \ref{fig:raw la tri order 6 faces} and the
multiset $\left\{ 1,1,1,1,2,2,2,2,3,3,3,3\right\} $. The SIN in the
derived woven board is $\left\{ 0,4\right\} $, as every warp asterism
intersects in four points every weft one. Comparing these results
with those in Example \ref{exa:If-we-latinize} we notice the role
of the multiset in the balance and intersection pattern of the resulting
boards.
\end{example}

\begin{example}
Fig. \ref{fig:latin tri order 7 vtces 9Lx4PSmall} (left) shows a
Latin triangle with the board in Fig. \ref{fig:raw la tri order 7 vtcesSmall}
latinized with numbers $1$ to $9$, i.e.: with an orthogonal warp
class in the derived woven board. 
\end{example}

\begin{example}
Fig. \ref{fig:latin tri order 7 vtces 9Lx4PSmall} (right) shows a
Latin triangle resulting from the latinization of the board in Fig.
\ref{fig:raw la tri order 7 vtcesSmall} with the multiset $\left\{ 1,1,1,2,2,2,3,3,3\right\} $.
Each warp asterism in the derived woven board intersects in four points
every weft asterism.
\end{example}

\begin{example}
Fig. \ref{fig:latin tri order 4 edges 3Lx10P} (left) shows a Latin
triangle resulting from the latinization of the board in Fig. \ref{fig:raw lat tri order 4 edges 3x10P}
with numbers $0$ to $9$, i.e.: with an orthogonal warp class in
its woven board. 
\end{example}

\begin{example}
Fig. \ref{fig:latin tri order 4 edges 3Lx10P} (right) shows a Latin
triangle with the board in Fig.\ref{fig:raw lat tri order 4 edges 3x10P}
and the multiset $\left\{ 0,0,1,1,2,2,3,3,4,4\right\} $. Here each
warp asterism in the woven board intersects in four points every weft
one.
\begin{figure}
\begin{centering}
\includegraphics[height=5cm]{latin_triangle_faces_Small}~~~~~~~~~~\includegraphics[height=5cm]{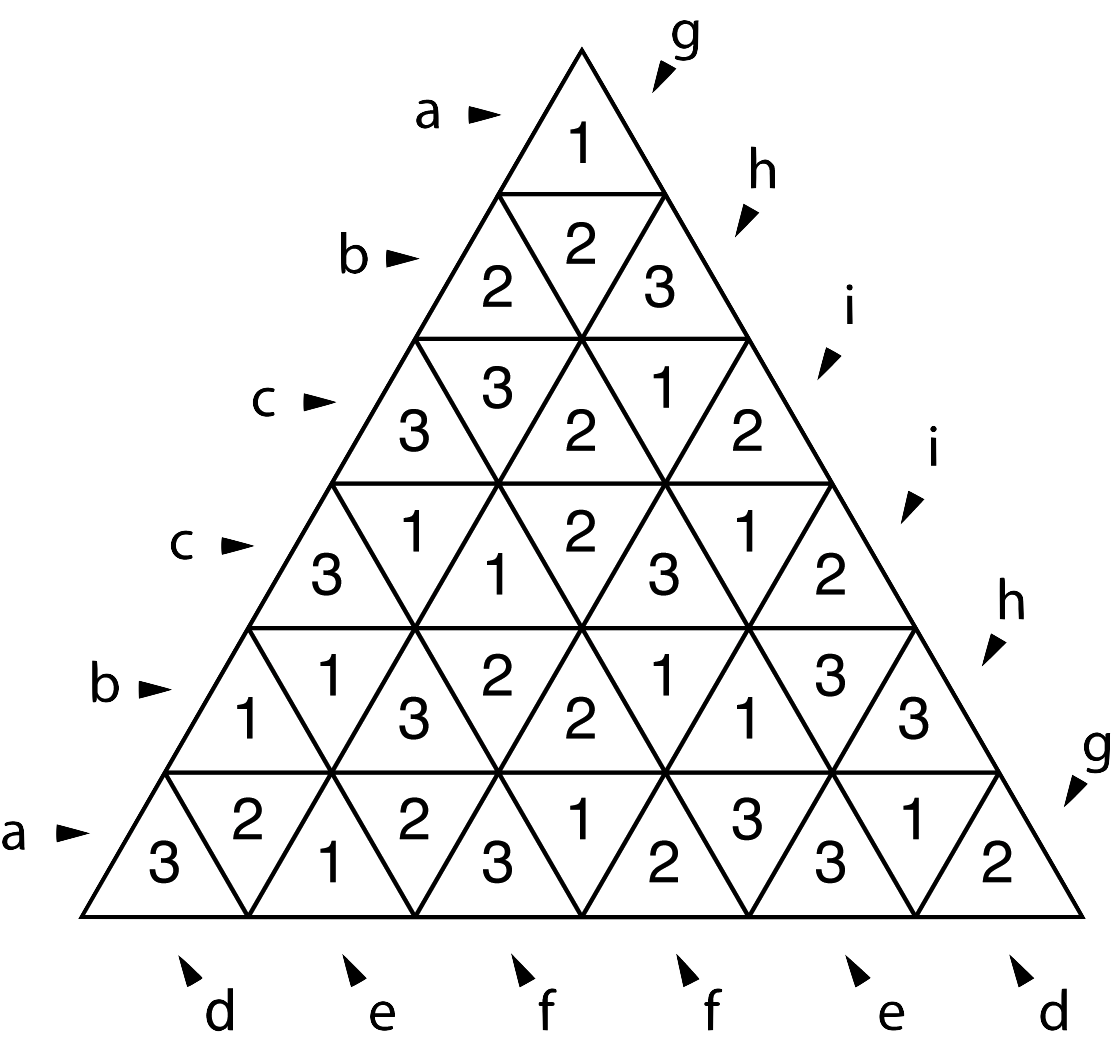}
\par\end{centering}

\protect\caption{\label{fig:latin triangle faces Small}Latin triangles with points
at the barycenters}
\end{figure}
\begin{figure}
\begin{centering}
\includegraphics[height=5.2cm]{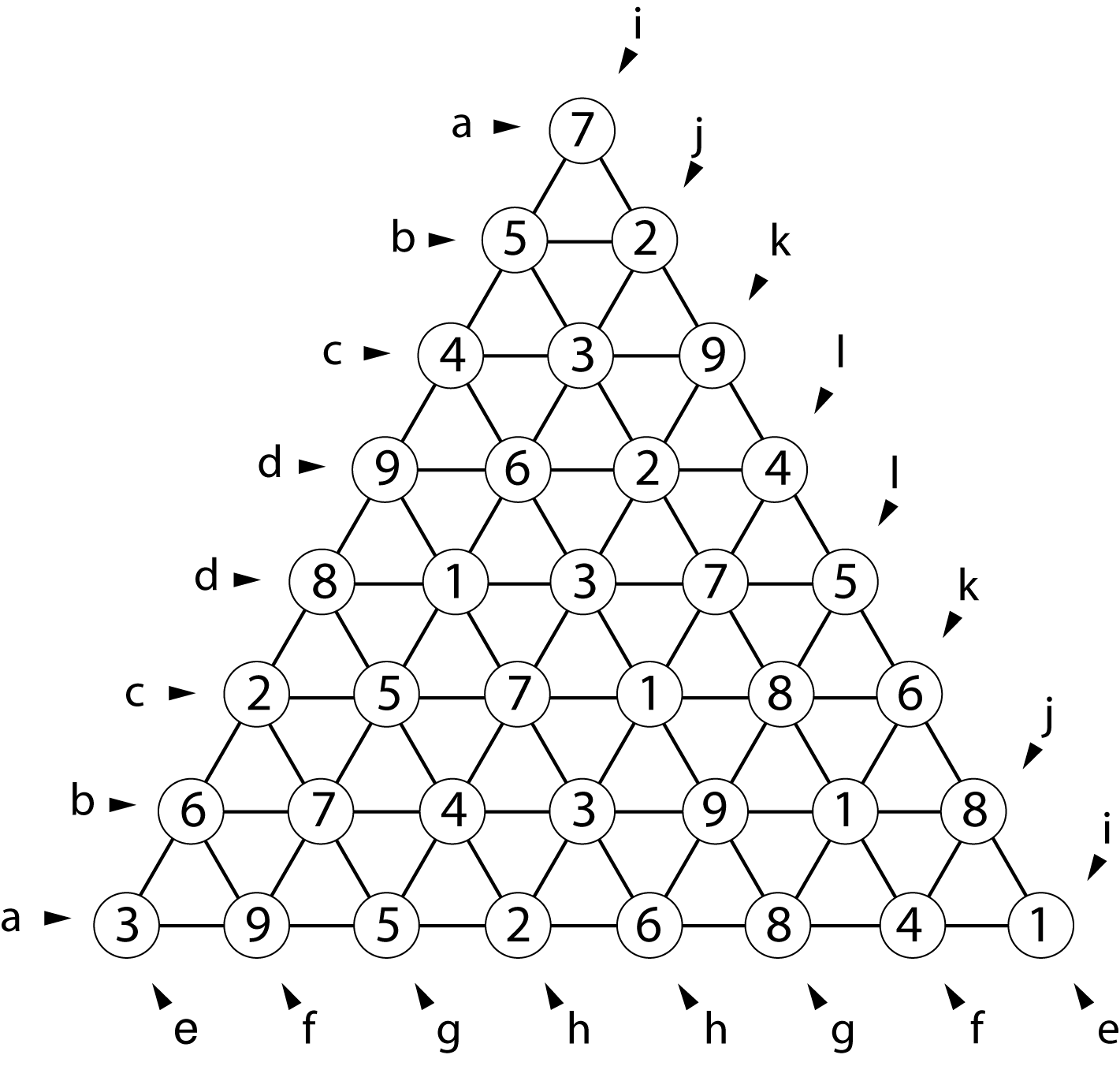}~~~~~\includegraphics[height=5.2cm]{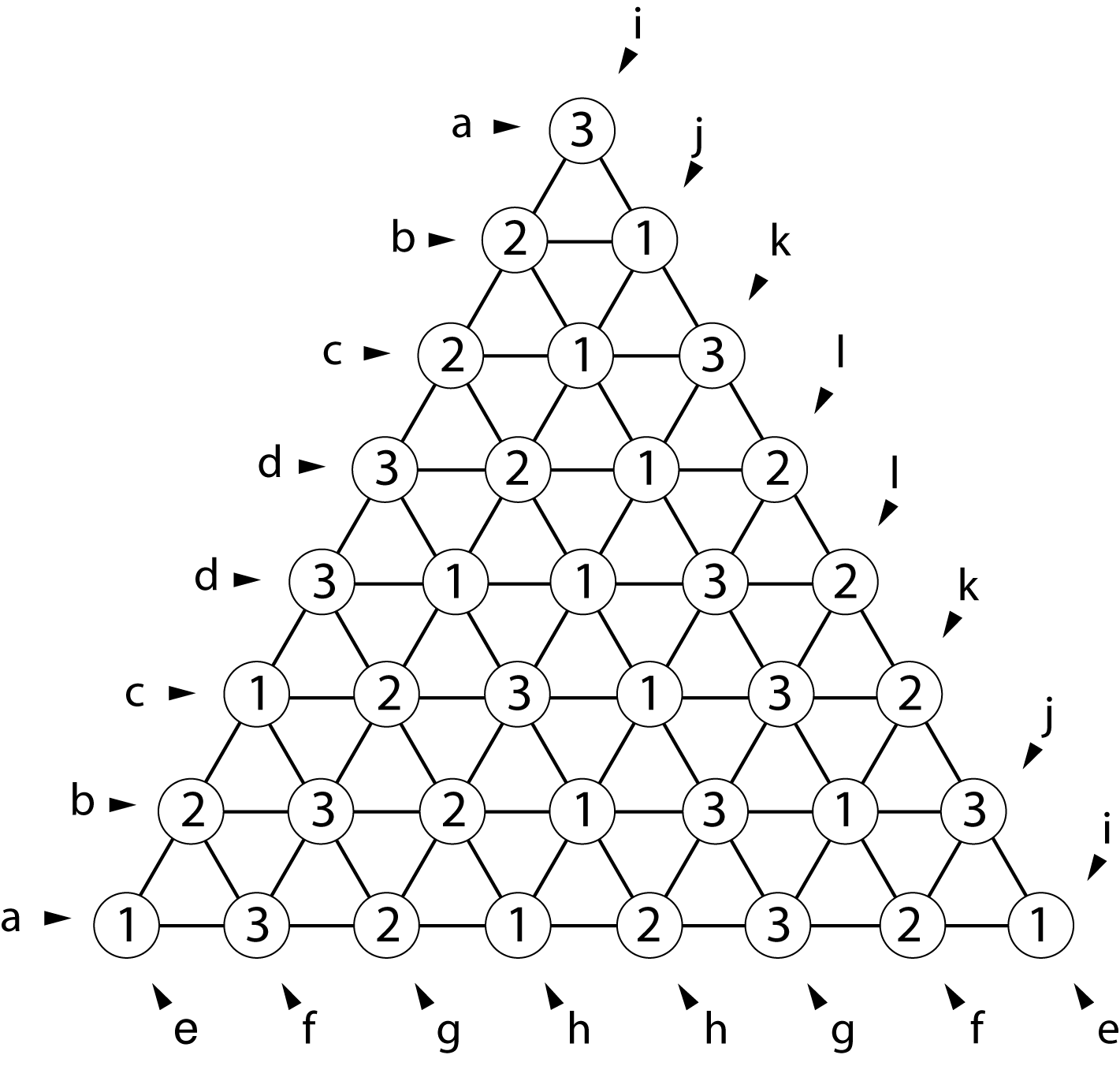}
\par\end{centering}

\protect\caption{\label{fig:latin tri order 7 vtces 9Lx4PSmall}Latin triangles with
points at the vertices}
\end{figure}
\begin{figure}
\begin{centering}
\includegraphics[height=5.5cm]{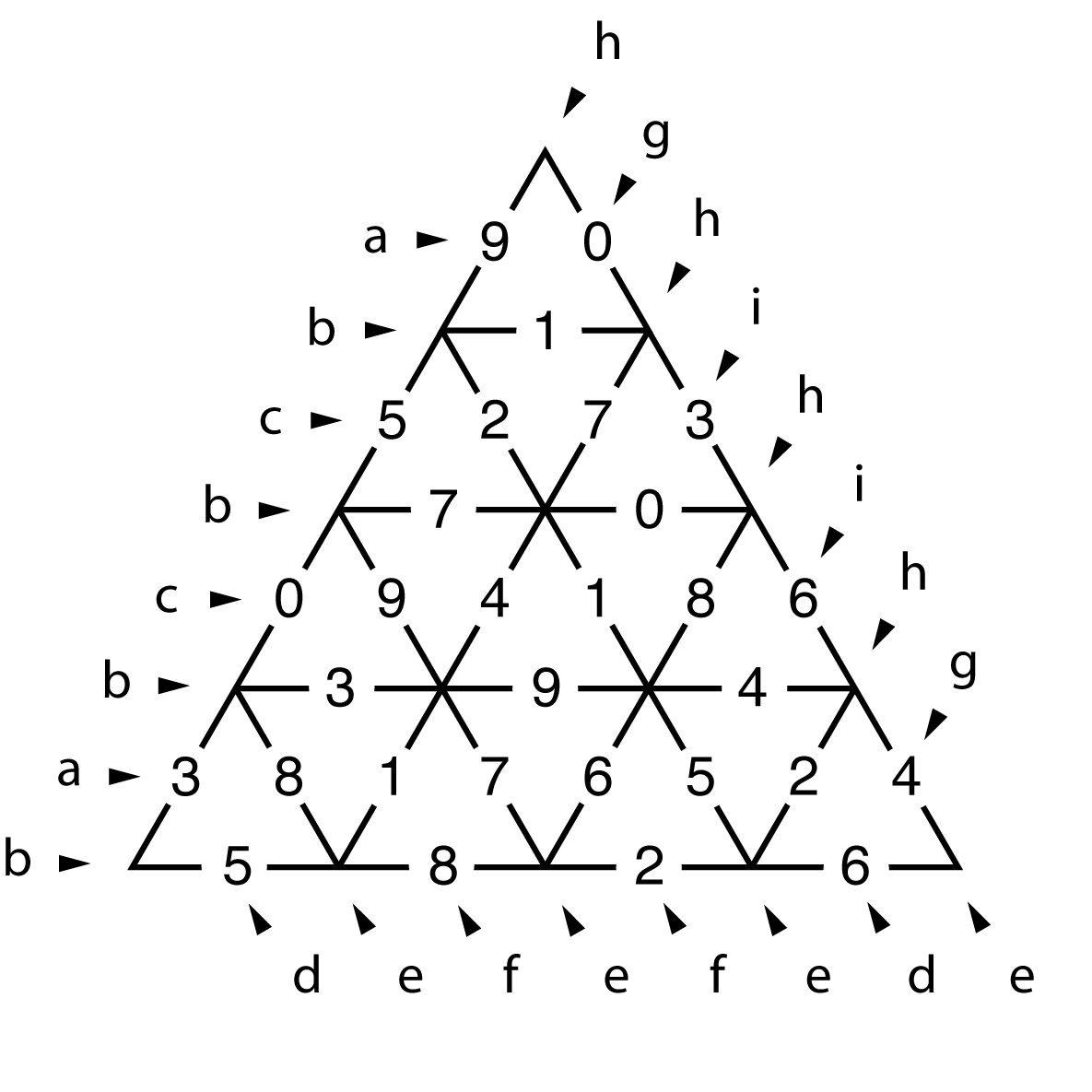}~~~~~\includegraphics[height=5.7cm]{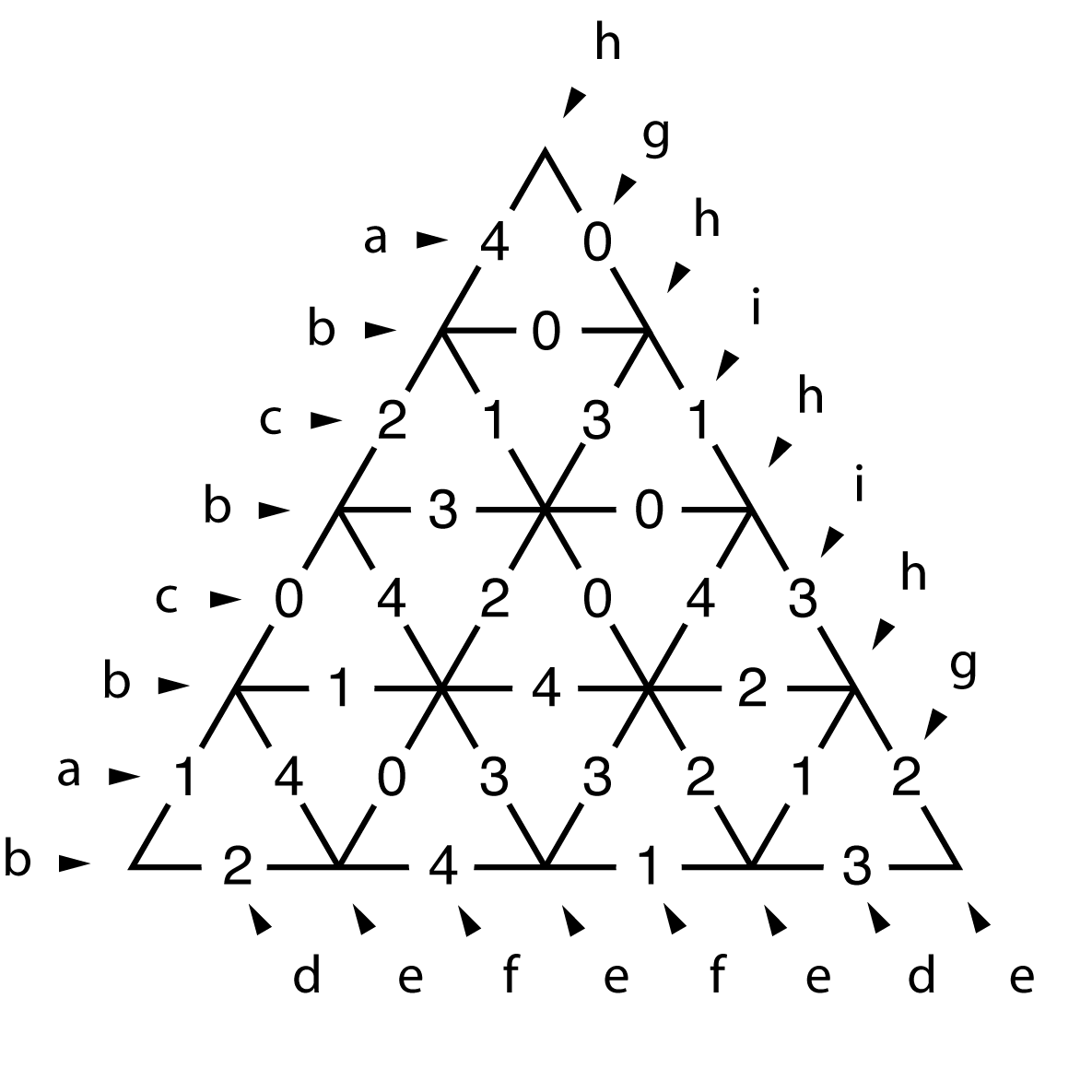}
\par\end{centering}

\protect\caption{\label{fig:latin tri order 4 edges 3Lx10P}Latin triangles with points
at the edge midpoints}
\end{figure}

\end{example}

\subsubsection{Free Latin Squares\label{sub:Free-Latin-Squares}}
\begin{defn}
A \textit{free Latin square }is a Latin polytope with $D_{4}$ symmetry. \end{defn}
\begin{example}
Fig. \ref{fig:free Latin square example} shows a free Latin square
resulting from the latinization of the board in Fig. \ref{fig:board Free Latin Square}
with all numbers from $1$ to $16$. Its derived woven board has SIN
$\left\{ 0,1,4\right\} $ because its warp asterisms constitute an
orthogonal class.
\end{example}

\begin{example}
Fig. \ref{fig:Latin sequence of Free Latin squares} shows three free
Latin square with the board in Fig. \ref{fig:board Free Latin Square}
and the multiset $\left\{ 1,1,1,1,2,2,2,2,3,3,3,3,4,4,4,4\right\} $.
Their derived woven boards have also SIN $\left\{ 0,4\right\} $ because
every warp asterism intersects in four points every weft asterism.
\end{example}

\begin{example}
\label{exa:Tartan board}Fig. \ref{fig:tartan board} shows a free
Latin square with the board in Fig. \ref{fig:tartan board-1} and
numbers $1$ to $16$. The warp class in the derived woven board is
orthogonal to every weft class, so its SIN is $\left\{ 0,1,4\right\} $.
\begin{figure}
\centering{}\includegraphics[height=4.7cm]{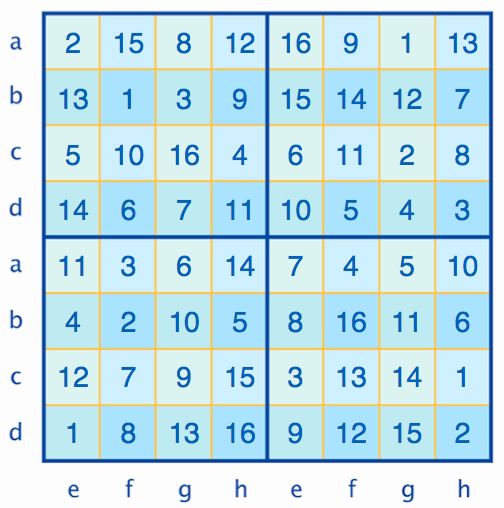}\protect\caption{\label{fig:tartan board}Free Latin square}
\end{figure}

\end{example}

\begin{example}
The base board for \emph{Sudoku} (see Fig. \ref{fig:sudoku design pts-1})
is clearly $D_{4}$ symmetric although not class-symmetric, because
there is no way to take the rows or columns to the $3\times3$ sub-square
asterisms with an element of the group. Latinizations of this board
with different multisets, like the ones shown on the right in Figs.
\ref{fig:Sudoku}, \ref{fig:Sudoku Ripeto} and \ref{fig:Custom Sudoku}
are also free Latin squares.
\end{example}

\begin{figure}
\begin{centering}
\includegraphics[height=5cm]{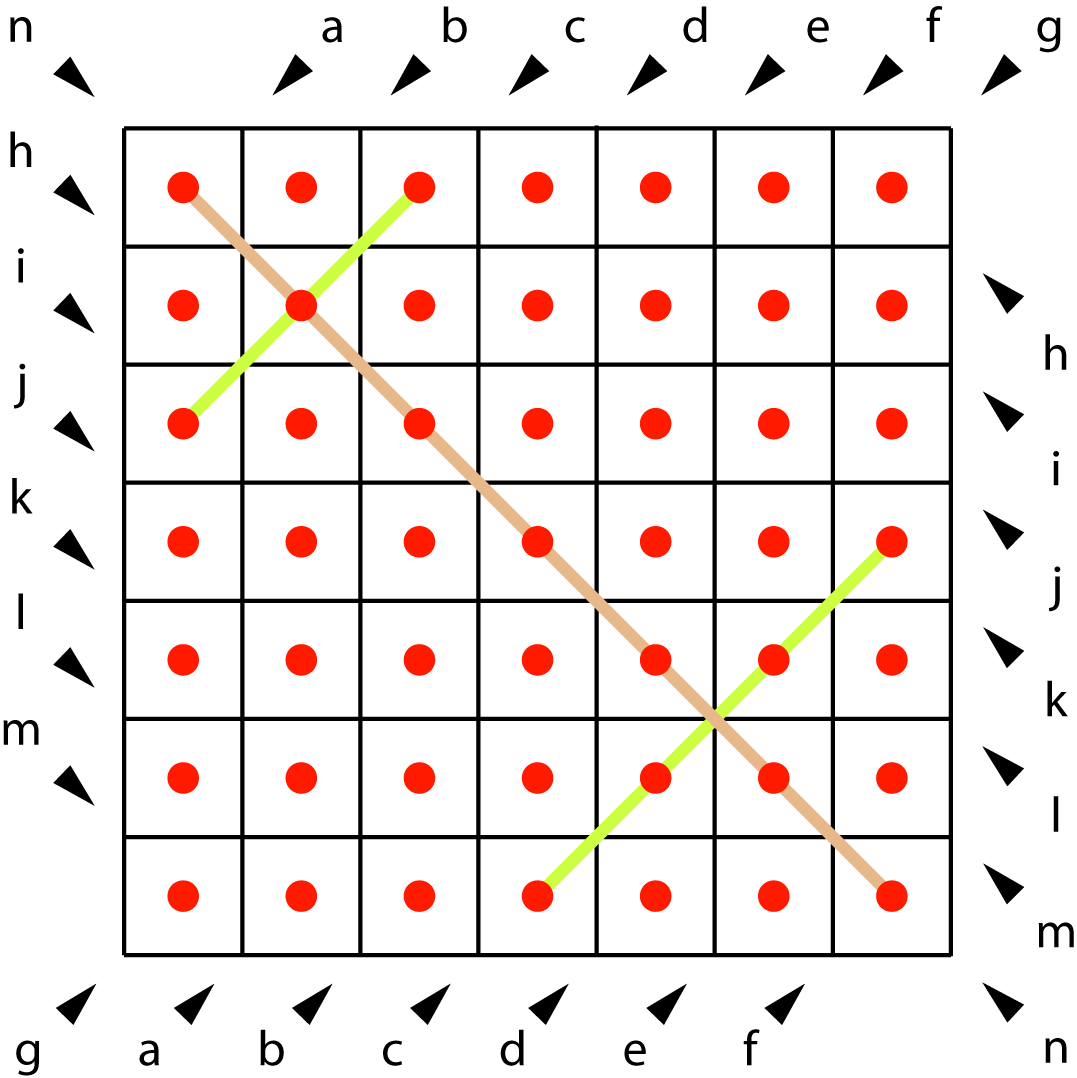}~~~~~~~~~~\includegraphics[height=5cm]{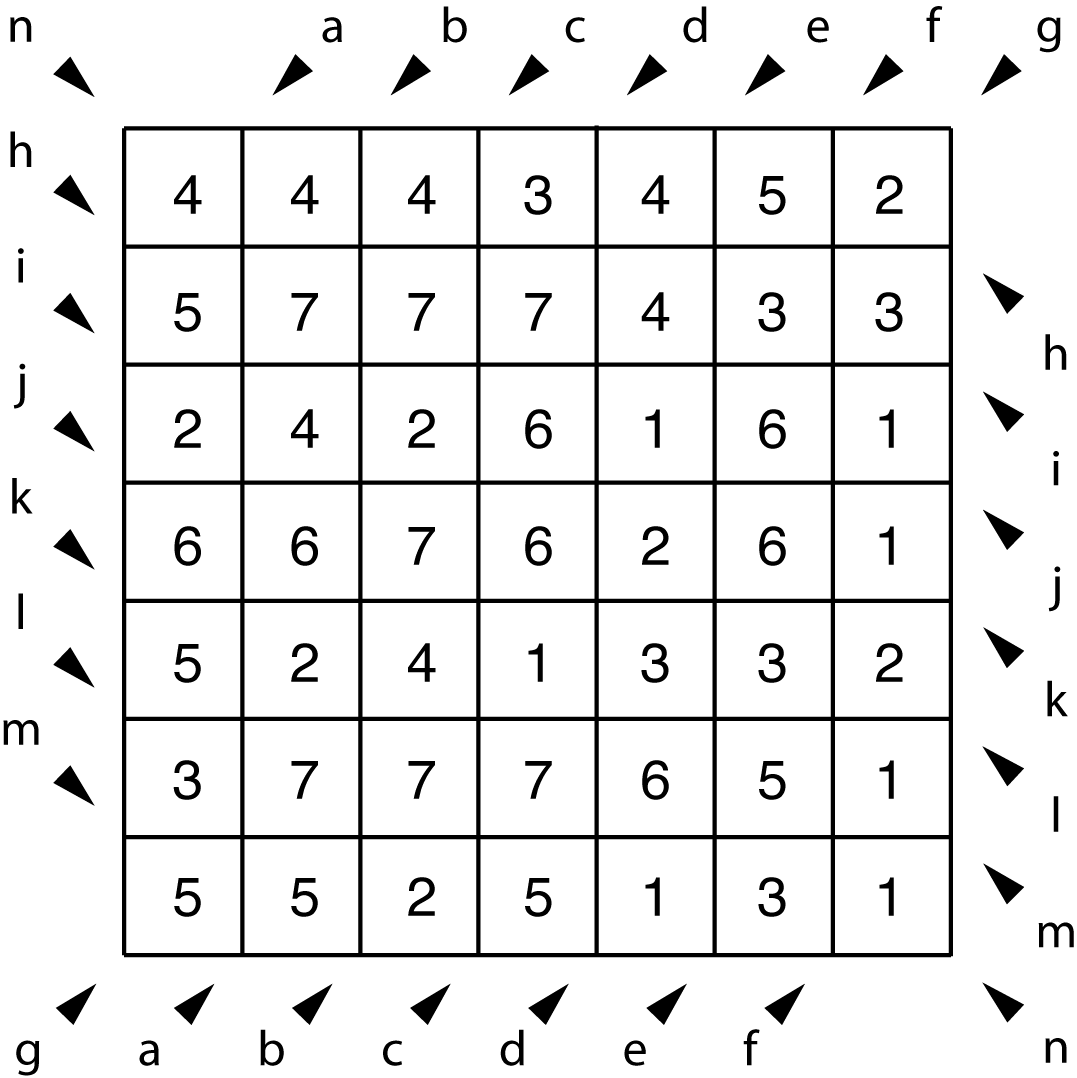}
\par\end{centering}

\protect\caption{\label{fig:slantedSolBoardNbrs}Square board with SIN $\left\{ 0,1\right\} $
and free Latin square}
\end{figure}

\begin{example}
A conventional Latin square is a Latin polygon having\end{example}
\begin{enumerate}
\item \noindent a biregular square (see Sect. \ref{exa:Boards-from-Birregular-squares})
as a source
\item \noindent face centers as board points
\item \noindent a class-symmetric board with \textit{$D_{4}$} symmetry
and SIN $\left\{ 0,1\right\} $
\item a derived woven board with SIN $\left\{ 0,1\right\} $ too
\end{enumerate}
\noindent This description does not univocally determines conventional
Latin squares, as we see next. Fig. \ref{fig:slantedSolBoardNbrs}
(left) shows a biregular square of order seven whose face centers
have been chosen as board points. Fourteen asterisms of seven points
each have been formed, each one coming from either one or two slanted
sets of points pointed to by the same letter (for example asterism
$c$ in green and asterism $n$ in brown). There are two parallel
classes here: $\left\{ a,b,c,d,e,f,g\right\} $ and $\left\{ h,i,j,k,l,m,n\right\} $.
As \textit{D}\textsubscript{4} acts transitively on the classes the
board is class-symmetric, and has SIN $\left\{ 0,1\right\} $. Fig.
\ref{fig:slantedSolBoardNbrs} (right) shows a free Latin square with
this board labeled with numbers $1$ to $7$. Its derived woven board
has SIN $\left\{ 0,1\right\} $ because the warp asterisms constitute
an orthogonal class. This free Latin square has then the main features
of conventional Latin squares without being itself one of them.

\begin{example}
Fig. \ref{fig:punctSlantSquare} shows a free Latin square with the
board in Fig. \ref{fig:punctSquareBoard} latinized with numbers $1$
to $8$. The SIN in the derived woven board is $\left\{ 0,1,2\right\} $
because each warp asterism is orthogonal to every weft one. This example
shows that the number of points in a free Latin square need not be
a perfect square.
\begin{figure}
\begin{centering}
\includegraphics[height=5cm]{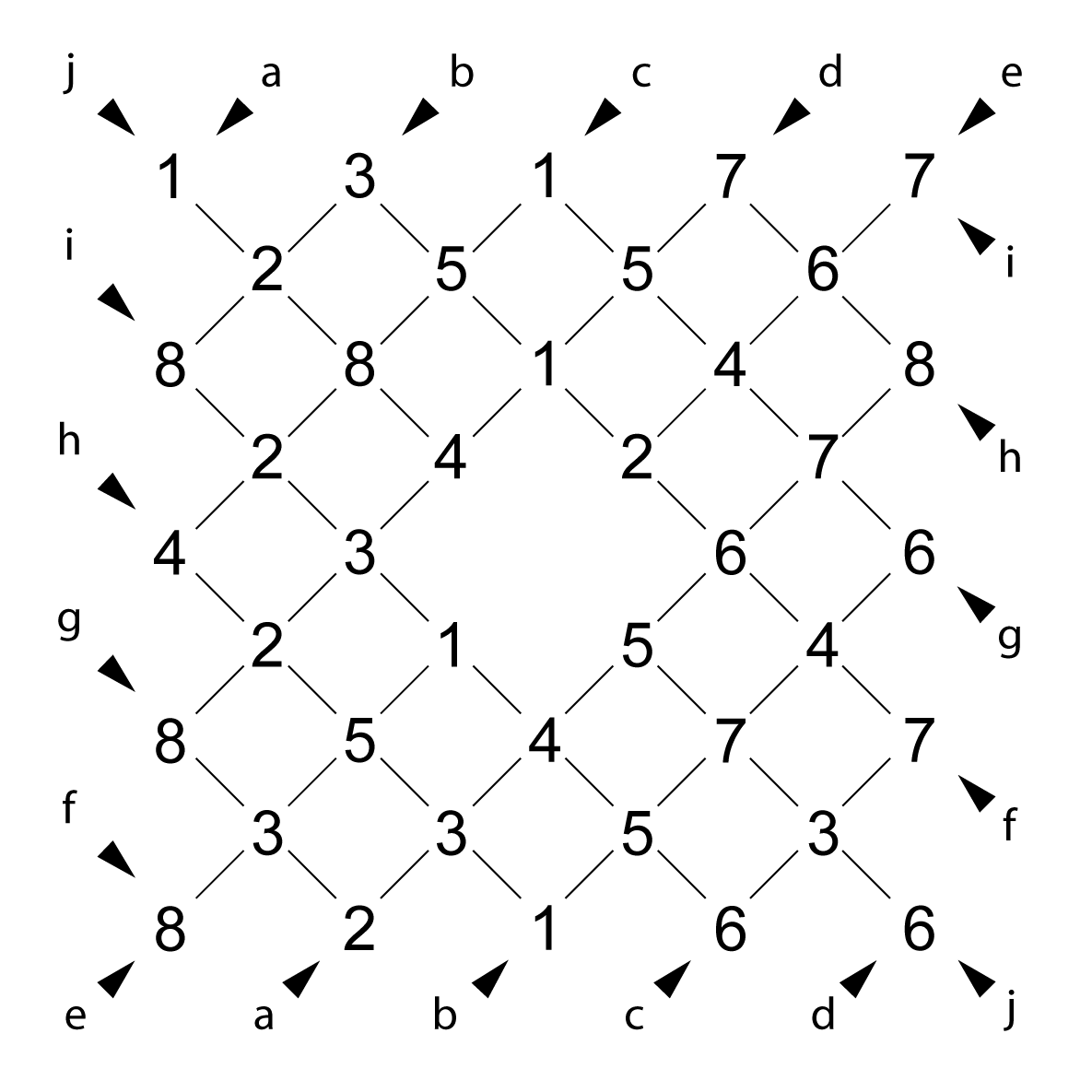}
\par\end{centering}

\protect\caption{\label{fig:punctSlantSquare}Free Latin square}
\end{figure}

\end{example}

\subsubsection{Latin Hexagons}
\begin{defn}
A \textit{Latin hexagon }is a Latin polytope with $D_{6}$ symmetry. \end{defn}
\begin{example}
\noindent Fig. \ref{fig:latinHexFacesOrthoOnSide} (left) shows a
Latin hexagon featuring the board in Fig. \ref{fig:hexFacesSymmetricBoard}
labeled with numbers $1$ to $18$. The derived woven board has SIN
$\left\{ 0,1,6,8\right\} $ since the warp asterisms constitute an
orthogonal class here.
\end{example}

\begin{example}
\noindent Fig. \ref{fig:latinHexFacesOrthoOnSide} (right) shows a
Latin hexagon with the board in Fig. \ref{fig:hexFacesSymmetricBoard}
and the multiset $\left\{ 1,1,1,1,1,1,2,2,2,2,2,2,3,3,3,3,3,3\right\} $.
This time the derived woven board has SIN $\left\{ 0,6,8\right\} $
because each warp asterism intersects in six points every weft asterism.
\begin{figure}
\begin{centering}
\includegraphics[height=5cm]{latinHexFacesOrthoOnSide}~~~~~~~\includegraphics[height=5cm]{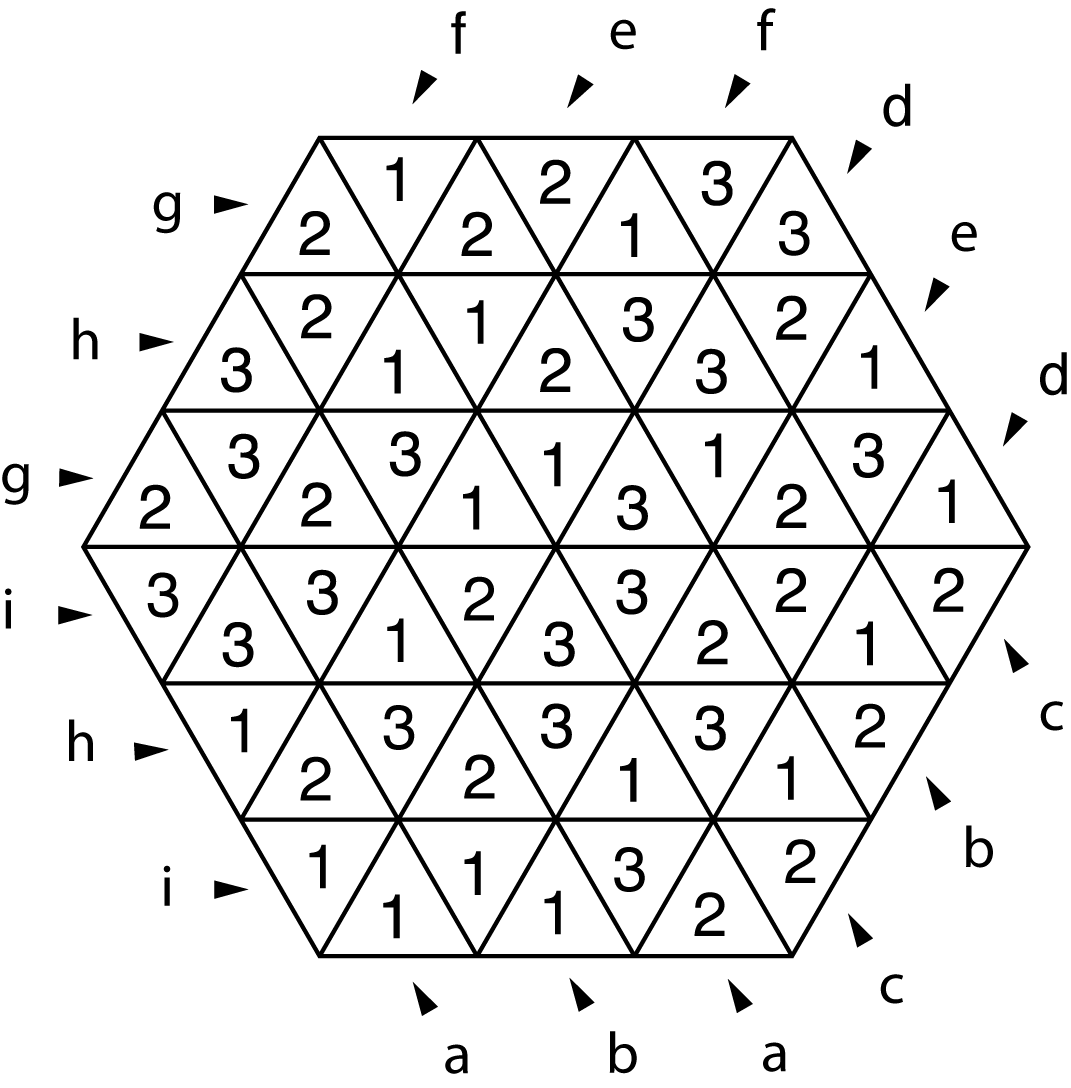}
\par\end{centering}

\protect\caption{\label{fig:latinHexFacesOrthoOnSide}Latin hexagons with points at
the barycenters}
\end{figure}

\end{example}

\begin{example}
Fig. \ref{fig:latinHexEdgesOrthoOnSide} (left) shows a Latin hexagon
with the board in Fig. \ref{fig:hexSymBoardEdgesOnSide} and all numbers
from $1$ to $14$. The derived woven board has an orthogonal warp
class because every warp asterism intersect in one point every weft
asterism.
\end{example}

\begin{example}
Fig. \ref{fig:latinHexEdgesOrthoOnSide} (right) shows another Latin
hexagon with the board in Fig. \ref{fig:hexFacesSymmetricBoard} latinized
this time with multiset $\left\{ 1,1,2,2,3,3,4,4,5,5,6,6,7,7\right\} $.
Here every warp asterism in the derived woven board intersects in
two points every weft asterism.
\begin{figure}
\begin{centering}
\includegraphics[height=5.4cm]{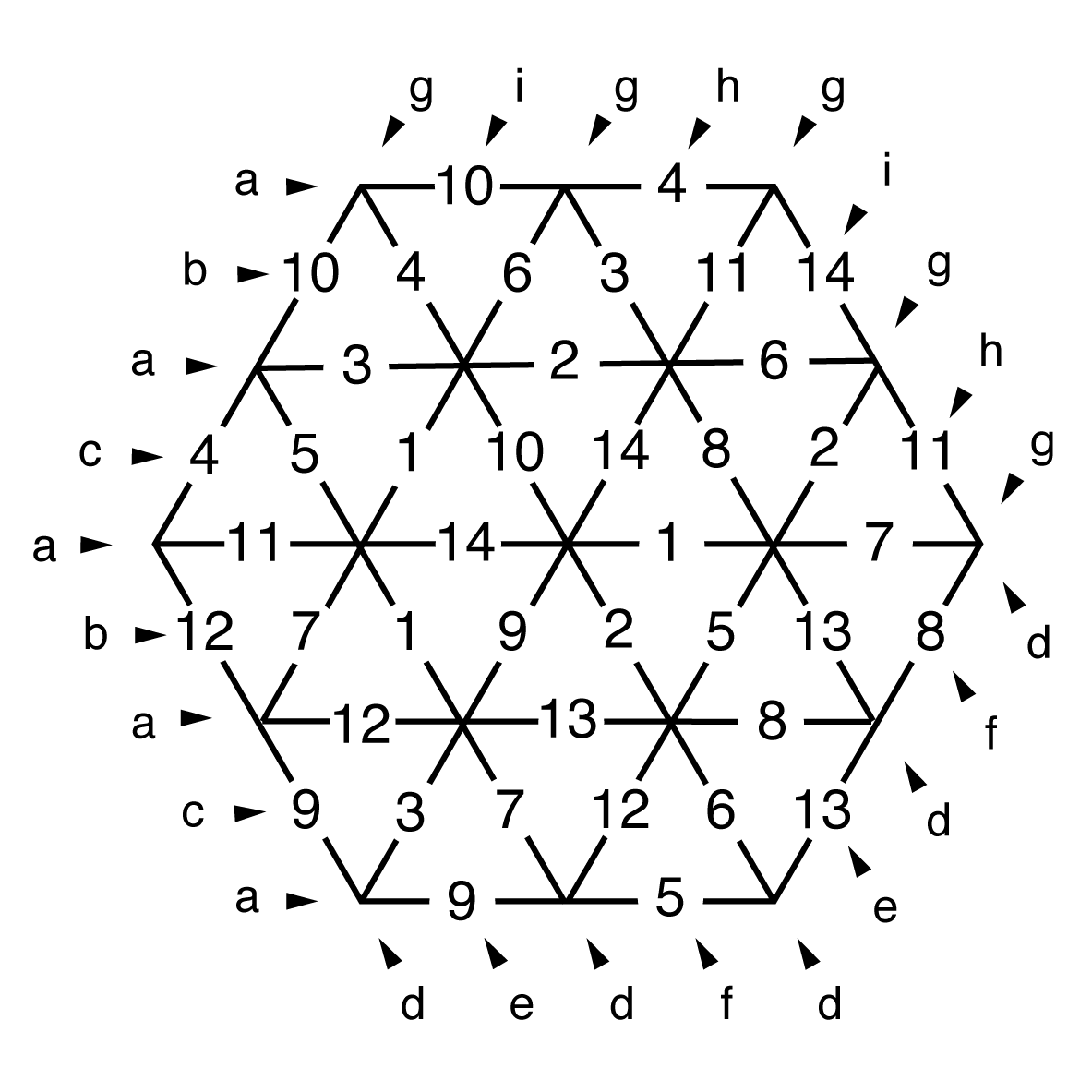}~~~~\includegraphics[height=5.5cm]{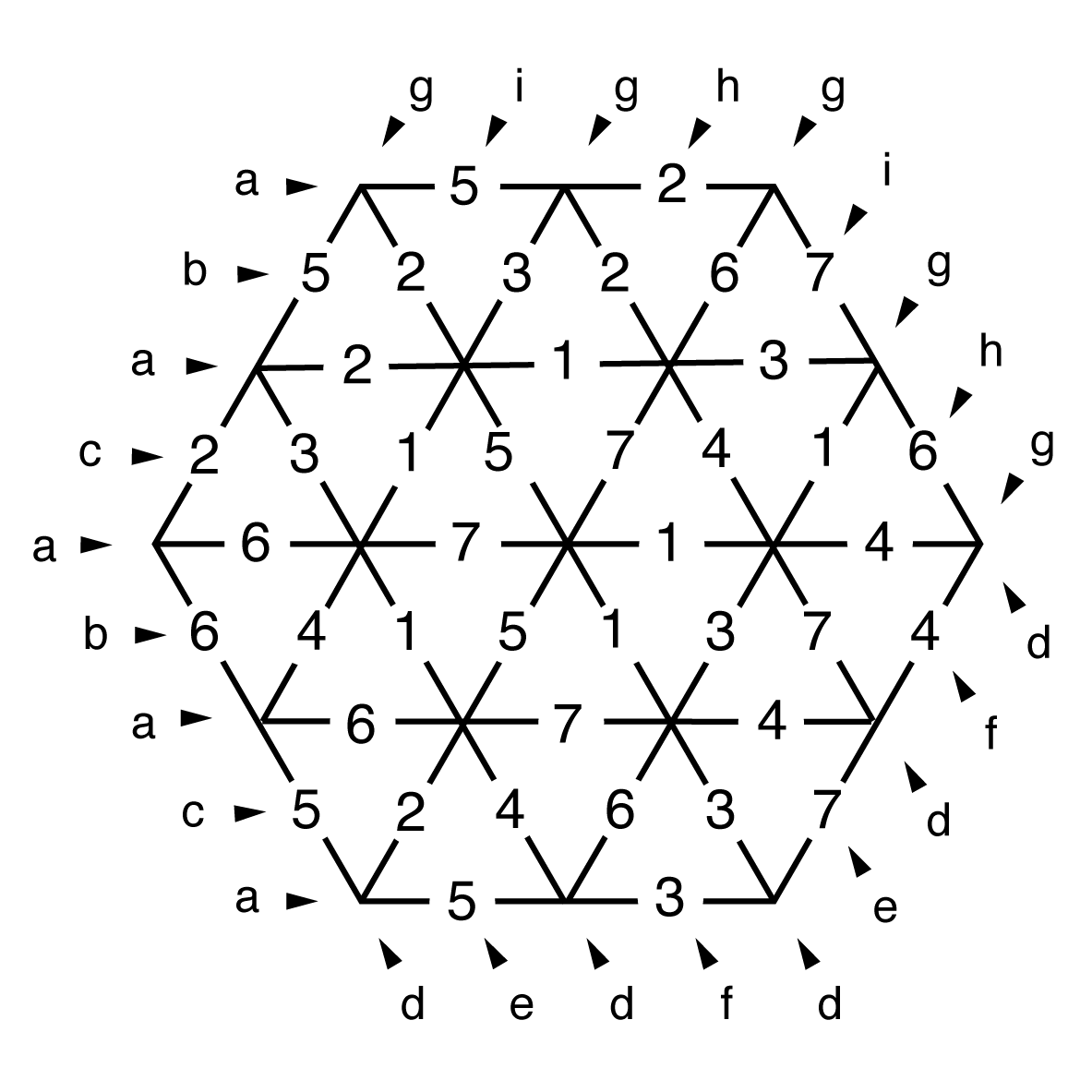}
\par\end{centering}

\protect\caption{\label{fig:latinHexEdgesOrthoOnSide}Latin hexagons with points at
the edge midpoints}
\end{figure}

\end{example}

\subsubsection{Latin Dodecagons}
\begin{defn}
A \textit{Latin dodecagon }is a Latin polytope with $D_{12}$ symmetry. \end{defn}
\begin{example}
\label{exa:Helios board} Fig. \ref{fig:HeliosSmall} shows a\emph{
}Latin dodecagon with the board in Fig. \ref{fig:emptyBrdHelios}
latinized with the set \{H, E, L, I, O, S\}. The SIN of the derived
woven board is $\left\{ 0,1\right\} $ because the warp asterisms
constitute a parallel class orthogonal to the weft asterisms.
\begin{figure}
\begin{centering}
\includegraphics[height=5cm]{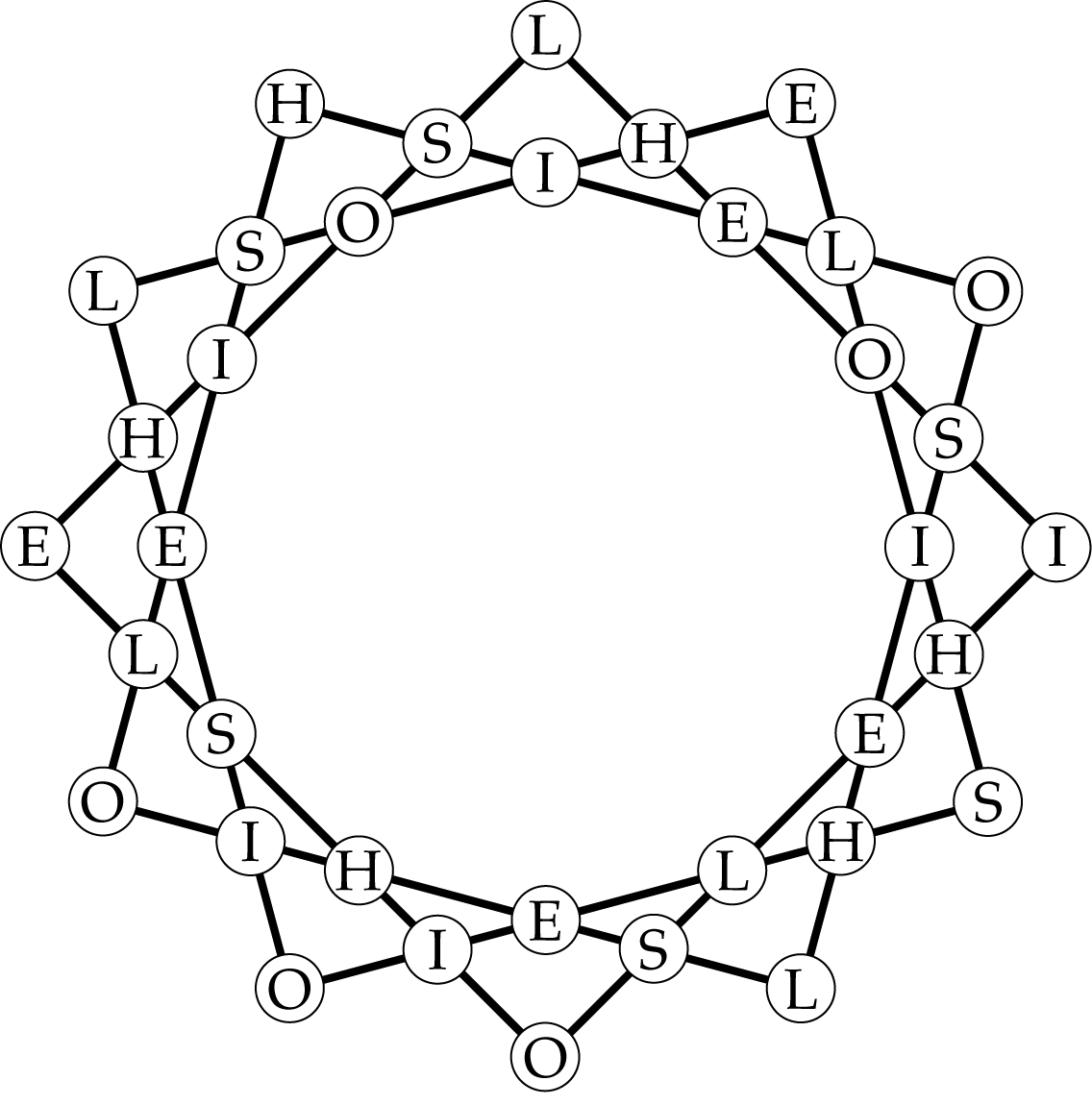}
\par\end{centering}

\protect\caption{\label{fig:HeliosSmall}Latin dodecagon}
\end{figure}

\end{example}

\subsection{Latin Polyhedra\label{sub:Latin-Polyhedra}}
\begin{defn}
A \textit{Latin polyhedron }is a Latin polytope with the symmetry
of a finite regular polyhedron.
\end{defn}
\noindent There are five finite real convex regular polyhedra, and
three polyhedral symmetry groups (see Table \ref{tab:symmetryGroupsPlatonicSolids-1}).
We may name the resulting Latin polyhedron either after its polyhedral
group or after the polyhedron that is closer to the source of symmetry
used. In any case the underlying symmetry group should be clear. Hereafter
we will adopt the second naming convention.

\subsubsection{About the Example Latin Polyhedra}

\noindent Each of the following examples presents a Latin polyhedron
and a corresponding Latin puzzle (see \cite{Latin Puzzles paper}).
All pictures are meant to be folded along the inner black lines and
glued along the outer ones. With the exception of those in the example
Latin dodecahedron, every asterism is made up of points enclosed by
two lines of the same color. The asterism becomes a closed band \textendash so
to speak\textendash{} in the polyhedron surface when the picture is
folded. The examples were created by the author using the latinization
techniques described in \cite{Latin Puzzles paper}.

\subsubsection{Latin Tetrahedra}

\begin{figure}
\begin{centering}
\includegraphics[height=5cm]{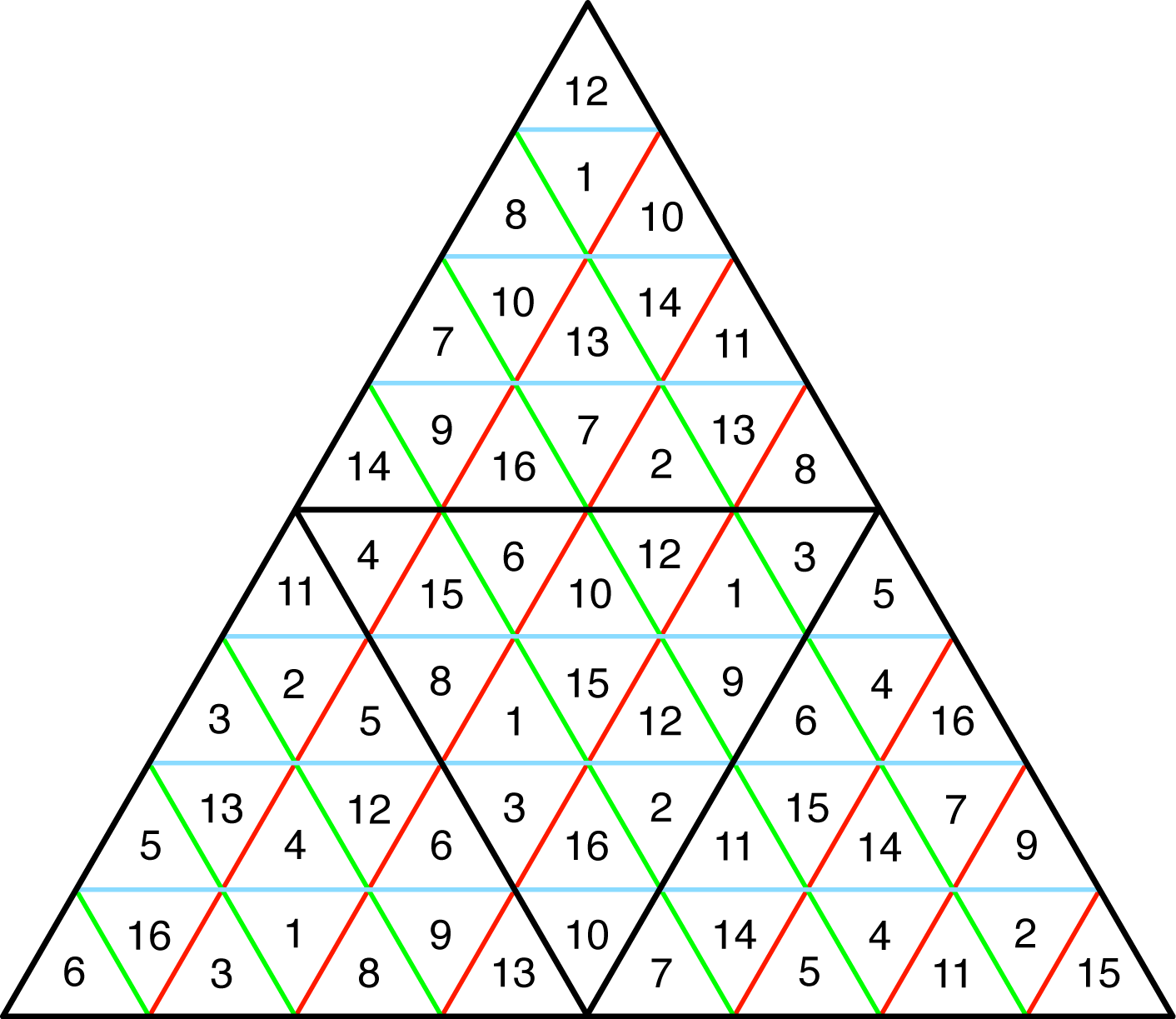}~~~~~~~~~~\includegraphics[height=5cm]{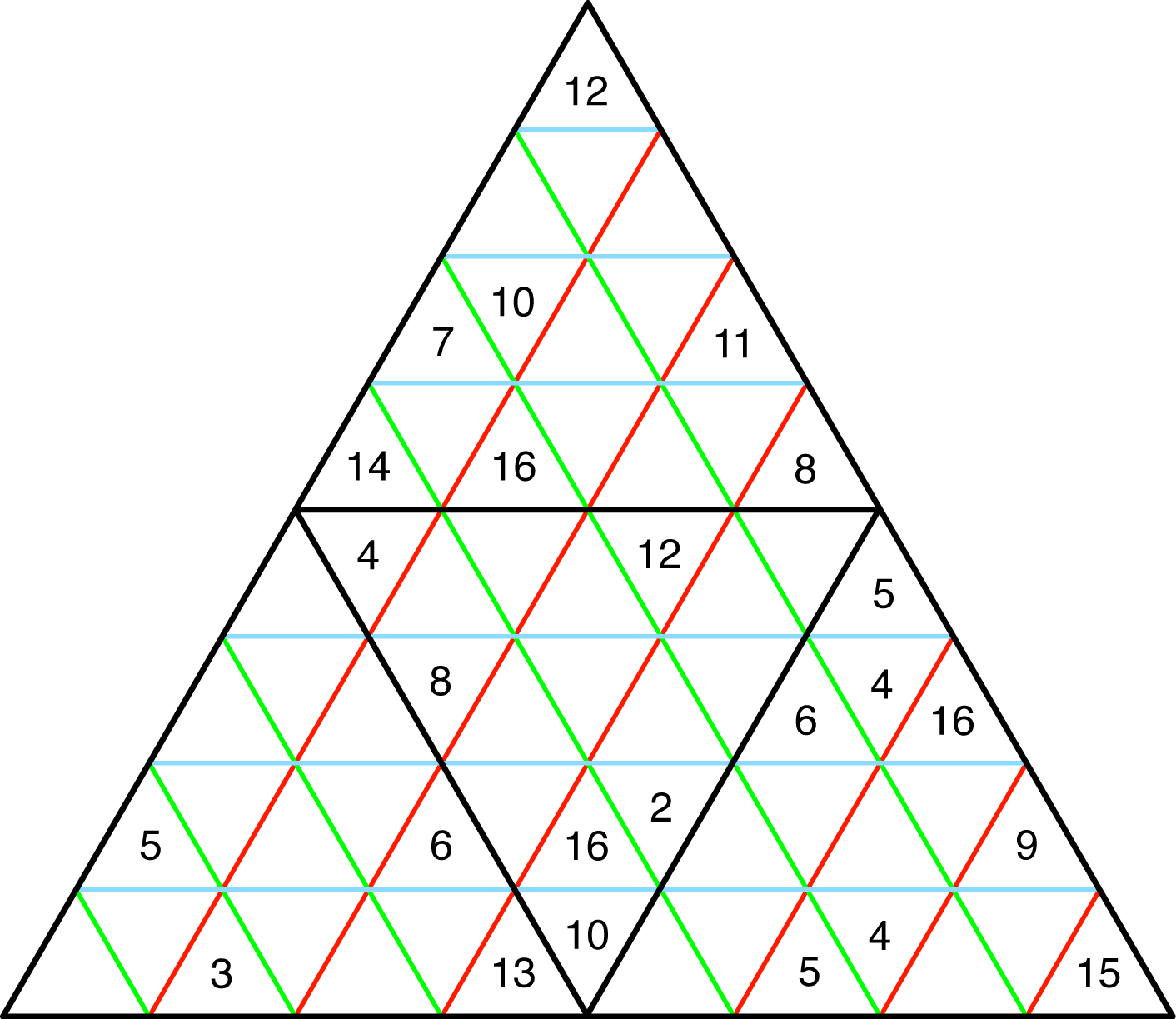}
\par\end{centering}

\protect\caption{\label{fig:hardLatinTetrahedron}Latin tetrahedron and Latin tetrahedron
puzzle}

\end{figure}

Fig. \ref{fig:hardLatinTetrahedron} (left) shows the board in Fig.
\ref{fig:biregularTetrahedron} latinized with all numbers from $1$
to $16$, i.e.: a \emph{Latin tetrahedron}\footnote{If we similarly fold and glue the Latin triangle in Example \ref{exa:latin triangle faces}
we also obtain a (smaller) Latin tetrahedron.}. The derived woven board has SIN $\left\{ 0,1,4\right\} $, since
it has an orthogonal warp class with sixteen asterisms with four points
each. On the right there is a corresponding Latin puzzle.

\subsubsection{Free Latin Cubes}

\begin{figure}
\begin{centering}
\includegraphics[height=8cm]{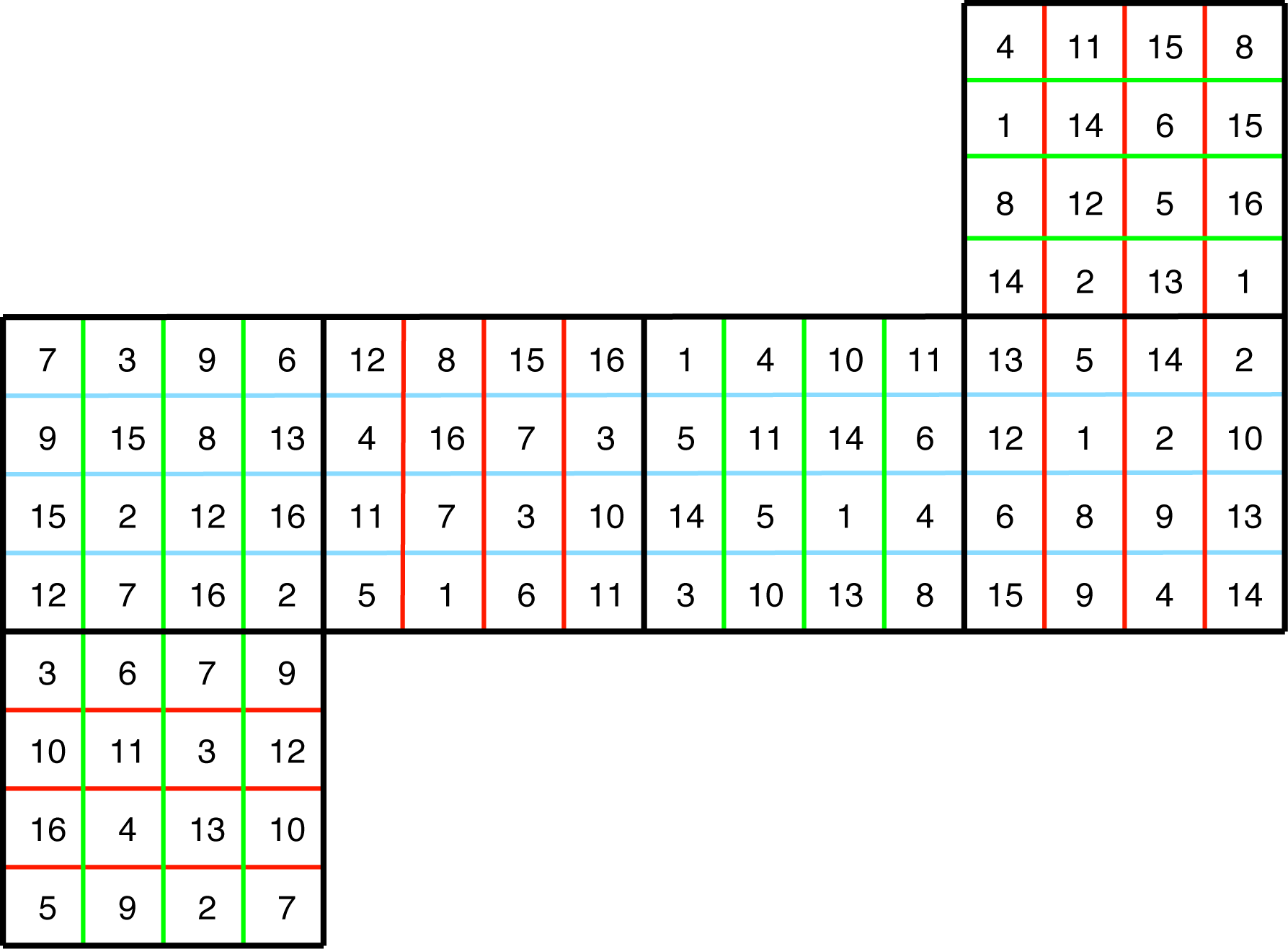}
\par\end{centering}

\protect\caption{\label{fig:latinCube}Free Latin cube}
\end{figure}

\begin{figure}
\begin{centering}
\includegraphics[height=8cm]{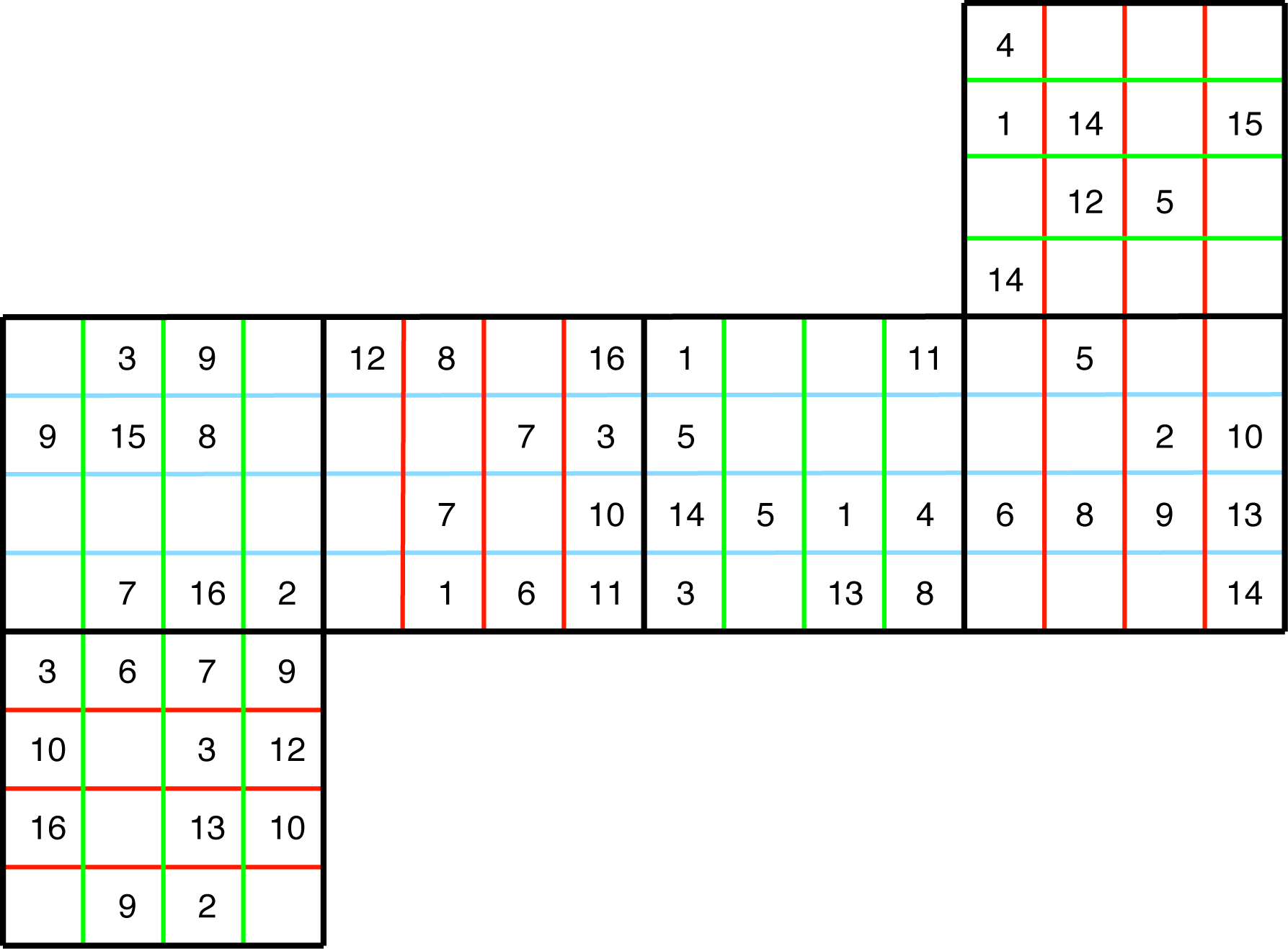}
\par\end{centering}

\protect\caption{\label{fig:partialLatinCube}Free Latin cube puzzle}
\end{figure}

Fig. \ref{fig:latinCube} shows a \emph{free} \emph{Latin cube}\footnote{\noindent Or \emph{Latin octahedron} if we follow the first naming
criterion in Sect. \ref{sub:Latin-Polyhedra}. We use ``free Latin
cube'' instead of ``Latin cube'' because this last name already
identifies two different combinatorial objects. One is the natural
generalization of Latin squares: an $n\times n\times n$ array in
which each of a set of $n$ elements occurs once in each row, once
in each column and once in each file. The other is a different object
used by statisticians. This is why the first meaning cited, the generalization
of Latin squares, is usually called ``permutation cube'' \cite{Latin squares and their applications}.} with the board in Fig. \ref{fig:latin cube symmetric board} latinized
with all numbers from $1$ to $16$. The derived woven board has SIN
$\left\{ 0,1,8\right\} $ because it has an orthogonal warp class
with sixteen asterism with six points each. Fig. \ref{fig:partialLatinCube}
shows a corresponding puzzle.

\subsubsection{Latin Octahedra}

\begin{figure}
\begin{centering}
\includegraphics[height=8cm]{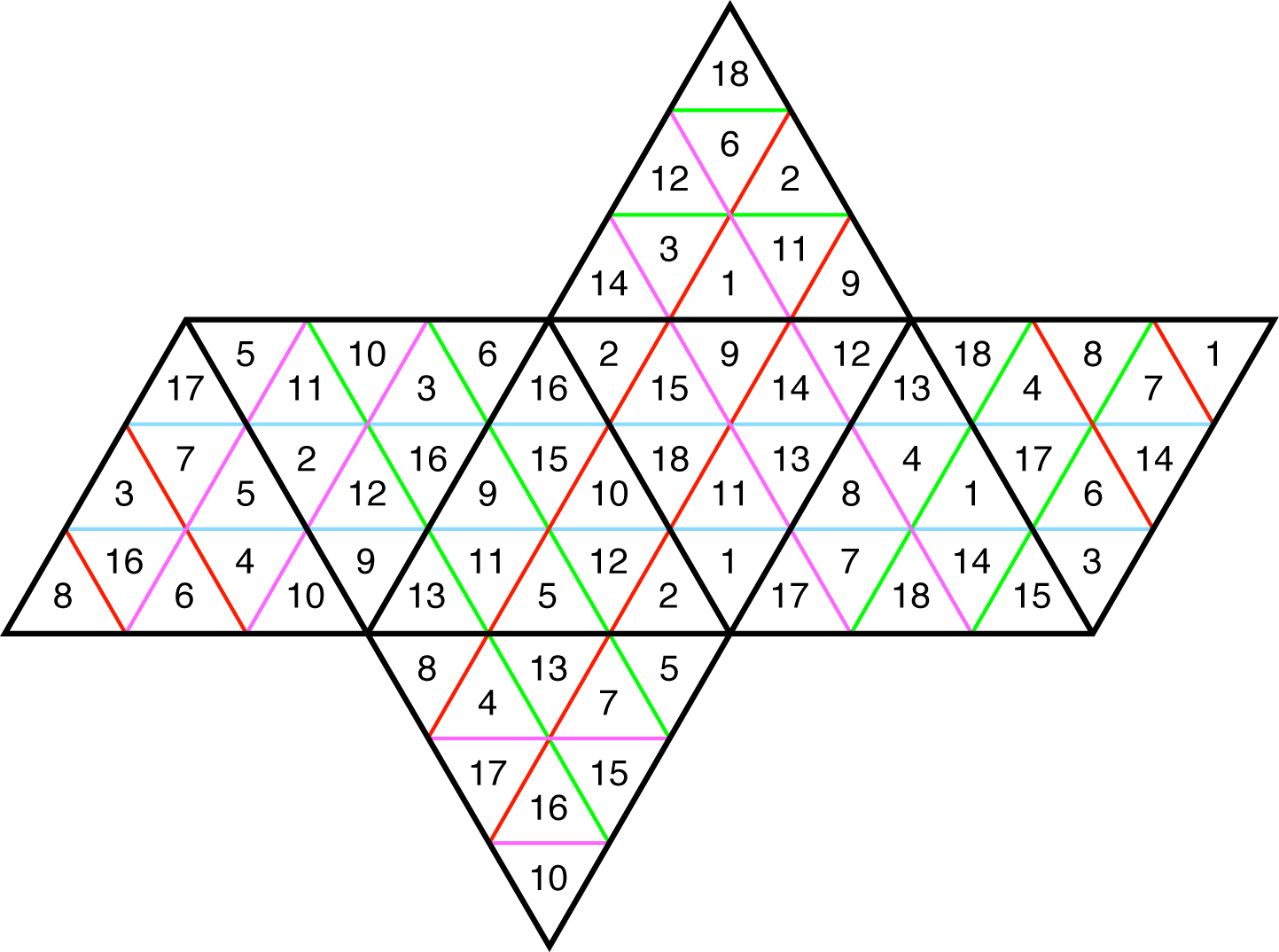}
\par\end{centering}

\protect\caption{\label{fig:latinOctahedron}Latin octahedron}
\end{figure}

\begin{figure}
\begin{centering}
\includegraphics[height=8cm]{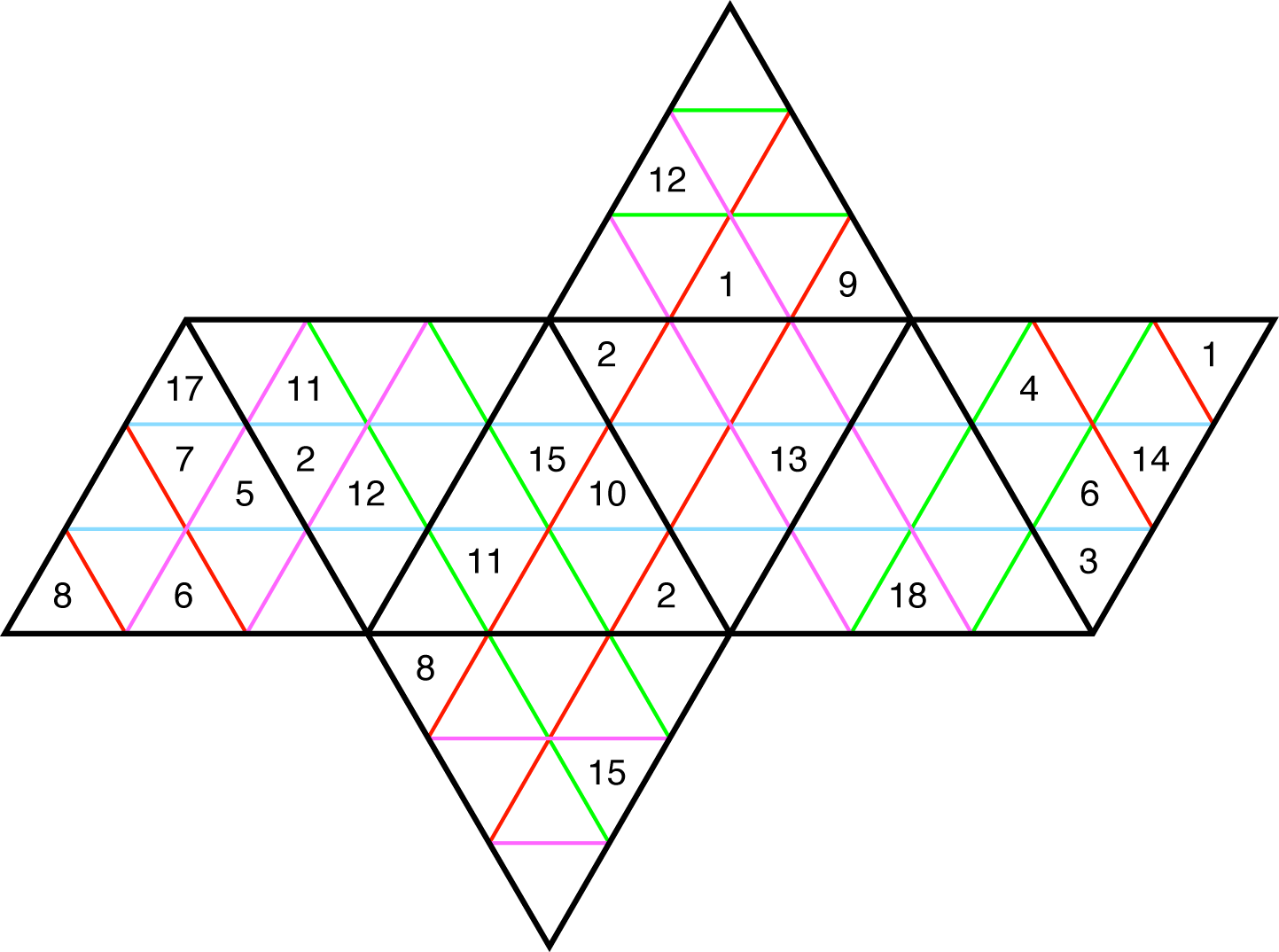}
\par\end{centering}

\protect\caption{\label{fig:partialLatinOctahedron}Latin octahedron puzzle}
\end{figure}

Fig. \ref{fig:latinOctahedron} shows a \emph{Latin octahedron} resulting
from the latinization of the board in Fig. \ref{fig:latin octahedron symmetric board}
with all numbers from $1$ to $18$. The derived woven board has SIN
$\left\{ 0,1,4\right\} $ because it has an orthogonal warp class
with eighteen asterisms with four points each. Fig. \ref{fig:partialLatinOctahedron}
shows a corresponding Latin puzzle.

\subsubsection{Latin Icosahedra}

\begin{figure}
\begin{centering}
\includegraphics[height=5.8cm]{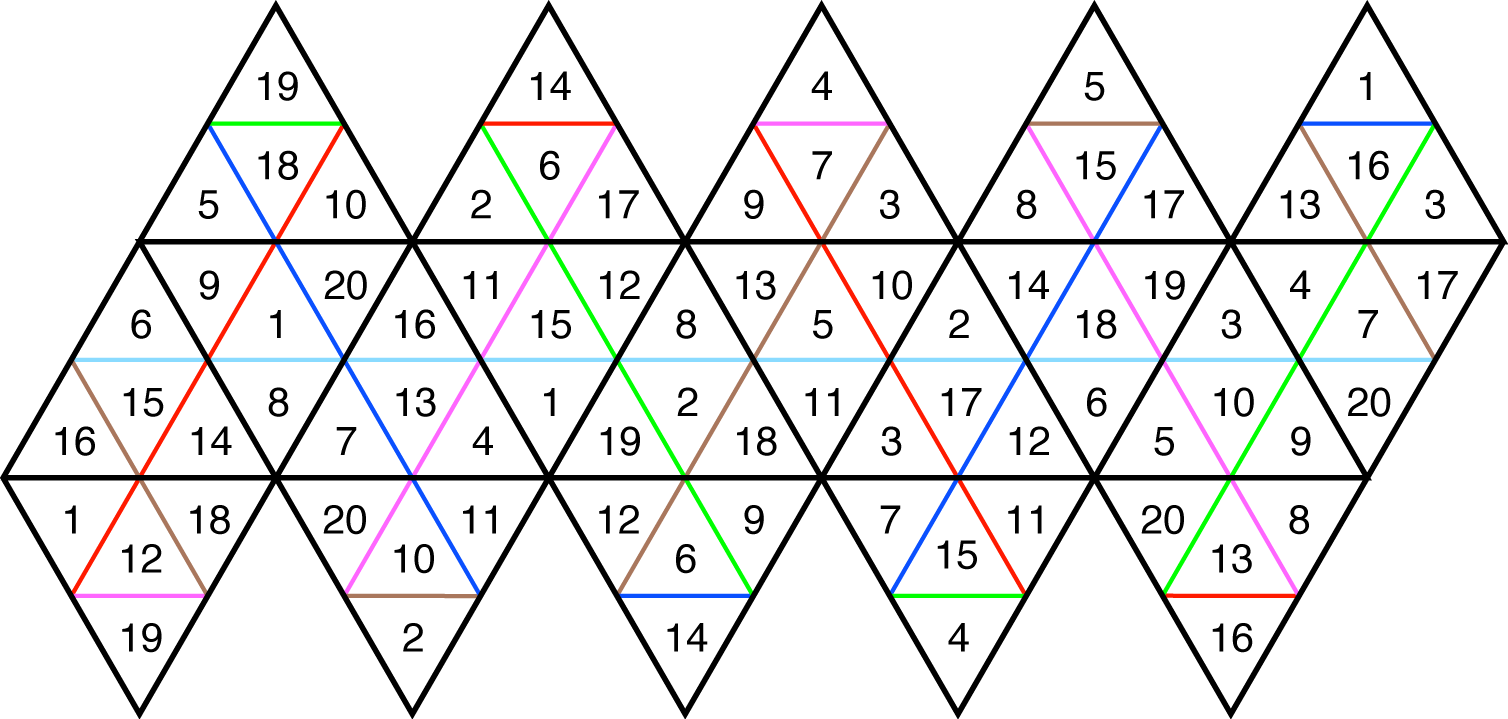}
\par\end{centering}

\protect\caption{\label{fig:latinIcosahedron}Latin icosahedron}
\end{figure}

\begin{figure}
\begin{centering}
\includegraphics[height=5.8cm]{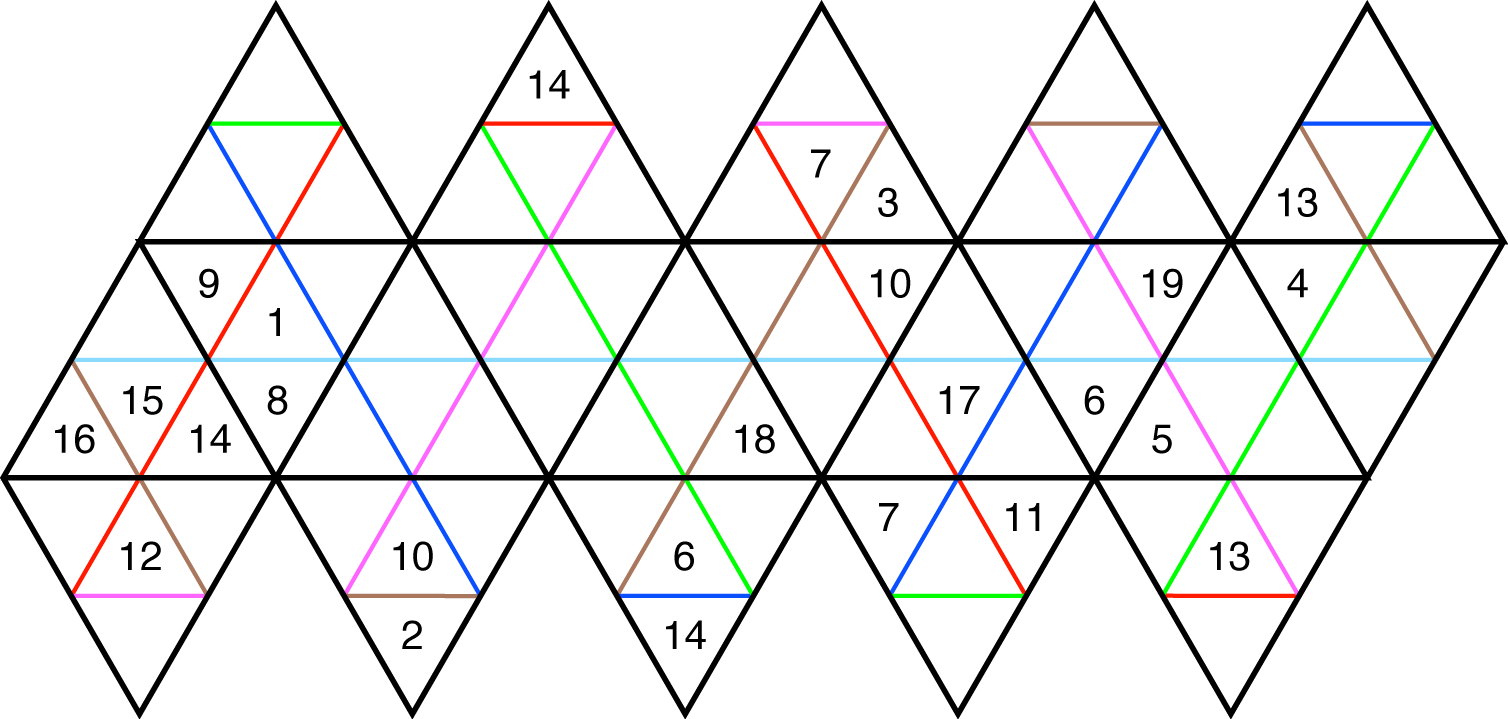}
\par\end{centering}

\protect\caption{\label{fig:partialLatinIcosahedron}Latin icosahedron puzzle}
\end{figure}

Fig. \ref{fig:latinIcosahedron} shows a \emph{Latin icosahedron}
resulting from the latinization of the board in Fig. \ref{fig:latin icosahedron symmetric board}
with all numbers from $1$ to $20$. The derived woven board has SIN
$\left\{ 0,1,4\right\} $ because it has an orthogonal warp class
with twenty asterisms with four points each. Fig. \ref{fig:partialLatinIcosahedron}
shows a corresponding Latin puzzle.

\subsubsection{Latin Dodecahedra}

\begin{figure}
\begin{centering}
\includegraphics[height=6cm]{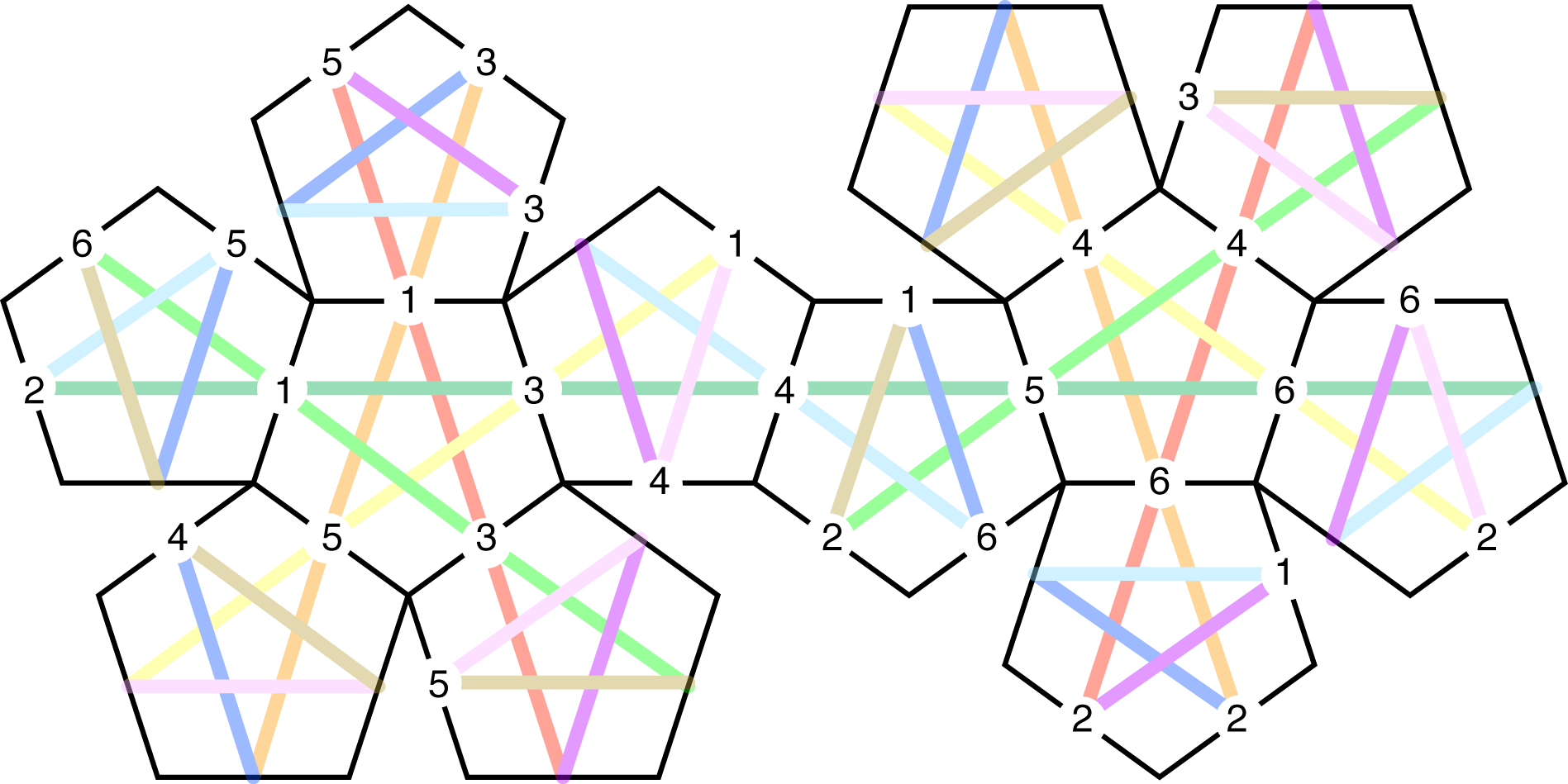}
\par\end{centering}

\protect\caption{\label{fig:hardLatinDodecahedronEdges}Latin dodecahedron}
\end{figure}

\begin{figure}
\begin{centering}
\includegraphics[height=6cm]{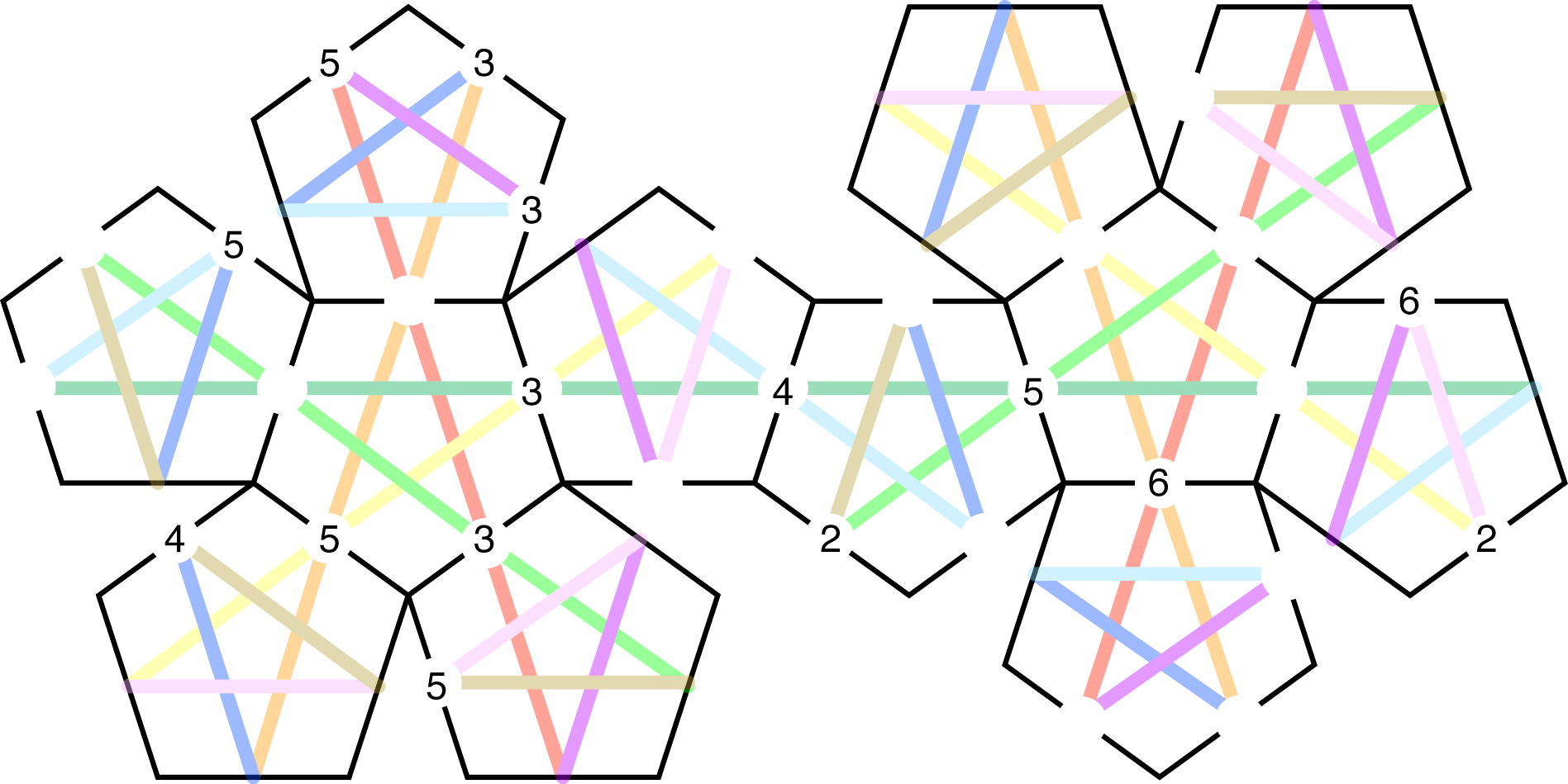}
\par\end{centering}

\protect\caption{\label{fig:criticalSetHardLatinDodecahedronEdge}Latin dodecahedron
puzzle}
\end{figure}

Fig. \ref{fig:hardLatinDodecahedronEdges} shows a \emph{Latin dodecahedron}\footnote{Or Latin icosahedron if we follow the first naming criterion in Sect.
\ref{sub:Latin-Polyhedra}.} resulting from the latinization of the board in Fig. \ref{fig:latinBaseDodecahedronEdgesAndLines}
with all numbers from $1$ to $6$. The derived woven board has SIN
$\left\{ 0,1,2\right\} $ because there is an orthogonal warp class
with six asterisms with five points each. Fig. \ref{fig:criticalSetHardLatinDodecahedronEdge}
shows a corresponding Latin puzzle.

\subsection{Equivalence Classes of Latin Polytopes\label{sub:Equivalence-Classes-of-Latin-Polytopes}}

We saw in Sect. \ref{sub:4.1-Equivalence-classes-Latin boards} that
a set of operations that easily generate Latin boards from a given
one serves to define a notion of equivalence among boards that classify
them into equivalence classes,. We also saw how this was useful when
enumerating or counting Latin boards.

On the other hand, asterisms in a Latin polytope are taken to other
asterisms by the action of a symmetry group. This not only gives us
some automorphisms of the board \textendash see Relation \ref{Relation three groups }\textendash{}
it also moves the labels so as to produce another, not necessarily
different, Latin polytope\footnote{Not all board automorphisms produce Latin polytopes though. In the
case of the board of a Latin square for example, swapping just the
first row and first column of points is a board automorphism that
does not produce a Latin square in general. }. This means that, in absence of larger groups, geometric symmetry
may be used to derive equivalence classes.

As commented in Sect. \ref{sub:4.1-Equivalence-classes-Latin boards},
the larger the symmetry group the less equivalence classes there will
be and the more boards per class will result\footnote{The precise number of classes generated by the group can be calculated
using \emph{Burnside's Lemma} \cite{Taking Sudoku Seriously}.}. This is illustrated by the next example. We know that the board
of a Latin square has \textit{D}\textsubscript{4} symmetry. Some
elements of this group have the same effect as an isotopism \textendash for
example a reflection across the vertical symmetry axis\textendash{}
while others produce a conjugate Latin square \textendash for example
a reflection across the top-left to bottom-right diagonal, that swaps
rows and columns.

\begin{figure}
\begin{centering}
\includegraphics[height=3cm]{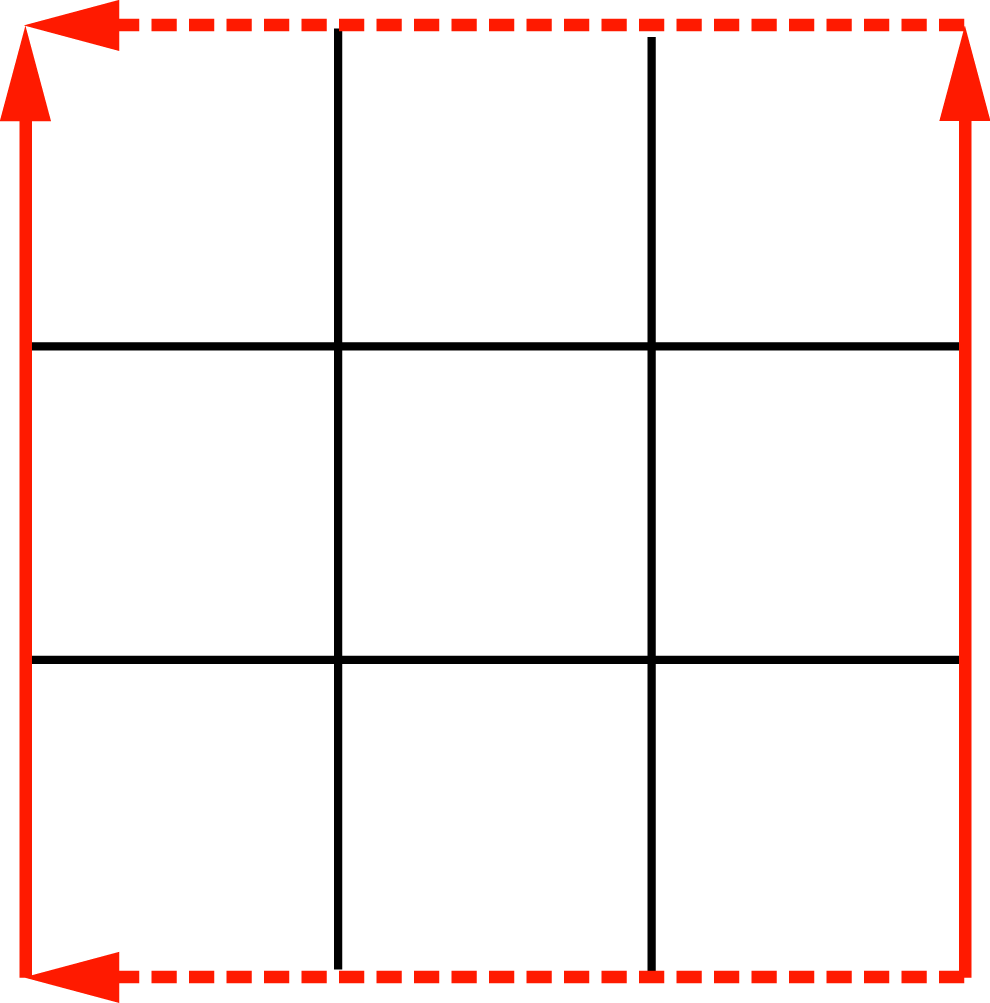}
\par\end{centering}

\protect\caption{\label{fig:torus board}Square-tiled toroidal source}

\end{figure}

The symmetry group of the square-tiled torus of order $n$, whose
fundamental polygon features single horizontal and vertical borders
(see Fig. \ref{fig:torus board} for the $n=3$ case), has \textit{D}\textsubscript{4}
as a subset. In a square-tiled torus of order $n$ the group also
includes rotations of $\frac{360\text{º}}{n}\times k,k=1,\ldots,n$,
around both its vertical and internal axis. These rotations induce
circular permutations in rows and columns; these permutations are
a subset of the full set of permutations of rows and columns in the
Latin square. The symmetry group of the square-tiled torus of order
$n$ will produce then less equivalence classes of Latin squares of
order $n$ than $D_{4}$, but still more than the number of isotopic
classes (see Sect. \ref{sub:Equivalence-classes-of-Latin squares}).

\section{Conclusions}

Taking Latin squares as a departure point, we have introduced the
concept of latinization and defined boards, Latin boards, symmetric
boards and Latin polytopes. We have also provided examples that prove
their existence for given values of the relevant parameters. 

These results are accompanied by corresponding generalizations of
other concepts related to Latin squares, like the orthogonal intersection
property in sets of mutually orthogonal Latin squares (MOLS). The
generalization of uniquely completable partial Latin squares to Latin
puzzles has also been mentioned. This topic is studied in the companion
paper \emph{Latin Puzzles} \cite{Latin Puzzles paper}, where several
latinization techniques and their computational complexity are also
discussed.

The new defined combinatorial  objects offer varied compromises between
balance and symmetry. This is especially true in the case of symmetric
boards and Latin polytopes, in which geometric symmetry translates
naturally into both combinatorial symmetry (resolvable boards with
simple intersection patterns) and algebraic symmetry (automorphisms
of the board and equivalence classes of Latin polytopes). On a related
note, we have seen how the use of multisets \textendash rather than
sets\textendash{} for latinization offers the possibility to fine-tuning
these properties in the resulting Latin polytopes and woven boards.

We have proposed several ways to construct symmetric boards from sources
of symmetry (orthogonal projections of sources onto lower dimensions,
kaleidoscopic sets generated by the action of symmetry groups, etc.).
We have also discussed ways to enhance the symmetry of a given source
(composition, truncation, tessellation, etc.), with an emphasis on
tessellations based on biregular polytopes, a method that can easily
generate whole families of sources of symmetry, and hence whole families
of symmetric boards and Latin polytopes. 

Finally, we have also touched upon the problem of deriving symmetric
boards from existing non-symmetric ones.

\section{Suggestions for Future Work}

Future work on Latin polytopes could proceed along the lines mentioned
below.

\subsubsection*{\noindent New Latin Boards and Polytopes}
\begin{itemize}
\item based on the different sources of symmetry listed in Sect. \ref{sub:Sources of Symmetry}
\item in spaces other than $\mathbb{R}^{n}$ (elliptic, hyperbolic or complex
spaces \cite{Regular Complex Polytopes,Regular polytopes})
\item Latin abstract polytopes (from sources derived from abstract regular
polytopes \cite{Abstract Regular Polytopes})
\item class-symmetric polyhedral boards ( in Sect. \ref{exa:Boards-from-Polyhedra with bireg faces}
only the tetrahedral board was class-symmetric)
\end{itemize}

\subsubsection*{Study of Properties}
\begin{itemize}
\item classification and enumeration 
\item notions of isotopism and conjugacy when applicable
\item relation with hypergraphs (automorphism group, coloring, etc.)
\item relation with combinatorial designs (extension of the notion of orthogonality,
resolvability, etc.)
\item notions of equivalence and derived equivalence classes
\item properties of related partial objects like partial Latin boards \cite{Latin Puzzles paper}
\end{itemize}

\subsubsection*{Applications}
\begin{itemize}
\item design of experiments
\item coding
\item scheduling
\item timetabling
\item resource allocation
\item mathematical puzzles
\end{itemize}

\section*{Acknowledgements\addcontentsline{toc}{section}{Acknowledgements}}

The author thanks Roberto Álvarez Chust, Prof. Francisco Ballesteros
Olmo, Latin squares expert Prof. Raúl Falcón Ganfornina and Manuel
Valderrama Conde for their thoughtful comments upon reading this paper.
The author also thanks the members of the \textit{Tertulia de Matemáticas}
\cite{Tertulia de Matem=0000E1ticas} for all the great conversations.

\paragraph*{\noindent About the author}

\noindent \textbf{\addcontentsline{toc}{section}{About the author}}{\small{}The
author is a telecommunication engineer from the Polytechnical University
of Madrid.}{\small \par}

\pagebreak{}

\tableofcontents{}\label{TOC}\addcontentsline{toc}{section}{Contents}

\pagebreak{}

\listoffigures

\addcontentsline{toc}{section}{List of Figures}

\bigskip{}

\noindent 

\begin{thebibliography}{10}
\bibitem{The history of Latin squares Andersen} \textbf{\addcontentsline{toc}{section}{References}}Andersen
L.D.: Chapter on \href{http://vbn.aau.dk/files/13649565/R-2007-32.pdf}{\textquotedblleft The history of Latin squares\textquotedblright},
intended for publication in \emph{The history of combinatorics}. Preliminary
version, Department of Mathematical Sciences, Aalborg University,
Research Report Series; No. R-2007-32, 2008. Accessed on 5 January
2016.

\bibitem{Combinatorics of Finite Sets} Anderson, I.: \emph{Combinatorics
of Finite Sets}, Clarendon Press, Oxford, 1987.

\bibitem{Principles of Constraint Programming Apt} Apt K.R.: \emph{Principles
of Constraint Programming}, Cambridge University Press, Cambridge,
2003.

\bibitem{full automorphism groups of edge-transitive hypergraphs}
Babai, L., Cameron P.J.: Most primitive groups are full automorphism
groups of edge-transitive hypergraphs, \emph{Journal of algebra} 421:
512-523, 2015, doi: 10.1016/j.jalgebra.2014.09.002.

\bibitem{The Fano Plane} Baez J.: The Octonions, \emph{Bulletin of
the American Mathematical Society} 39 (2): 145\textendash 205, 2002.
DOI:10.1090/S0273-0979-01-00934-X.

\bibitem{hypergraphs} Berge, C. \textit{Hypergraphs}, North-Holland
Mathematical Library, Amsterdam, 1989.

\bibitem{Space Groups for Scientists and Engineers} Burns G., Glazer
A. M.: \emph{Space Groups for Scientists and Engineers}, 2nd edition,
Academic Press, Inc., Boston, 1990.

\bibitem{Handbook of Combinatorial Designs} Colbourn C., Dinitz J.,
Jeffrey H. (eds.):\textit{ Handbook of Combinatorial Designs}, 2nd
edition, CRC Press Inc., Boca Raton, Florida, 2007.

\bibitem{Regular Complex Polytopes} Coxeter, H. S. M.: \emph{Regular
Complex Polytopes,} 2nd edition, Cambridge University Press, Cambridge,
1991.

\bibitem{Regular polytopes}  Coxeter, H. S. M.: \emph{Regular Polytopes},
3rd edition, Dover Publications Inc., New York, 1973.

\bibitem{Area requirement and symmetry display of planar upward drawings}
Di Battista G., Tamassia R., Tollis, I.G.: Area requirement and symmetry
display of planar upward drawings, \emph{Discrete and Computational
Geometry} 7 (1): 381\textendash 401, 1992, doi:10.1007/BF02187850.

\bibitem{De quadratis magicis} Euler L.: De quadratis magicis, \emph{Commentationes
arithmeticae} 2: 593-602, 1849.

\bibitem{On chromatic number of graphs and set systems} Erdös, P.,
Hajnal, A.: On chromatic number of graphs and set systems. \emph{Acta
Mathematica, }Academy of Sciences, Hungary 17: 61-99, 1966.

\bibitem{Enumerating possible Sudoku grids} Felgenhauer B., Jarvis
F.: Enumerating possible Sudoku grids. Available from \href{http://www.afjarvis.staff.shef.ac.uk/sudoku/sudoku.pdf}{http://www.afjarvis.staff.shef.ac.uk/sudoku/sudoku.pdf}.

\bibitem{Latin triangles} Halberstam F.Y., Hoffman D.G., Richter
R.B.: Latin triangles. \emph{Ars Combinatoria} 21: 51-58, 1986. MR
87h:05047.

\bibitem{Drawing graphs symmetrically in three dimensions} Hong,
S.H.: Drawing graphs symmetrically in three dimensions, \emph{Proc.
9th Int. Symp. Graph Drawing} (GD 2001): 106\textendash 108, 2002,
doi:10.1007/3-540-45848-4\_16.

\bibitem{Latin squares and their applications} Keedwell A.D., Dénes
J.: \textit{Latin squares and their Applications}, 2nd edition, North
Holland, Amsterdam, 2015.

\bibitem{Some NP-complete problems similar to graph isomorphism}
Lubiw, A.: Some NP-complete problems similar to graph isomorphism,
\emph{SIAM Journal on Computing} 10 (1): 11\textendash 21, 1981, doi:10.1137/0210002.

\bibitem{Abstract Regular Polytopes} McMullen P., Schulte E.: \emph{Abstract
Regular Polytopes}, Cambridge University Press, 2002.

\bibitem{Number of Paratopy classes of Latin Squares} OEIS. \textit{Number
of Paratopy classes of Latin Squares of order n}.\\
Available from \href{http://oeis.org/A003090}{http://oeis.org/A003090}.

\bibitem{Latin Puzzles paper} Palomo M.G.: Latin Puzzles. Avail.
from \href{http://arxiv.org/abs/1602.06946}{http://arxiv.org/abs/1602.06946}.

\bibitem{Latin Puzzles site} Palomo, M.G.: Latin Puzzles site. \href{http://www.latinpuzzles.com}{http://www.latinpuzzles.com}.
Accessed on 25 January 2016.

\bibitem{Taking Sudoku Seriously} Rosenhouse J., Taalman L.: \textit{Taking
Sudoku Seriously}, Oxford University Press, New York, 2011.

\bibitem{Stinson Combinatorial Designs} Stinson D.R.: \textit{Combinatorial
Designs, Constructions and Analysis}, Springer-Verlag, New York, 2004.

\bibitem{An illustrated philatelic introduction} Styan G.P.H.: \href{http://www.math.mcgill.ca/styan/Smolenice-beamer-16june09}{An illustrated philatelic introduction to 4\texttimes 4 Latin squares in Europe}:
1283\textendash 1788 (with some comments about le fauteuil Dagobert
and about philatelic Latin squares in Europe: 1984\textendash 2007),
\emph{Report 2009-02, }Department of Mathematics and Statistics, McGill
University, Montréal, 2009. Accessed on 5 January 2016.

\bibitem{Some comments on Latin squares and on Graeco-Latin squares}
Styan G.P.H., Boyer C., Chu K.L.: Some comments on Latin squares and
on Graeco-Latin squares, illustrated with postage stamps and old playing
cards, \emph{Statistical Papers} 50: 917\textendash 941, 2009. DOI
10.1007/00362-009-0261-5.2009.

\bibitem{Tertulia de Matem=0000E1ticas} Tertulia de Matemáticas.
\href{http://www.tertuliadematematicas.miguelpalomo.com/}{www.tertuliadematematicas.miguelpalomo.com}.
Accessed on 25 January 2016.

\bibitem{Coloring Mixed Hypergraphs} Voloshin V.I.: \emph{Coloring
Mixed Hypergraphs: Theory, Algorithms and Applications}, Fields Institute
Monographs, American Mathematical Society, 2002.\pagebreak{}

\end{thebibliography}
\end{document}